\documentclass[12pt,a4paper]{amsart}

\usepackage{amsmath,amsfonts,amssymb,amsthm, tikz-cd, mathtools,stackengine}
\usepackage{graphicx} 
\usepackage{xy}
\xyoption{all}
\usepackage[colorlinks, linkcolor=blue!72,anchorcolor=orange,
    citecolor=red,urlcolor=Emerald, bookmarksopen,bookmarksdepth=2]{hyperref}
\usetikzlibrary{matrix,positioning,decorations.markings,arrows,decorations.pathmorphing,backgrounds,fit,positioning,shapes.symbols,chains,shadings,fadings,calc, shadows,shapes.geometric}
\usepackage[dvipsnames]{xcolor}
\usepackage{tikz}
\usepackage{cleveref}
\usepackage{url}
\usepackage{stmaryrd}
\usepackage{lipsum}
\usepackage{tcolorbox}

\tcbuselibrary{breakable}
\tcbset{
  width=0.7\textwidth,
  halign=justify,
  center,
  breakable,
  colback=white    
}

\usetikzlibrary{matrix,positioning,decorations.markings,arrows.meta,decorations.pathmorphing,backgrounds,fit,positioning,shapes.symbols,chains,shadings,fadings,calc}
\tikzset{->-/.style={decoration={  markings,  mark=at position #1 with
    {\arrow{>}}},postaction={decorate}}}
\tikzset{-<-/.style={decoration={  markings,  mark=at position #1 with
    {\arrow{<}}},postaction={decorate}}}
\theoremstyle{plain}
\newtheorem{theorem}{Theorem}[section]
\newtheorem{lemma}[theorem]{Lemma}
\newtheorem{proposition}[theorem]{Proposition}
\newtheorem{corollary}[theorem]{Corollary}
\theoremstyle{definition}
\newtheorem{definition}[theorem]{Definition}

\newtheorem{remark}[theorem]{Remark}
\newtheorem{example}[theorem]{Example}

\newtheorem*{theorem*}{Theorem}

\newtheorem*{definition*}{Definition}
\newtheorem*{corollary*}{Corollary}
\numberwithin{equation}{section}

\newcommand{\Hom}{\operatorname{Hom}}
\newcommand{\Eb}{\mathbb{E}}
\newcommand{\End}{\operatorname{End}}
\newcommand{\Ext}{\operatorname{Ext}}
\newcommand{\thick}{\operatorname{thick}}

\newcommand{\add}{\operatorname{add}}
\newcommand{\proj}{\operatorname{proj}}
\newcommand{\inj}{\operatorname{inj}}
\newcommand{\Filt}{\operatorname{Filt}}
\renewcommand{\ker}{\operatorname{ker}}
\newcommand{\coker}{\operatorname{coker}}
\newcommand{\im}{\operatorname{im}}

\newcommand{\ssim}{\operatorname{sim}}
\renewcommand{\mod}{\operatorname{mod}}
\newcommand{\Fac}{\operatorname{Fac}}
\newcommand{\rad}{\operatorname{rad}}
\newcommand{\Sub}{\operatorname{Sub}}

\newcommand{\emod}[2]{#1\text{-}\mod #2}

\newcommand{\AR}[1]{\tau_{[#1]}}
\newcommand{\id}{\operatorname{id}}

\def\K{\mathrm{K}}

\def\B{\mathcal{B}}
\def\C{\mathcal{C}}
\def\D{\mathcal{D}}
\def\E{\mathcal{E}}

\def\T{\mathcal{T}}
\def\F{\mathcal{F}}
\def\H{\mathcal{H}}

\def\SS{\mathcal{S}}

\def\W{\mathcal{W}}
\def\X{\mathcal{X}}
\def\Y{\mathcal{Y}}
\def\Z{\mathcal{Z}}
\def\La{\Lambda}
\def\LW{W_L}
\def\LWlen{W_L^{len}}

\def\Ho{\mathrm{H}}

\usepackage[maxbibnames=99]{biblatex}
\newcommand{\po}{\ar@{}[dr]|{\text{\pigpenfont R}}}
\newcommand{\pb}{\ar@{}[dr]|{\text{\pigpenfont J}}}

\newcommand\Rightarrowtail[2][]{\ensurestackMath{\mathrel{%
  \stackengine{1pt}{%
    \stackengine{0pt}{\rightarrowtail}{\scriptstyle#2}{O}{c}{F}{F}{S}%
  }{\scriptstyle#1}{U}{c}{F}{F}{S}%
}}}
\newcommand\Twoheadrightarrow[2][]{\ensurestackMath{\mathrel{%
  \stackengine{1pt}{%
    \stackengine{0pt}{\twoheadrightarrow}{\scriptstyle#2}{O}{c}{F}{F}{S}%
  }{\scriptstyle#1}{U}{c}{F}{F}{S}%
}}}

\newsavebox{\myboxa}
\sbox{\myboxa}{%
   \begin{tikzpicture}
   \begin{scope}[on grid,node distance=2cm]
    \node (A) {$P_2$};
    \node (B)[right = 2cm of A] {$I_1$};
    \node (C)[above right = 0.9cm and 1cm of A] {$P_1$};
    \node (D)[right = 2cm of C] {$P_2[1]$};
    \node (E)[right = 2cm of B] {$P_1[1]$};
    \node (F)[right = 2cm of D] {$I_1[1]$};
    \path   (A) edge[->]    (C)
        (C) edge[->]   (B)
        (B) edge[->]  (D);
      \path  (D) edge[->]  (E)
        (E) edge[->]  (F);
    \end{scope}
   \end{tikzpicture}
 }
 \newsavebox{\myboxb}
\sbox{\myboxb}{%
   \begin{tikzpicture}
   \begin{scope}[on grid,node distance=2cm]
    \node (A) {$P_2$};
    \node (B)[right = 2cm of A] {$I_1$};
    \node (C)[above right = 0.9cm and 1cm of A] {$P_1$};
    \node (D)[right = 2cm of C] {$P_2[1]$};
    \node (E)[right = 2cm of B] {$P_1[1]$};
    \node (F)[right = 2cm of D] {$I_1[1]$};
    \path   (A) edge[->]    (C)
        (C) edge[->]   (B)
        (B) edge[->]  (D);
      \path  (D) edge[->]  (E)
        (E) edge[->]  (F);
    \end{scope}
    \draw[fill=none](D) circle (12 pt);
     \begin{scope}[on background layer]\draw[fill=gray!20, draw=none] (D) circle (15pt);
     \end{scope}
   \end{tikzpicture}
 }
\newsavebox{\myboxc}
\sbox{\myboxc}{%
   \begin{tikzpicture}
   \begin{scope}[on grid,node distance=2cm]
    \node (A) {$P_2$};
    \node (B)[right = 2cm of A] {$I_1$};
    \node (C)[above right = 0.9cm and 1cm of A] {$P_1$};
    \node (D)[right = 2cm of C] {$P_2[1]$};
    \node (E)[right = 2cm of B] {$P_1[1]$};
    \node (F)[right = 2cm of D] {$I_1[1]$};
    \path   (A) edge[->]    (C)
        (C) edge[->]   (B)
        (B) edge[->]  (D);
      \path  (D) edge[->]  (E)
        (E) edge[->]  (F);
    \end{scope}
    \draw[fill=none](F) circle (12 pt);
     \begin{scope}[on background layer]\draw[fill=gray!20, draw=none] (F) circle (15pt);
     \end{scope}
   \end{tikzpicture}
 }
\newsavebox{\myboxd}
\sbox{\myboxd}{%
   \begin{tikzpicture}
   \begin{scope}[on grid,node distance=2cm]
    \node (A) {$P_2$};
    \node (B)[right = 2cm of A] {$I_1$};
    \node (C)[above right = 0.9cm and 1cm of A] {$P_1$};
    \node (D)[right = 2cm of C] {$P_2[1]$};
    \node (E)[right = 2cm of B] {$P_1[1]$};
    \node (F)[right = 2cm of D] {$I_1[1]$};
    \path   (A) edge[->]    (C)
        (C) edge[->]   (B)
        (B) edge[->]  (D);
      \path  (D) edge[->]  (E)
        (E) edge[->]  (F);
    \end{scope}
    \draw[fill=none](A) circle (10 pt);
     \begin{scope}[on background layer]\draw[fill=gray!20, draw=none] (A) circle (15pt);
     \end{scope}
   \end{tikzpicture}
 }
\newsavebox{\myboxe}
\sbox{\myboxe}{%
   \begin{tikzpicture}
   \begin{scope}[on grid,node distance=2cm]
    \node (A) {$P_2$};
    \node (B)[right = 2cm of A] {$I_1$};
    \node (C)[above right = 0.9cm and 1cm of A] {$P_1$};
    \node (D)[right = 2cm of C] {$P_2[1]$};
    \node (E)[right = 2cm of B] {$P_1[1]$};
    \node (F)[right = 2cm of D] {$I_1[1]$};
    \path   (A) edge[->]    (C)
        (C) edge[->]   (B)
        (B) edge[->]  (D);
      \path  (D) edge[->]  (E)
        (E) edge[->]  (F);
    \end{scope}
     \draw[fill=none](E) circle (13 pt);
     \begin{scope}[on background layer]\draw[rounded corners,fill=gray!20, draw=none] (E) circle (15 pt);;
     \end{scope}
   \end{tikzpicture}
 }
\newsavebox{\myboxf}
\sbox{\myboxf}{%
   \begin{tikzpicture}
   \begin{scope}[on grid,node distance=2cm]
    \node (A) {$P_2$};
    \node (B)[right = 2cm of A] {$I_1$};
    \node (C)[above right = 0.9cm and 1cm of A] {$P_1$};
    \node (D)[right = 2cm of C] {$P_2[1]$};
    \node (E)[right = 2cm of B] {$P_1[1]$};
    \node (F)[right = 2cm of D] {$I_1[1]$};
    \path   (A) edge[->]    (C)
        (C) edge[->]   (B)
        (B) edge[->]  (D);
      \path  (D) edge[->]  (E)
        (E) edge[->]  (F);
    \end{scope}
    \draw[fill=none](B) circle (10 pt);
     \begin{scope}[on background layer]\draw[rounded corners,fill=gray!20, draw=none] (B) circle (15pt);
     \end{scope}
   \end{tikzpicture}
 }
\newsavebox{\myboxg}
\sbox{\myboxg}{%
   \begin{tikzpicture}
   \begin{scope}[on grid,node distance=2cm]
    \node (A) {$P_2$};
    \node (B)[right = 2cm of A] {$I_1$};
    \node (C)[above right = 0.9cm and 1cm of A] {$P_1$};
    \node (D)[right = 2cm of C] {$P_2[1]$};
    \node (E)[right = 2cm of B] {$P_1[1]$};
    \node (F)[right = 2cm of D] {$I_1[1]$};
    \path   (A) edge[->]    (C)
        (C) edge[->]   (B)
        (B) edge[->]  (D);
      \path  (D) edge[->]  (E)
        (E) edge[->]  (F);
    \end{scope}
    \draw[fill=none](D) circle (13 pt);
     \draw[fill=none](F) circle (13 pt);
     \begin{scope}[on background layer]\draw[rounded corners,fill=gray!20, draw=none] (2.7,1.2)--(4,-0.1)--(5.3,1.2)--cycle;
     \end{scope}
   \end{tikzpicture}
 }
\newsavebox{\myboxh}
\sbox{\myboxh}{%
   \begin{tikzpicture}
   \begin{scope}[on grid,node distance=2cm]
    \node (A) {$P_2$};
    \node (B)[right = 2cm of A] {$I_1$};
    \node (C)[above right = 0.9cm and 1cm of A] {$P_1$};
    \node (D)[right = 2cm of C] {$P_2[1]$};
    \node (E)[right = 2cm of B] {$P_1[1]$};
    \node (F)[right = 2cm of D] {$I_1[1]$};
    \path   (A) edge[->]    (C)
        (C) edge[->]   (B)
        (B) edge[->]  (D);
      \path  (D) edge[->]  (E)
        (E) edge[->]  (F);
    \end{scope}
    \draw[fill=none](A) circle (10 pt);
     \draw[fill=none](F) circle (13 pt);
     \begin{scope}[on background layer]\draw[rounded corners,fill=gray!20, draw=none](A) circle (15pt) (F) circle (15pt);
     \end{scope}
   \end{tikzpicture}
 }
\newsavebox{\myboxi}
\sbox{\myboxi}{%
   \begin{tikzpicture}
   \begin{scope}[on grid,node distance=2cm]
    \node (A) {$P_2$};
    \node (B)[right = 2cm of A] {$I_1$};
    \node (C)[above right = 0.9cm and 1cm of A] {$P_1$};
    \node (D)[right = 2cm of C] {$P_2[1]$};
    \node (E)[right = 2cm of B] {$P_1[1]$};
    \node (F)[right = 2cm of D] {$I_1[1]$};
    \path   (A) edge[->]    (C)
        (C) edge[->]   (B)
        (B) edge[->]  (D);
      \path  (D) edge[->]  (E)
        (E) edge[->]  (F);
    \end{scope}
    \draw[fill=none](B) circle (10 pt);
    \draw[fill=none](E) circle (13 pt);
     \begin{scope}[on background layer]\draw[rounded corners,fill=gray!20, draw=none] (1.8,-0.2)--(3,1.1)--(4.2,-0.2)--cycle;
     \end{scope}
   \end{tikzpicture}
 }
\newsavebox{\myboxj}
\sbox{\myboxj}{%
   \begin{tikzpicture}
   \begin{scope}[on grid,node distance=2cm]
    \node (A) {$P_2$};
    \node (B)[right = 2cm of A] {$I_1$};
    \node (C)[above right = 0.9cm and 1cm of A] {$P_1$};
    \node (D)[right = 2cm of C] {$P_2[1]$};
    \node (E)[right = 2cm of B] {$P_1[1]$};
    \node (F)[right = 2cm of D] {$I_1[1]$};
    \path   (A) edge[->]    (C)
        (C) edge[->]   (B)
        (B) edge[->]  (D);
      \path  (D) edge[->]  (E)
        (E) edge[->]  (F);
    \end{scope}
    \draw[fill=none](C) circle (10 pt);
     \begin{scope}[on background layer]\draw[rounded corners,fill=gray!20, draw=none] (C) circle (15 pt);
     \end{scope}
   \end{tikzpicture}
 }
\newsavebox{\myboxk}
\sbox{\myboxk}{%
   \begin{tikzpicture}
   \begin{scope}[on grid,node distance=2cm]
    \node (A) {$P_2$};
    \node (B)[right = 2cm of A] {$I_1$};
    \node (C)[above right = 0.9cm and 1cm of A] {$P_1$};
    \node (D)[right = 2cm of C] {$P_2[1]$};
    \node (E)[right = 2cm of B] {$P_1[1]$};
    \node (F)[right = 2cm of D] {$I_1[1]$};
    \path   (A) edge[->]    (C)
        (C) edge[->]   (B)
        (B) edge[->]  (D);
      \path  (D) edge[->]  (E)
        (E) edge[->]  (F);
    \end{scope}
    \draw[fill=none](C) circle (10 pt);
     \draw[fill=none](D) circle (13 pt);
     \begin{scope}[on background layer]\draw[rounded corners,fill=gray!20, draw=none] (1,1.1)--(2,-0.2)--(3.2,1.1)--cycle;
     \end{scope}
   \end{tikzpicture}
 }
\newsavebox{\myboxl}
\sbox{\myboxl}{%
   \begin{tikzpicture}
   \begin{scope}[on grid,node distance=2cm]
    \node (A) {$P_2$};
    \node (B)[right = 2cm of A] {$I_1$};
    \node (C)[above right = 0.9cm and 1cm of A] {$P_1$};
    \node (D)[right = 2cm of C] {$P_2[1]$};
    \node (E)[right = 2cm of B] {$P_1[1]$};
    \node (F)[right = 2cm of D] {$I_1[1]$};
    \path   (A) edge[->]    (C)
        (C) edge[->]   (B)
        (B) edge[->]  (D);
      \path  (D) edge[->]  (E)
        (E) edge[->]  (F);
    \end{scope}
    \draw[fill=none](A) circle (10 pt);
     \draw[fill=none](B) circle (10 pt);
     \begin{scope}[on background layer]\draw[rounded corners,fill=gray!20, draw=none] (-0.2,-0.2)--(1,1.1)--(2.2,-0.2)--cycle;
     \end{scope}
   \end{tikzpicture}
 }

\title{Semibricks and wide subcategories in extended module categories}

\dedicatory{Dedicated to Claus Michael Ringel on the occasion of his 80th birthday}

\author{Esha Gupta}
\address{Université Paris-Saclay, UVSQ, CNRS, Laboratoire de Mathématiques de Versailles, 78000, Versailles, France.}
\email{esha.gupta2@uvsq.fr}

\author{Yu Zhou}
\address{School of Mathematical Sciences,
Beijing Normal University,
100875 Beijing,
China}
\email{yuzhoumath@gmail.com}

\date{}



\thanks{The work was supported by National Natural Science Foundation of China (Grant No. 12671045).}

\begin{document}
\begin{abstract}
For $d \geq 1$, we introduce semibricks and wide subcategories in the $d$-extended heart of a bounded $t$-structure on a triangulated category. We establish a bijection between semibricks and length wide subcategories. Furthermore, we investigate three types of torsion-like classes, namely $d$-FAE closed subcategories, positive torsion classes, and $s$-torsion classes. We show that these notions, while coinciding in the classical case, play distinct roles in the higher setting. When specialized to the extended module category of a finite-dimensional algebra, we define left/right-finite semibricks and wide subcategories. We then establish bijections between these classes and $(d+1)$-term simple-minded collections as well as functorially finite positive torsion classes, which generalizes the known correspondences between $2$-term simple-minded collections and the respective finite classes in module categories. Finally, by exploiting a link between semibricks and silting complexes, we provide a criterion to characterize which mutations of $(d+1)$-term silting complexes remain $(d+1)$-term.
\end{abstract}
\maketitle
\tableofcontents{}

\section{Introduction}

Let $\La$ be a finite-dimensional algebra over a field $K$. In the representation theory of finite-dimensional algebras, simple modules play a fundamental role as they generate the category $\mod\La$ of finite-dimensional modules under extensions. The notion of simple modules can be generalised to the class of bricks, which are modules whose endomorphism algebra is a division algebra. These have been extensively studied in the literature. Generalising the fact that there are no non-zero morphisms between two non-isomorphic simple modules, one considers the class of semibricks, i.e., a collection of pairwise Hom-orthogonal bricks. A classical result of Ringel \cite[\S~1.2]{R} shows that semibricks are in bijection with wide subcategories of $\mod\La$, which are full subcategories that are closed under kernels, cokernels, and extensions.

In the context of $\tau$-tilting theory, Asai \cite{A} introduced a finiteness condition on semibricks by defining a semibrick to be left-finite if the smallest torsion class containing it is functorially finite. He established that these left-finite semibricks are in bijection with functorially finite torsion classes in $\mod\La$. By utilizing Ringel's bijection, he subsequently extended the notion of left-finiteness to wide subcategories, thereby recovering the bijection between left-finite wide subcategories and functorially finite torsion classes previously studied by Marks and Šťovíček \cite{MS}. This correspondence was first established for hereditary algebras by Ingalls and Thomas \cite{IT}. 

By viewing simple modules as objects in the bounded derived category $\D^b(\mod\La)$, one generalizes them to simple-minded collections. 
A simple-minded collection is called 2-term if the cohomology of all its elements is concentrated in degrees -1 and 0. Asai established that left-finite semibricks in the module category are also in bijection with these 2-term simple-minded collections. Furthermore, these collections were shown by Brüstle and Yang \cite{BY} to be in bijection with 2-term silting objects in the bounded homotopy category $\K^b(\proj\La)$ of projective modules.

Asai also established a direct bijection between left-finite semibricks and support $\tau$-tilting modules by considering the quotient of a support $\tau$-tilting module by its radical over its endomorphism algebra. Support $\tau$-tilting modules were introduced by Adachi, Iyama, and Reiten \cite{AIR} as a mutation completion of tilting modules. These modules were shown to be in bijection with functorially finite torsion classes in the module category and with 2-term silting objects. The connection between 2-term silting objects and support $\tau$-tilting modules can also be found in the work of Derksen and Fei \cite{DF}. In the work of Pauksztello and Zvonareva \cite{PZ}, functorially finite torsion classes were shown to be in bijection with complete cotorsion classes in the full subcategory $\K^{[-1,0]}(\proj\La)$ of $\K^b(\proj\La)$ consisting of complexes concentrated in degrees -1 and 0. A summary of these results is provided in Figure~\ref{1}. Furthermore, a bijection was established by Garcia \cite{Ga1} between left-finite wide subcategories and thick subcategories with enough injectives in $\K^{[-1,0]}(\proj\La)$. Various relations between torsion classes, cotorsion classes, 2-term silting objects, and intermediate $t$-structures have also been extensively studied \cite{BY, IJY, PZ}.

\begin{figure}[htpb]
\centering
\begin{tikzpicture}[>=stealth,thick,font=\sffamily, elli/.style args={#1 and #2}{ellipse,minimum width=#1cm, minimum height=#2cm,align=center}]
\useasboundingbox (-1.5,-4) rectangle (12.5,4.3);
\begin{scope}[local bounding box=ells,fill opacity=0.4,text opacity=1]
    \node[elli=1 and 1] (E1) {\parbox{30mm}{\centering left-finite \\ semibricks\\ in $\mod\La$}};
    \node[elli=1 and 1, right = 0.3cm of E1.east] (E2) {\parbox{35mm}{\centering functorially finite torsion classes\\ in $\mod\La$}};
    \node[elli=1 and 1, right = 0.3cm of E2.east] (E3) {\parbox{35mm}{\centering basic $2$-term \\ silting  objects \\ in $\K^{b}(\proj\La)$}};
    \node[elli=1 and 1, above = 0.7cm of E3.north west] (E4) {\parbox{45mm}{\centering complete \\ cotorsion classes \\ in $\K^{[-1,0]}(\proj\La)$}};
    \node[elli=1 and 1, below = 0.7cm of E3.south west] (E5) {\parbox{45mm}{\centering basic support \\ $\tau$-tilting modules \\ in $\mod\La$}};
    \node[elli=1 and 1, below = 0.7cm of E1.south east] (E6) {\parbox{45mm}{\centering left-finite \\ wide subcategories \\ in $\mod\La$}};
    \node[elli=1 and 1, above = 0.7cm of E1.north east] (E7) {\parbox{45mm}{\centering $2$-term \\ simple-minded collections \\ in $\D^b(\mod\La)$}};
\end{scope}
  \draw ([xshift=-1cm]E1.east)--([xshift=1em]E2.west) node [midway, below] (E) {\cite{A}};
  \draw ([xshift=-1em]E2.east)--([xshift=.7cm]E3.west) node [midway, below] (F) {\cite{AIR}};
  \draw ([xshift=-1em]E2.north east)--([xshift=3em]E4.south west) node [midway, below, sloped] (G) {\cite{PZ}};
  \draw ([xshift=-1em]E2.south east)--([xshift=3em]E5.north west) node [midway, below, sloped] (H) {\cite{AIR}};
  \draw ([xshift=1em]E2.south west)--([xshift=-3em]E6.north east) node [midway, above, sloped] (I) {\cite{IT,MS}};
  \draw ([xshift=1em]E2.north west)--([xshift=-3em]E7.south east) node [midway, below, sloped] (J) {\cite{BY}};
  \draw ([xshift=-1em,yshift=.8em]E3.south)--([xshift=1em,yshift=-.9em]E5.north) node [midway, above, sloped] (K) {\cite{AIR,DF}};
  \draw ([xshift=-1em,yshift=-.8em]E3.north)--([xshift=1em,yshift=.9em]E4.south) node [midway, below, sloped] (L) {\cite{AT}};
  \draw ([xshift=1em,yshift=-.8em]E1.north)--([xshift=-1em,yshift=.9em]E7.south) node [midway, below, sloped] (M) {\cite{A}};
  \draw ([xshift=1em,yshift=.8em]E1.south)--([xshift=-1em,yshift=-.9em]E6.north) node [midway, above, sloped] (N) {\cite{R,A}};
  \draw ([xshift=-1em]E7.north east)to[out=15,in=180]node[below]{\cite{BY}} ([yshift=1em]E4.north)to[out=0,in=105] ([xshift=-2em]E3.north east);
  \draw ([xshift=1em]E5.south west)to[out=-170,in=0]node[above]{\cite{A}} ([yshift=-1em]E6.south)to[out=180,in=-75] ([xshift=2em]E1.south west);
 \end{tikzpicture}
     \caption{Bijections for $d=1$.}
     \label{1}
 \end{figure}
 
Our goal in this paper is to generalize the previous results to a higher setting. Specifically, we aim to extend these classes to establish bijections with arbitrary $(d+1)$-term simple-minded collections for any $d \geq 1$. This is achieved by considering subcategories of the form
$$\emod{d}{\La} := \D^{[-d+1,0]}(\mod\La)$$
which consist of objects in the bounded derived category $\D^b(\mod\La)$ whose cohomology is concentrated in degrees from $-d+1$ to $0$. These subcategories are referred to as extended module categories or truncated derived categories. In a previous work \cite{Gu1}, the first author introduced a notion of torsion classes in these categories and showed that functorially finite, positive torsion classes are in bijection with $(d+1)$-term silting objects and complete, hereditary cotorsion classes. These extended module categories were also investigated by the second author \cite{Z}, who established an Auslander-Reiten theory on these extriangulated categories (see also \cite{MP} for a generalisation). This structure involves an Auslander-Reiten translation denoted by $\tau_{[d]}$. Furthermore, the second author introduced the notion of a $\tau_{[d]}$-tilting pair and demonstrated its bijection with $(d+1)$-term silting objects. Recently, Iyama and Jin proved that $(d+1)$-term silting complexes are also in bijection with $(d+1)$-term simple-minded collections \cite{IJ}.

In this paper, we introduce semibricks and wide subcategories in extended module categories, specifically focusing on those with the length property, and establish a bijection between them. These concepts are detailed in Definitions~\ref{semibrick}~and~\ref{wide}. We subsequently define left-finite semibricks and left-finite wide subcategories as those for which the smallest positive torsion class containing them is functorially finite. This framework allows us to establish the comprehensive set of bijections illustrated in Figure~\ref{2}.

\begin{figure}[htpb]
\centering
\begin{tikzpicture}[>=stealth,thick,font=\sffamily, elli/.style args={#1 and #2}{ellipse,minimum width=#1cm, minimum height=#2cm,align=center}]
\useasboundingbox (-1.5,-3.5) rectangle (12.5,4.3);
\begin{scope}[local bounding box=ells,fill opacity=0.4,text opacity=1]
    \node[elli=1 and 1] (E1) {\parbox{30mm}{\centering left-finite semibricks\\ in $\emod{d}{\La}$}};
    \node[elli=1 and 1, right = 0.3cm of E1.east] (E2) {\parbox{35mm}{\centering functorially finite,\\ positive torsion\\ classes in $\emod{d}{\La}$}};
    \node[elli=1 and 1, right = 0.3cm of E2.east] (E3) {\parbox{35mm}{\centering basic $(d+1)$-term silting objects \\ in $\K^{b}(\proj\La)$}};
    \node[elli=1 and 1, above = 0.7cm of E3.north west] (E4) {\parbox{45mm}{\centering complete, hereditary \\ cotorsion classes \\ in $\K^{[-d,0]}(\proj\La)$}};
    \node[elli=1 and 1, below = 0.7cm of E3.south west] (E5) {\parbox{45mm}{\centering basic \\ $\AR{d}$-tilting pairs \\ in $\emod{d}{\La}$}};
    \node[elli=1 and 1, below = 0.7cm of E1.south east] (E6) {\parbox{45mm}{\centering left-finite \\ wide subcategories \\ in $\emod{d}{\La}$}};
    \node[elli=1 and 1, above = 0.7cm of E1.north east] (E7) {\parbox{45mm}{\centering $(d+1)$-term \\ simple-minded collections \\ in $\D^b(\mod\La)$}};
\end{scope}
  \draw ([xshift=-1cm]E1.east)--([xshift=1em]E2.west) (E);
  \draw ([xshift=1em]E2.north west)--([xshift=-3em]E7.south east) (J);
  \draw ([xshift=1em,yshift=-.8em]E1.north)--([xshift=-1em,yshift=1em]E7.south) (M);
  \draw ([xshift=1em,yshift=-4em]E7)node{Theorem~\ref{main}};
  \draw ([xshift=1em]E2.south west)--([xshift=-3em]E6.north east) (I);
  \draw ([xshift=1em,yshift=.8em]E1.south)--([xshift=-1em,yshift=-1em]E6.north) (N);
  \draw ([xshift=1em,yshift=4em]E6)node{Theorem~\ref{main2}};
  \draw ([xshift=-1em]E2.east)--([xshift=.7cm]E3.west) (F);
  \draw ([xshift=-1em]E2.north east)--([xshift=3em]E4.south west) (G);
  \draw ([xshift=-1em,yshift=-.8em]E3.north)--([xshift=1em,yshift=1em]E4.south) (L);
  \draw ([xshift=-1em,yshift=-4em]E4)node{\cite{Gu1}};
  
  \draw ([xshift=-1em]E2.south east)--([xshift=3em]E5.north west) node [midway, above, sloped] (H) {\cite{Z}};
  \draw ([xshift=-1em,yshift=.8em]E3.south)--([xshift=1em,yshift=-1em]E5.north) node [midway, above, sloped] (K) {\cite{Z}};
  \draw ([xshift=-1em]E7.north east)to[out=15,in=180]node[below]{\cite{IJ}} ([yshift=1em]E4.north)to[out=0,in=105] ([xshift=-2em]E3.north east);
 \end{tikzpicture}
     \caption{Bijections for arbitrary $d$.}
     \label{2}
 \end{figure}
 
Although the relevant notions and results formally parallel those in the module category, the absence of tools of abelian categories in the extended case requires new methods for both the construction of the bijections and the underlying proof techniques.

First, in the category of finitely generated modules over a finite-dimensional algebra, torsion classes admit several equivalent characterizations \cite{D}. In the general $d$-setting, however, these characterizations are no longer equivalent (see Figure~\ref{fig:inclusions}). Such equivalences are used crucially in many constructions and proofs within the classical module category. For example, the smallest torsion class containing a given subcategory can be obtained by taking its closure under extensions and factors, a fact which fails for arbitrary $d$. This breakdown necessitates new methods and reveals that the resulting non-equivalent notions of torsion-like classes play distinct roles.


In the classical setting, a wide subcategory is defined as a subcategory closed under kernels, cokernels, and extensions. In the extended module category, this approach is no longer available as there is no suitable notion of kernels and cokernels. Instead, we construct the wide subcategory associated to a $d$-FAE closed subcategory via the exact heart construction. Notably, wide subcategories in extended module categories remain abelian but are not necessarily of finite length, representing a key point of divergence from the classical case.

Finally, in contrast to the classical case, a $(d+1)$-term simple-minded collection does not necessarily decompose into a disjoint union of its objects lying in the torsion and torsion-free parts. This phenomenon is closely related to the fact that the corresponding $(d+1)$-term silting complex may admit both left and right mutations at a given direct summand simultaneously. These features are detailed in Section~\ref{subsec:sb} and Section~\ref{S8}. As a result, the proofs in the general $d$-setting become substantially more involved than their counterparts in the module category.

We note that one of the primary components of our proofs is the set of bijections between silting objects, simple-minded collections, bounded $t$-structures with length heart, and bounded co-$t$-structures. These correspondences were established by K\"{o}nig and Yang \cite{KY}, extending results first proven by Keller and Nicolás \cite{KN} for homologically smooth non-positive differential graded algebras. Furthermore, the correspondence between silting objects and co-$t$-structures was also investigated by Mendoza, Sáenz, Santiago, and Souto \cite{MSSS}.

The paper is organised as follows. In Section~\ref{S2}, we fix a triangulated category equipped with a bounded $t$-structure and introduce its $d$-extended heart. We also present the concepts of silting objects and simple-minded collections along with their basic properties. In Section~\ref{S3}, we collect important properties of three types of torsion-like classes, namely, $d$-FAE closed subcategories, positive torsion classes, and $s$-torsion classes. In particular, we establish a bijection between $d$-FAE closed subcategories and bounded suspended subcategories. Section~\ref{S4} is devoted to the study of semibricks and wide subcategories in the $d$-extended heart, where we prove that semibricks are in bijection with length wide subcategories. In Section~\ref{S5}, we provide a systematic way to associate a wide subcategory with any $d$-FAE closed subcategory by introducing the notion of the exact heart. Section~\ref{S6} focuses on the properties of two specific maps, denoted by $T$ and $\phi$, which send a semibrick to its smallest positive torsion class and its smallest $d$-FAE closed subcategory, respectively. We also show that the original semibrick can be recovered from the exact hearts of these subcategories. In Section~\ref{S7}, we specialize to the case of extended module categories of finite-dimensional algebras. We define left-finite and right-finite semibricks as well as wide subcategories, and establish their bijections with functorially finite positive torsion classes and $(d+1)$-term simple-minded collections. We also demonstrate that these wide subcategories can be realized as module categories of certain finite-dimensional algebras. Finally, in Section~\ref{S8}, we provide a criterion to determine which mutations of a $(d+1)$-term silting complex remain $(d+1)$-term.

\vspace{1em}
\noindent \textbf{Acknowledgements.}
The first author would like to thank the second author for his hospitality during the visit to Beijing Normal University. She would also like to thank her PhD supervisor Pierre-Guy Plamondon for his support and helpful discussions. The second author would like to thank Xiao-Wu Chen for many insightful discussions.

The authors used ChatGPT Pro to find the example of a $d$-FAE closed subcategory that is not a positive torsion class (Example \ref{d-FAE-not-positive}). All the remaining results and examples are completely attributed to the authors.

\section{Notation and background}\label{S2}

Throughout this work, $\La$ denotes a finite-dimensional algebra over a field $K$. All subcategories are assumed to be full, closed under isomorphisms, and contain the zero object. For morphisms $f: X \to Y$ and $g: Y \to Z$, we denote by $g\circ f$ the composition $X \xrightarrow{f} Y \xrightarrow{g} Z$. For any category $\D$, we write $X \in \D$ to indicate that $X$ is an object of $\D$, and $\C\subseteq\D$ to indicate that $\C$ is a subcategory of $\D$. For any $X\in\D$, we use $\add X$ to denote the additive hull of $X$ in $\D$, i.e., the subcategory of $\D$ consisting of direct summands of finite direct sums of copies of $X$. For two subcategories $\X$ and $\Y$ of $\D$, the notation $\Hom(\X,\Y)=0$ means that $\Hom(X,Y)=0$ for all $X\in \X$ and $Y\in \Y$. For any $\C\subseteq\D$, a morphism $f:X\to C$ in $\D$ is called a \emph{left $\C$-approximation} of $X$ if $C\in\C$ and for any morphism $h:X\to C'$ with $C'\in\C$, there exists a morphism $g:C\to C'$ such that $h=g\circ f$. Such a left $\C$-approximation $f$ is \emph{minimal} if any endomorphism $g:C\to C$ satisfying $g\circ f=f$ is an isomorphism. A subcategory $\C\subseteq\D$ is called \emph{covariantly finite} if every object $X\in\D$ admits a left $\C$-approximation. The notions of \emph{(minimal) right $\C$-approximations} and \emph{contravariantly finite} subcategories are defined dually. A subcategory $\C\subseteq\D$ is called \emph{functorially finite} if it is both contravariantly finite and covariantly finite.

The shift functor in a triangulated category $\D$ is denoted by $[1]$. In a triangle
\[X\to Y\to Z\to X[1],\]
we often omit the last term and the final arrow when they are irrelevant to the discussion. For a morphism $f$ in $\D$, $C(f)$ denotes the cone of $f$. For any two subcategories $\X$ and $\Y$ of $\D$, we denote by $\X*\Y$ the subcategory of $\D$ consisting of objects $Z$ such that there exists a triangle 
\[X \to Z \to Y,\]
with $X\in \X$ and $Y\in \Y$. Note that the operation $*$ is associative due to the octahedral axiom. A subcategory $\X$ of $\D$ is called \emph{closed under extensions} if $\X\ast\X\subseteq\X$. For objects $X,Y\in\D$, we simplify $\add X*\add Y$ to $X*Y$. For any subcategory $\X$ of $\D$, $\thick\X$ denotes the smallest triangulated subcategory containing $\X$ and closed under direct summands.

Recall that an exact category is a pair $(\C,\E)$ with $\C$ an additive category and $\E$ a collection of kernel-cokernel pairs in $\C$ satisfying certain axioms. See \cite{Bu} for a detailed exposition on exact categories. Elements of $\E$ are called \emph{conflations}, denoted by $A\rightarrowtail B \twoheadrightarrow C$. The first map of a conflation is called an \emph{inflation} and the second a \emph{deflation}. A morphism $f\in \C$ is called \emph{admissible} if $f=g\circ h$ with $h$ a deflation and $g$ an inflation. For a subcategory $\X$ of $\C$, $\Filt(\X)$ denotes the subcategory of objects $X\in \C$ such that there exists a sequence of inflations
\[0=X_0\rightarrowtail X_1\rightarrowtail\cdots\rightarrowtail X_n=X,\] 
with the cokernel of $X_{i-1} \rightarrowtail X_i$ belonging to $\X$ for each $i$. Note that the category $\Filt(\X)$ is precisely the extension closure of $\X$, the smallest subcategory of $\C$ containing $\X$ that is closed under extensions. A non-zero object $X\in \C$ is called \emph{simple} if there is no conflation $L\rightarrowtail X \twoheadrightarrow M$ with $L,M\neq 0$. Note that this coincides with the classical definition for $\C$ being abelian. We denote the set of isoclasses of simple objects in $\C$ by $\ssim \C$. An exact category $\C$ is called a \emph{length} category if $\C=\Filt(\ssim \C)$.

\subsection{\texorpdfstring{$t$}{t}-structures and extended hearts}
Let $\D$ be a triangulated category. A pair of full subcategories $(\C^{\leq 0}, \C^{\geq 0})$ of $\D$ is called a \emph{$t$-structure} \cite{BBDG} if the following conditions hold:
\begin{enumerate}
    \item $\C^{\leq 0}[1]\subseteq\C^{\leq 0}$ and $\C^{\geq 0}[-1]\subseteq\C^{\geq 0}$;
    \item $\Hom(\C^{\leq 0},\C^{\geq 0}[-1])=0$; and
    \item $\D=\C^{\leq 0}*\C^{\geq 0}[-1]$.
\end{enumerate}
The full subcategory $\C^{\leq 0}\cap\C^{\geq 0}$ is called the \emph{heart} of $(\C^{\leq 0}, \C^{\geq 0})$, which is known to be an abelian category. A $t$-structure $(\C^{\leq 0}, \C^{\geq 0})$ is called \emph{bounded} if for every object $X\in\D$, there exists an integer $n>0$ such that $X[n]\in\C^{\leq 0}$ and $X[-n]\in\C^{\geq 0}$. A heart is called \emph{bounded} if it arises from a bounded $t$-structure. We refer to \cite[Lemma~3.2]{Bri} for an equivalent characterization of bounded hearts without using the notion of $t$-structures.

Throughout this paper, we assume that $\D$ is a triangulated category equipped with a bounded $t$-structure $(\D^{\leq 0}, \D^{\geq 0})$. This assumption implies, in particular, that $\D$ is idempotent complete \cite{LC}. Let $\H=\D^{\leq 0}\cap \D^{\geq 0}$ be the corresponding heart. For any integers $n$ and $m$, we set
\[\D^{\leq n}:=\D^{\leq 0}[-n],\ \D^{\geq m}:=\D^{\geq 0}[-m],\text{ and }\D^{[m,n]}:=\D^{\leq n}\cap\D^{\geq m}.\]
We denote by
\[\sigma_{\leq n}\colon\D\to \D^{\leq n}\text{ and }\sigma_{\geq m}\colon\D\to \D^{\geq m}\]
the truncation functors which are right and left adjoint to the inclusions $\D^{\leq n}\hookrightarrow\D$ and $\D^{\geq m}\hookrightarrow\D$, respectively. For any $X\in \D$ and $n\in\mathbb{Z}$, there exists a truncation triangle 
\[\sigma_{\leq n}(X)\to X\to \sigma_{\geq n+1}(X).\]

Fix $d\geq 1$. We refer to $\D^{[-d+1,0]}$ as the \emph{$d$-extended heart} of the $t$-structure. Note that $\D^{[0,0]}=\H$, and, more generally,
\[\D^{[-d+1,0]}=\H[d-1]*\H[d-2]*\cdots*\H,\]
which is also denoted as $d$-$\H$. Since $\D^{[-d+1,0]}$ is closed under extensions, it inherits the extriangulated structure in the sense of \cite{NP} from the triangulated structure of $\D$, where the extension bifunctor is given by 
\[\Eb(X,Y) := \Hom(X,Y[1]) \quad \text{for } X,Y \in \D^{[-d+1,0]}.\]
Moreover, it is an extriangulated category with negative first extensions \cite{AET}, defined by $\Eb^{-1}(X,Y):=\Hom(X,Y[-1])$.

For a finite-dimensional $K$-algebra $\La$, we denote by $\mod\La$ the category of finitely generated right $\La$-modules. We denote
\begin{itemize}
    \item the bounded derived category of $\mod\La$ by $\D^b(\mod\La)$,
    \item the category of finitely generated projective (resp. injective) right $\La$-modules by $\proj\La$ (resp. $\inj\La$), 
    \item the homotopy category of bounded complexes of finitely generated projective (resp. injective) modules over $\La$ by $\K^b(\proj\La)$ (resp. $\K^b(\inj\La)$).
\end{itemize}
The category $\D^b(\mod\La)$ is equipped with the standard $t$-structure defined by the pair $(\D^{\leq 0}(\mod\La),\D^{\geq 0}(\mod\La))$, where
\begin{equation*}
    \begin{split}
        \D^{\leq n}(\mod\La)&:=\{X\in \D^b(\mod\La)\mid \Ho^i(X) = 0 \text{ for all } i > n\},\\
        \D^{\geq m}(\mod\La)&:=\{X\in \D^b(\mod\La)\mid \Ho^i(X) = 0 \text{ for all } i < m\},\\
        \D^{[m,n]}(\mod\La)&:=\D^{\leq n}(\mod\La)\cap\D^{\geq m}(\mod\La).
    \end{split}
\end{equation*}
For $d\geq 1$, we refer to $\D^{[-d+1,0]}(\mod\La)$ as the \emph{$d$-extended module category} and denote it as $\emod{d}{\La}$.

\subsection{Silting objects}\label{sec:sil}

We briefly recall the definition of silting objects in $\K^b(\proj\La)$.

\begin{definition}
    Let $P\in \K^b(\proj\La)$. 
    \begin{enumerate}
        \item $P$ is called a \emph{silting object} if 
        \begin{enumerate}
        \item $\Hom(P,P[i])=0$ for all $i>0$;
        \item $\thick P=\K^b(\proj\La)$.
        \end{enumerate}
        \item A silting object is said to be \emph{$(d+1)$-term} if it is isomorphic to a complex in $\K^b(\proj\La)$ whose terms are concentrated in degrees $-d, \dots, 0$.
    \end{enumerate}
\end{definition}

The set of isomorphism classes of basic silting objects in $\K^b(\proj\La)$ is equipped with a partial order $\leq$ defined by \cite[Theorem~2.11]{AI} as follows:
\[P \leq Q :\iff  \Hom (Q, P[i]) = 0 \text{ for all } i>0.\]
The cover relation in this poset can be described via mutation. Let $P = \bigoplus_{i=1}^n P_i$ be a basic silting object with each $P_i$ indecomposable. The \emph{left mutation} of $P$ at $P_i$ is defined as
\[\mu_i^+(P):=\left(\bigoplus_{j\neq i}P_j\right)\oplus P'_i,\]
where $P'_i$ is the cone of a minimal left $\left(\add \bigoplus_{j\neq i}P_j\right)$-approximation of $P_i$. The \emph{right mutation} $\mu_i^-(P)$ is defined dually using a minimal right approximation.

\begin{theorem}[{\cite[Theorem~2.31, Theorem~2.35]{AI}}]
    Both $\mu_i^+(P)$ and $\mu_i^-(P)$ are again silting objects. Moreover, they satisfy the cover relation $\mu_i^-(P)\gtrdot P\gtrdot\mu_i^+(P)$, where $x\gtrdot  y$ means that $x> y$ and there exists no $z$ with $x > z > y$.
\end{theorem}

We say that the left (resp. right) mutation of a $(d+1)$-term silting object $P$ at a direct summand $P_i$ \emph{exists} if $\mu_i^+(P)$ (resp. $\mu_i^-(P)$) is again $(d+1)$-term.

\subsection{Simple-minded collections}\label{sec:smc}

We briefly recall the definition of simple-minded collections in $\D^b(\mod\La)$.

\begin{definition}\label{defn:smc}
    A collection $\X$ of objects in $\D^b(\mod\La)$ is called a \emph{simple-minded collection} if the following conditions hold for all $X,Y\in \X$.
    \begin{enumerate}
        \item[(a)] $\End(X)$ is a division algebra and $\Hom(X, Y)=0$ for $X\neq Y$;
        \item[(b)] $\Hom(X,Y[i]) = 0$ for all $i < 0$;
        \item[(c)] $\thick\X=\D^b(\mod\La)$. 
    \end{enumerate}
    A simple-minded collection $\X$ is called \emph{$(d+1)$-term} if $\X\subseteq \D^{[-d,0]}(\mod\La).$
\end{definition}

The cardinality of any simple-minded collection equals the rank $n$ of $\La$, i.e., the number of pairwise non-isomorphic indecomposable finitely generated projective $\La$-modules \cite{KY}.

We consider simple-minded collections up to isomorphism. Similar to silting objects, simple-minded collections admit a notion of mutation \cite{KY}. Let $\X=\{X_i\}_{1\leq i\leq n}$ be a simple-minded collection. The \emph{left mutation} of $\X$ at $X_i$ is defined as
\[\mu_i^+(\X)=\{X'_j\}_{j\neq i}\cup \{X_i[1]\},\]
where $X'_j$ is the cone of a minimal left $W(X_i)$-approximation of $X_j[-1]$, and $W(X_i)$ is the extension closure of $X_i$ in $\D^b(\mod\La)$. The \emph{right mutation} $\mu_i^-(\X)$ is defined dually. We say that the left (resp. right) mutation of a $(d+1)$-term simple-minded collection $\X$ at $X_i$ \emph{exists} if $\mu_i^+(\X)$ (resp. $\mu_i^-(\X)$) is again $(d+1)$-term.

The following theorem summarizes the connections between the key notions of this section.

\begin{theorem}[{\cite[Theorems~6.1~and~7.12]{KY} and \cite[Theorem~3.3 and Corollary~3.4]{IJ}}]\label{silting-simple}
    There are bijections between the following sets:
    \begin{itemize}
        \item[(1)] basic silting objects in $\K^b(\proj\La)$;
        \item[(2)] bounded $t$-structures on $\D^b(\mod\La)$ with length hearts; 
        \item[(3)] length bounded hearts in $\D^b(\mod\La)$;
        \item[(4)] simple-minded collections in $\D^b(\mod\La)$.
    \end{itemize}
    These bijections commute with left and right mutations. Specifically, a silting object $P$ corresponds to the $t$-structure $(\D^{\leq 0}(P), \D^{\geq 0}(P))$ where
    \[\begin{array}{c}
        \D^{\leq 0}(P)=\{X\in \D^b(\mod\La)\mid \Hom(P,X[i])=0\text{ for all }i>0\},  \\
        \D^{\geq 0}(P)=\{X\in \D^b(\mod\La)\mid \Hom(P,X[i])=0\text{ for all }i<0\}. 
    \end{array}\]
    A $t$-structure is mapped to its heart, and a heart is mapped to its set of isoclasses of simple objects. Moreover, these bijections restrict to:
    \begin{itemize}
        \item[(1')] basic $(d+1)$-term silting objects in $\K^b(\proj\La)$,
        \item[(2')] bounded $t$-structures $(\C^{\leq 0},\C^{\geq 0})$ on $\D^b(\mod\La)$ with length hearts satisfying $\D^{\leq -d}(\mod\La)\subseteq\C^{\leq 0}\subseteq\D^{\leq 0}(\mod\La)$, 
        \item[(3')] length bounded hearts lying in $\D^{[-d,0]}(\mod\La)$, and
        \item[(4')] $(d+1)$-term simple-minded collections in $\D^b(\mod\La)$.
    \end{itemize}
\end{theorem}

\section{Three types of torsion-like classes in extended hearts}\label{S3}

In this section, we establish the basic properties of three notions of torsion-like classes in the $d$-extended heart $\D^{[-d+1,0]}$, namely $d$-FAE closed subcategories (where $d$-FAE stands for $d$-factors and extensions), positive torsion classes, and $s$-torsion classes. When $\D=\D^b(\mod\La)$ and $d=1$, so that $\D^{[-d+1,0]}\simeq\mod\La$, these three notions coincide and recover the usual notion of a torsion class in an abelian category. In general, however, they are distinct (see Figure~\ref{fig:inclusions}).

\subsection{\texorpdfstring{$d$}{d}-FAE and \texorpdfstring{$d$}{d}-SAE closed subcategories}

In this subsection, we study subcategories of $\D^{[-d+1,0]}$ defined by specific closure properties. We begin by recalling the following notions, which provide a natural generalization of the classical concepts of factor objects and subobjects in abelian categories to a higher setting.

\begin{definition}[{\cite[Definition~1.13]{Z}}]
    Let $\X\subseteq \D^{[-d+1,0]}$ and let $m$ be a positive integer. An object $Z\in\D^{[-d+1,0]}$ is called an \emph{$m$-factor} of $\X$ if there exist $m$ triangles
    \[Z_i\to X_i\to Z_{i-1},\ 1\leq i\leq m,\]
    such that $Z_0=Z$, $Z_1,\dots,Z_m\in\D^{[-d+1,0]}$, and $X_1,\dots,X_m\in\X$. We denote by $\Fac_m(\X)$ the full subcategory of $\D^{[-d+1,0]}$ consisting of $m$-factors of $\X$. For convenience, we set $\Fac_0(\X)=\D^{[-d+1,0]}$. We say that $\X$ is \emph{closed under $d$-factors} if $\Fac_d(\X)\subseteq\X$. 

    Dually, we define \emph{$m$-subobjects} of $\X$ and denote the corresponding subcategory by $\Sub_m(\X)$. We say that $\X$ is \emph{closed under $d$-subobjects} if $\Sub_d(\X)\subseteq \X$.
\end{definition}

Note that we have inclusions
\[\begin{array}{c}
     \X\subseteq\cdots \subseteq\Fac_2(\X)\subseteq \Fac_1(\X)\subseteq \Fac_0(\X)=\D^{[-d+1,0]},  \\
     \X\subseteq\cdots \subseteq\Sub_2(\X)\subseteq \Sub_1(\X)\subseteq \Sub_0(\X)=\D^{[-d+1,0]}. 
\end{array}\]

\begin{remark}[{\cite[Example~1.14]{Z}}]\label{exm:d-fac}
    Let $Y\in\D^{[-d+1,0]}$, and fix a non-negative integer $l$ such that $Y[l]\in\D^{[-d+1,0]}$. Let $\X\subseteq\D^{[-d+1,0]}$, and fix a non-negative integer $m$ such that $Y\in\Fac_m(\X)$. Then $Y[l]\in\Fac_{m+l}(\X)$.
\end{remark}

The following construction of $d$-factors will be useful.

\begin{lemma}\label{1.15update}
    Let $\X\subseteq \D^{[-d+1,0]}$. For any $X\in\X$ and any integer $1\leq j\leq d-1$, we have 
    \[\sigma_{\geq -d+1+j}(X)\in \Fac_{d-j+1}(\X)\quad\text{and}\quad\sigma_{\geq -d+1}(X[j])\in \Fac_{d+1}(\X).\]
\end{lemma}

\begin{proof}
    Since $X\in\X\subseteq \D^{[-d+1,0]}$, we have 
    $\sigma_{\leq -d+j}(X)=(\sigma_{\leq 0}(X[-d+j]))[d-j],$ 
    by the definition of truncation and shift. It then follows from Remark~\ref{exm:d-fac} that
    \[\sigma_{\leq -d+j}(X)\in \Fac_{d-j}(\X).\]
    Consider the truncation triangle of $X$
    \[\sigma_{\leq -d+j}(X)\to X\to \sigma_{\geq -d+j+1}(X).\]
    It follows that 
    $\sigma_{\geq -d+1+j}(X)\in \Fac_{d-j+1}(\X).$ 
    Applying Remark~\ref{exm:d-fac} once again, we conclude that
    \[\sigma_{\geq -d+1}(X[j])=(\sigma_{\geq -d+1+j}(X))[j]\in\Fac_{d+1}(\X).\]
\end{proof}

We record the following basic properties of subcategories that are closed under $d$-factors.

\begin{lemma}\label{summand}
    Let $\X\subseteq\D^{[-d+1,0]}$ be closed under $d$-factors. Then:
    \begin{enumerate}
        \item[(1)] $\X$ is closed under direct summands.
        \item[(2)] For any morphism $f\colon X\to Y$ in $\X$, the truncation $\sigma_{\geq -d+1}(C(f))\in\X$.
    \end{enumerate}
\end{lemma}

\begin{proof}
    For (1), let $X=Y\oplus Z\in \X$. Then there are triangles
    \begin{gather*}
        Z\to X\to Y,\\
        Y\to X\to Z.
    \end{gather*}
    Iterating these triangles, we see that $Y\in \Fac_{d}(\X)\subseteq \X$.
    
    For (2), from the triangle $X\xrightarrow{f}Y\xrightarrow{g}C(f)$, we have $C(f)\in\D^{[-d,0]}$. Consider the truncation triangle of $C(f)$:
    \[\sigma_{\leq -d}(C(f))\to C(f)\xrightarrow{h} \sigma_{\geq -d+1}(C(f)),\]
    where $\sigma_{\leq -d}(C(f))\in\D^{[-d,-d]}=\H[d]$ and $\sigma_{\geq -d+1}(C(f))\in\D^{[-d+1,0]}$. Applying the octahedral axiom to the composition $h\circ g$, we obtain a diagram of triangles
    \begin{center}
    \begin{tikzcd}
    & & \sigma_{\geq -d+1}(C(f))[-1] \arrow[r, Rightarrow, no head] \arrow[d] & \sigma_{\geq -d+1}(C(f))[-1] \arrow[d]   \\
    \sigma_{\leq -d}(C(f))[-1]\arrow[r] \arrow[d] & X \arrow[r] \arrow[d, Rightarrow, no head] & Z \arrow[r] \arrow[d] & \sigma_{\leq -d}(C(f)) \arrow[d] \\
    \sigma_{\leq -d}(C(f))[-1]\arrow[r] & X \arrow[r, "f"]& Y \arrow[r, "g"] \arrow[d, "h\circ g"]  & C(f) \arrow[d, "h"]\\
    & & \sigma_{\geq -d+1}(C(f)) \arrow[r, Rightarrow, no head] & \sigma_{\geq -d+1}(C(f))        
    \end{tikzcd}
    \end{center}
    From the triangle in the second row and the triangle in the third column, we obtain
    \[Z\in(X*\sigma_{\leq -d}(C(f)))\cap(\sigma_{\geq -d+1}(C(f))[-1]*Y)\subseteq\D^{\leq 0}\cap \D^{\geq -d+1}=\D^{[-d+1,0]}.\]
    Since $\sigma_{\leq -d}(C(f))[-1]\in\H[d-1]$, by Remark~\ref{exm:d-fac}, we have that
    $\sigma_{\leq -d}(C(f))[-1]\in\Fac_{d-1}(\X).$ 
    The triangle in the second row then yields $Z\in \Fac_d(\X)$. Consequently, the triangle in the third column shows that
    \[\sigma_{\geq -d+1}(C(f))\in \Fac_{d+1}(\X)\subseteq\Fac_d(\X)\subseteq\X.\] 
\end{proof}

We now introduce the first notion of torsion-like classes, namely $d$-FAE closed subcategories, together with their dual notion of $d$-SAE closed subcategories.

\begin{definition}\label{d-FAE}
    A subcategory $\X$ of $\D^{[-d+1,0]}$ is called \emph{$d$-FAE closed} (resp. \emph{$d$-SAE closed}) if it is closed under $d$-factors (resp. $d$-subobjects) and extensions.
\end{definition}

Recall from \cite{AMV} that a subcategory $\mathcal{C}$ of $\mathcal{D}$ is called \emph{suspended} (resp. \emph{cosuspended}) if it is closed under extensions and positive (resp. negative) shifts. Following \cite{Buan}, such subcategories are also referred to as \emph{coresolving} (resp. \emph{resolving}).

\begin{proposition}\label{prop:bi1}
    There is a bijection between
    \begin{itemize}
        \item the set of $d$-FAE closed subcategories of $\D^{[-d+1,0]}$, and
        \item the set of suspended subcategories $\C$ of $\D$ satisfying $\D^{\leq -d}\subseteq\C\subseteq\D^{\leq 0}$,
    \end{itemize}
    given by 
    $\X\mapsto \D^{\leq -d}*\X,$ 
    with inverse 
    $\C\mapsto\C\cap\D^{[-d+1,0]}.$
    
    Dually, there is a bijection between
    \begin{itemize}
        \item the set of $d$-SAE closed subcategories of $\D^{[-d+1,0]}$, and
        \item the set of cosuspended subcategories $\C$ of $\D$ satisfying $\D^{\geq 0}\subseteq\C\subseteq\D^{\geq -d}$,
    \end{itemize}
    given by 
    $\X\mapsto \X[1]*\D^{\geq 0},$ 
    with inverse 
    $\C\mapsto\C[-1]\cap\D^{[-d+1,0]}.$
\end{proposition}

\begin{proof}
    We only prove the first bijection, since the second one can be shown dually.

    We first show that for any suspended subcategory $\C$ of $\D$ satisfying $\D^{\leq -d}\subseteq\C\subseteq\D^{\leq 0}$, the subcategory $\C\cap\D^{[-d+1,0]}$ of $\D^{[-d+1,0]}$ is $d$-FAE closed. Since both $\C$ and $\D^{[-d+1,0]}$ are closed under extensions, their intersection $\C\cap\D^{[-d+1,0]}$ is also closed under extensions. Let $Z\in\Fac_d(\C\cap\D^{[-d+1,0]})$. By definition, we have
    \[Z\in(\C\cap\D^{[-d+1,0]})*(\C\cap\D^{[-d+1,0]})[1]*\cdots*(\C\cap\D^{[-d+1,0]})[d-1]*\D^{[-d+1,0]}[d].\]
    Since $\C$ is closed under positive shifts, $(\C\cap\D^{[-d+1,0]})[i]\subseteq\C$ for all $i\geq 0$. Note also that $\D^{[-d+1,0]}[d]\subseteq\D^{\leq -d}\subseteq\C$. It follows that $Z\in\C$ because $\C$ is closed under extensions. Therefore, $\C\cap\D^{[-d+1,0]}$ is closed under $d$-factors.

    We next show that for any $d$-FAE closed subcategory $\X$ of $\D^{[-d+1,0]}$, the subcategory $\D^{\leq -d}*\X$ of $\D$ is suspended. Note that $\D^{\leq -d}\subseteq\D^{\leq -d}*\X\subseteq\D^{\leq 0}$. To show that $\D^{\leq -d}*\X$ is closed under extensions, it is equivalent to prove
    \[(\D^{\leq -d}*\X)*(\D^{\leq -d}*\X)\subseteq\D^{\leq -d}*\X.\] Since both $\D^{\leq -d}$ and $\X$ are closed under extensions, it is enough to verify 
    \[\X*\D^{\leq -d}\subseteq \D^{\leq -d}*\X.\] 
    Let $Z\in \X*\D^{\leq -d}$. Then there exists a triangle 
    \[X\to Z\to Y\]
    with $X\in \X$ and $Y\in \D^{\leq -d}$. By the octahedral axiom, we obtain the following commutative diagram of triangles:
    \begin{center}
    \begin{tikzcd}
    & \sigma_{\geq -d+1}(Z)[-1] \arrow[r, Rightarrow, no head] \arrow[d] & \sigma_{\geq -d+1}(Z)[-1] \arrow[d] & \\
    Y[-1] \arrow[r] \arrow[d, Rightarrow, no head]& V \arrow[r] \arrow[d] & \sigma_{\leq -d}(Z) \arrow[r] \arrow[d] & Y \arrow[d, Rightarrow, no head] \\
    Y[-1] \arrow[r] & X \arrow[r] \arrow[d] & Z \arrow[r] \arrow[d] & Y \\
    & \sigma_{\geq -d+1}(Z) \arrow[r, Rightarrow, no head] & \sigma_{\geq -d+1}(Z)& 
    \end{tikzcd}
    \end{center}
    By the triangle in the second row, we have
    \[V\in Y[-1]*\sigma_{\leq -d}(Z)\subseteq \D^{\leq -d+1}*\D^{\leq -d}\subseteq \D^{\leq -d+1},\] 
    and by the triangle in the second column,
    \[V\in \sigma_{\geq -d+1}(Z)[-1]*X\subseteq\D^{\geq -d+2}*\D^{\geq -d+1}\subseteq \D^{\geq -d+1}.\]
    Consequently, $V\in \D^{[-d+1, -d+1]}=\H[d-1]$, which implies $V\in \Fac_{d-1}(\X)$ by Remark~\ref{exm:d-fac}. Hence, by the triangle in the second column, we obtain $\sigma_{\geq -d+1}(Z)\in \Fac_d(\X)\subseteq\X$. It follows that $Z\in \D^{\leq -d}*\X$, as desired.
    
    Now, let $Z\in \D^{\leq -d}*\X$. Then $\sigma_{\geq -d+1}(Z)\in\X$. By Lemma~\ref{1.15update}, we have 
    \[\sigma_{\geq -d+1}(Z[1])= \sigma_{\geq -d+1}(\sigma_{\geq -d+1}(Z)[1])\in\X.\]
    Consider the triangle obtained by truncating $Z[1]$:
    \[\sigma_{\leq -d}(Z[1])\to Z[1]\to \sigma_{\geq -d+1}(Z[1]).\]
    Since $\sigma_{\le -d}(Z[1]) \in \D^{\le -d}$ and $\sigma_{\ge -d+1}(Z[1]) \in \X$, it follows that $Z[1]\in \D^{\leq -d}*\X$. Therefore, $\D^{\leq -d}*\X$ is closed under positive shifts, and hence suspended.
    
    Finally, we show that these two maps are mutually inverse. For any $\X\subseteq\D^{[-d+1,0]}$, it is straightforward to see, by taking truncation, that 
    $(\D^{\leq -d}*\X)\cap\D^{[-d+1,0]}=\X.$ 
    Hence, it remains to show that 
    \[\D^{\leq -d}*(\C\cap\D^{[-d+1,0]})=\C\] 
    for any suspended subcategory $\C$ of $\D$ satisfying $\D^{\leq -d}\subseteq\C\subseteq\D^{\leq 0}$. Since both $\D^{\leq -d}$ and $\C\cap\D^{[-d+1,0]}$ are contained in $\C$, and $\C$ is closed under extensions, we have 
    \[\D^{\leq -d}*(\C\cap\D^{[-d+1,0]})\subseteq\C.\] 
    Now, let $Z\in\C\subseteq\D^{\leq 0}$. By taking truncation, there exists a triangle
    \[\sigma_{\leq -d}(Z)\to Z\to \sigma_{\geq -d+1}(Z)\to\sigma_{\leq -d}(Z)[1],\]
    where $\sigma_{\geq -d+1}(Z)\in\D^{[-d+1,0]}$. Since $\sigma_{\leq -d}(Z)\in\D^{\leq -d}\subseteq\C$ and $\C$ is closed under positive shifts, we have $\left(\sigma_{\leq -d}(Z)\right)[1]\in\C$. Moreover, since $\C$ is closed under extensions, it follows that
    \[\sigma_{\geq -d+1}(Z)\in Z*\left(\sigma_{\leq -d}(Z)\right)[1]\subseteq\C,\]
    Hence, $\sigma_{\geq -d+1}(Z)\in \C\cap\D^{[-d+1,0]}$, and consequently $Z\in\D^{\leq -d}*(\C\cap\D^{[-d+1,0]})$. This shows that $\C\subseteq\D^{\leq -d}*(\C\cap\D^{[-d+1,0]})$, and the proof is complete.
\end{proof}

\begin{remark}\label{rmk:smd}
For any $d$-FAE (resp. $d$-SAE) closed subcategory $\X$ of $\D^{[-d+1,0]}$, the full subcategory $\D^{\leq -d}*\X$ (resp. $\X[1]*\D^{\geq 0}$) is closed under direct summands. To see this, suppose that $X\oplus Y\in \D^{\leq -d}*\X$. Then 
\[\sigma_{\geq -d+1}(X)\oplus \sigma_{\geq -d+1}(Y)\cong\sigma_{\geq -d+1}(X\oplus Y)\in \X.\]
By Lemma~\ref{summand}(1), it follows that $\sigma_{\geq -d+1}(X)\in \X$. Hence, 
\[X\in\sigma_{\leq -d}(X)*\sigma_{\geq -d+1}(X)\subseteq \D^{\leq -d}*\X.\]
Combining this with Proposition~\ref{prop:bi1}, we conclude that any suspended (resp. cosuspended) subcategory $\C$ of $\D$ satisfying $\D^{\leq -d}\subseteq\C\subseteq\D^{\leq 0}$ (resp. $\D^{\geq 0}\subseteq\C\subseteq\D^{\geq -d}$) is closed under direct summands.
\end{remark}

\begin{remark}
    It is straightforward to verify that $d$-FAE closed subcategories are closed under intersections. Consequently, for any subcategory $\Y$ of $\D^{[-d+1,0]}$, there exists a smallest $d$-FAE closed subcategory containing $\Y$, denoted by $\phi(\Y)$.
\end{remark}

We have the following explicit description of the smallest $d$-FAE closed subcategory containing a given subcategory; this description will be used in Section~\ref{S6}. We define $\Phi^0(\Y)=\Y$, and for each $i\geq 1$ set \[\Phi^i(\Y)=\Fac_{d}\left(\Phi^{i-1}(\Y)\right)*\Fac_{d}\left(\Phi^{i-1}(\Y)\right).\]
Note that for all $i\geq 1$, we have 
\[\Phi^{i-1}(\Y)\subseteq \Fac_{d}\left(\Phi^{i-1}(\Y)\right)\subseteq \Phi^i(\Y).\]

\begin{lemma}\label{union}
    For any $\Y\subseteq\D^{[-d+1,0]}$, the smallest $d$-FAE closed subcategory $\phi(\Y)$ containing $\Y$ is $\bigcup_{i\geq 0}\Phi^i(\Y)$.
\end{lemma}

\begin{proof}
    Let $\B:=\bigcup_{i\geq 0}\Phi^i(\Y)$. By definition, $\B\subseteq\phi(\Y)$. To prove the reverse inclusion, it suffices to show that $\B$ is $d$-FAE closed.

    Consider a sequence of triangles
    \[Z_j\to X_j\to Z_{j-1},\ 1\leq j\leq d,\]
    with $Z_0,Z_1,\dots,Z_d\in\D^{[-d+1,0]}$ and $X_1,\dots, X_{d}\in \B$. Then there exists $s\geq 0$ such that $X_j\in \Phi^s(\Y)$ for all $j$. By definition, the resulting object $Z_0$ obtained from this sequence of triangles belongs to $\Fac_d(\Phi^s(\Y))\subseteq \Phi^{s+1}(\Y)\subseteq\B$. Hence, $\B$ is closed under $d$-factors.

    Finally, suppose there is a triangle $L\to M\to N$ with $L,N\in \B$. As before, there exists $s\geq 0$ such that $L,N\in \Phi^s(\Y)$. It follows that 
    \[M\in \Phi^s(\Y)*\Phi^s(\Y)\subseteq \Phi^{s+1}(\Y)\subseteq\B.\]
    Therefore, $\B$ is closed under extensions.
\end{proof}

\subsection{Positive torsion/torsion-free classes}

For any $\X\subseteq\D^{[-d+1,0]}$, we define the following two full subcategories of $\D^{[-d+1,0]}$:
\[{}^{\perp}\X=\{Z\in\D^{[-d+1,0]}\mid \Hom(Z,\X)=0\},\]
\[\X^{\perp}=\{Z\in\D^{[-d+1,0]}\mid \Hom(\X,Z)=0\}.\]
The second notion of torsion-like classes is that of positive torsion classes, introduced in \cite{Gu1}.

\begin{definition}[{\cite[Definition~3.14]{Gu1}}]
    A pair of subcategories $(\T,\F)$ of $\D^{[-d+1,0]}$ is called a \emph{torsion pair} if $\T={}^\perp\F$ and $\F=\T^\perp$. A torsion pair is called \emph{positive} if 
    \[\Hom(\T,\F[j])=0\text{ for all integers }j\leq 0.\]
    
    In a (positive) torsion pair $(\T,\F)$, the subcategory $\T$ is called a \emph{(positive) torsion class} and $\F$ is called a \emph{(positive) torsion-free class}.
\end{definition}

It follows directly from the definition that any torsion class and any torsion-free class are closed under extensions. The following lemma guarantees that any positive torsion (resp. torsion-free) class is $d$-FAE (resp. $d$-SAE) closed.

\begin{lemma}\label{positive}
Let $(\T,\F)$ be a torsion pair in $\D^{[-d+1,0]}$. The following are equivalent. 
\begin{enumerate}
    \item $\T$ is a positive torsion class. 
    \item For all morphisms $f$ in $\T$, the truncation $\sigma_{\geq -d+1}(C(f))$ belongs to $\T$.
    \item $\T$ is $d$-FAE closed.
    \item $\F$ is a positive torsion-free class.
    \item For all morphisms $g$ in $\F$, the truncation $\sigma_{\leq 0}\left(C(g)[-1]\right)$ belongs to $\F$.
    \item $\F$ is $d$-SAE closed.
    \item $\Hom(\T,\F[-1])=0$.
\end{enumerate}
\end{lemma}

\begin{proof}
   The equivalence of (1), (2), (4), (5), and (7) follows from \cite[Lemma~3.16]{Gu1}, noting that although the original proof is formulated for $\D^{[-d+1,0]}(\mod\La)$, it generalizes directly to the general case of $\D^{[-d+1,0]}$. The equivalence of (1), (3), (4), and (6) is then deduced from the arguments in the proof of \cite[Proposition~1.18]{Z}.
\end{proof}

\begin{remark}
It is known that torsion classes are closed under intersections \cite[Theorem~4.10]{Gu1}. Lemma~\ref{positive} then implies that positive torsion classes are closed under intersections as well. Consequently, for any subcategory $\Y$ of $\D^{[-d+1,0]}$, there exists a smallest positive torsion class containing $\Y$, denoted by $T(\Y)$.
\end{remark}

We now present several results that allow one to construct (positive) torsion and torsion-free classes from arbitrary subcategories; 
in particular, the smallest positive torsion class mentioned above can be constructed in a direct manner, rather than merely as the intersection of all positive torsion classes containing $\Y$.

\begin{lemma}\label{perp}
    Let $\Y$ be a full subcategory of $\D^{[-d+1,0]}$. Then the following hold.
    \begin{enumerate}
        \item $\Y^\perp$ is a torsion-free class.
        \item ${}^\perp\Y$ is a torsion class.
    \end{enumerate}
\end{lemma}

\begin{proof}
    We only prove (1), as the proof of (2) is dual. By definition, we have $\Y\subseteq{}^\perp(\Y^\perp)$. This immediately implies $\Y^\perp\supseteq({}^\perp(\Y^\perp))^\perp$. Conversely, for any $Y\in\Y^\perp$ and $V\in{}^\perp(\Y^\perp)$, we have $\Hom(V,Y)=0$, which shows that $Y\in ({}^\perp(\Y^\perp))^\perp$. Hence, $\Y^\perp\subseteq({}^\perp(\Y^\perp))^\perp$. Combining these two inclusions, we obtain $\Y^\perp=({}^\perp(\Y^\perp))^\perp$, and hence $\Y^\perp$ is a torsion-free class.
\end{proof}

For $\Y\subseteq \D^{[-d+1,0]}$, we define the following subcategories of $\D^{[-d+1,0]}$:
\[{}^{\perp_{\leq 0}}\Y=\{Z\in\D^{[-d+1,0]}\mid \Hom(Z,\Y[j])=0,\ \forall j\leq 0\},\]
\[\Y^{\perp_{\leq 0}}=\{Z\in\D^{[-d+1,0]}\mid \Hom(\Y,Z[j])=0,\ \forall j\leq 0\}.\]

\begin{lemma}\label{1.21update}
    For any subcategory $\Y$ of $\D^{[-d+1,0]}$, we have 
    \begin{equation}\label{eqs1}
        \Y^{\perp_{\leq 0}}=(\Fac_{d}\Y)^{\perp_{\leq 0}}=(\Fac_{d}\Y)^\perp,
    \end{equation}
    \begin{equation}\label{eqs2}
    {}^{\perp_{\leq 0}}\Y={}^{\perp_{\leq 0}}(\Sub_{d}\Y)={}^\perp(\Sub_{d}\Y).
    \end{equation}
\end{lemma}

\begin{proof}
    We only show the equalities in \eqref{eqs1}, since those in \eqref{eqs2} can be proved similarly. By \cite[Lemma~1.19]{Z}, we have $\Y^{\perp_{\leq 0}}\subseteq(\Fac_{d}\Y)^{\perp_{\leq 0}}$. It is clear that $(\Fac_{d}\Y)^{\perp_{\leq 0}}\subseteq(\Fac_{d}\Y)^\perp$. Thus, it suffices to show  $(\Fac_{d}\Y)^\perp\subseteq\Y^{\perp_{\leq 0}}$. 
    
    Let $X\in(\Fac_{d}\Y)^\perp$. For any $Y\in\Y$ and any $j\geq 0$, consider the truncation triangle of $Y[j]$:
    \[\sigma_{\leq -d}(Y[j])\to Y[j]\to \sigma_{\geq -d+1}(Y[j])\to \sigma_{\leq -d}(Y[j])[1].\]
    Applying $\Hom(-,X)$ to this triangle yields an isomorphism
    \[\Hom(Y[j],X)\cong\Hom(\sigma_{\geq -d+1}(Y[j]),X).\]
    By Lemma~\ref{1.15update}, we have $\sigma_{\geq -d+1}(Y[j])\in\Fac_{d+1}(\Y)\subseteq\Fac_{d}(\Y)$, and hence
    \[\Hom(\sigma_{\geq -d+1}(Y[j]),X)=0.\]
    It follows that $\Hom(Y[j],X)=0$ for all $j\geq 0$, so $X\in\Y^{\perp_{\leq 0}}$.
\end{proof}

We obtain the following description of the smallest positive torsion class containing a given subcategory, which will be used in Sections~\ref{S6}~and~\ref{S7}.

\begin{proposition}\label{perp:tor}
    For any subcategory $\Y$ of $\D^{[-d+1,0]}$, the following hold.
    \begin{enumerate}
        \item[(1)] The subcategory $^{\perp_{\leq 0}}\Y$ is a positive torsion class. Its corresponding torsion-free class $({}^{\perp_{\leq 0}}\Y)^\perp$ is the smallest positive torsion-free class containing $\Y$.
        \item[(2)] The subcategory $\Y^{\perp_{\leq 0}}$ is a positive torsion-free class. Its corresponding torsion class ${}^\perp(\Y^{\perp_{\leq 0}})$ is the smallest positive torsion class containing $\Y$.
    \end{enumerate}
\end{proposition}

\begin{proof}
    We only prove (1), while the proof of (2) is dual. By Lemma~\ref{1.21update}, we have ${}^{\perp_{\leq 0}}\Y={}^\perp(\Sub_{d}\Y)$, and hence Lemma~\ref{perp} implies that ${}^{\perp_{\leq 0}}\Y$ is a torsion class. 
    
    Let $Z\in\Fac_d({}^{\perp_{\leq 0}}\Y)$. By definition, there exist triangles
    \[Z_i\to X_i\to Z_{i-1},\ 1\leq i\leq d,\]
    with $Z_0=Z$, $Z_1,\dots,Z_d\in\D^{[-d+1,0]}$, and $X_1,\dots,X_d\in{}^{\perp_{\leq 0}}\Y$. For any $Y\in\Y$ and any $j\leq 0$, applying $\Hom(-,Y[j])$ to these triangles yields exact sequences for $1\leq i\leq d$:
    \[\Hom(X_i[1],Y[j])\to\Hom(Z_{i}[1],Y[j])\to\Hom(Z_{i-1},Y[j])\to\Hom(X_i,Y[j]).\]
    Since $X_i\in{}^{\perp_{\leq 0}}\Y$, the first and fourth terms vanish. It follows that  \[\Hom(Z_{i}[1],Y[j])\cong\Hom(Z_{i-1},Y[j]).\]
    Iterating this for $i = 1, \dots, d$, we obtain
    \[\Hom(Z,Y[j])\cong \Hom(Z_{d},Y[j-d])=0,\]
    which shows that $Z\in{}^{\perp_{\leq 0}}\Y$. Therefore, $^{\perp_{\leq 0}}\Y$ is closed under $d$-factors. By Lemma~\ref{positive}, it follows that ${}^{\perp_{\leq 0}}\Y$ is a positive torsion class, and hence $({}^{\perp_{\leq 0}}\Y,({}^{\perp_{\leq 0}}\Y)^\perp)$ is a positive torsion pair. 
    
    Finally, by definition, we have $\Y \subseteq ({}^{\perp_{\leq 0}}\Y)^\perp$. If $\F$ is any positive torsion-free class containing $\Y$, then ${}^\perp\F\subseteq{}^{\perp_{\leq 0}}\Y$, which implies 
    \[({}^{\perp_{\leq 0}}\Y)^\perp\subseteq ({}^\perp\F)^\perp=\F.\] Hence, $({}^{\perp_{\leq 0}}\Y)^\perp$ is indeed the smallest positive torsion-free class containing $\Y$.
\end{proof}

\begin{remark}\label{fac->torf}
    If $\Y$ is already closed under $d$-factors, Lemma~\ref{1.21update} yields the simplified identity $\Y^\perp=\Y^{\perp_{\leq 0}}$. Hence, by Proposition~\ref{perp:tor}, $\Y^\perp$ is a positive torsion-free class, and its corresponding torsion class ${}^\perp(\Y^\perp)$ is the smallest positive torsion class containing $\Y$.
\end{remark}

\subsection{\texorpdfstring{$s$}{s}-torsion/torsion-free classes}

The third notion of torsion-like classes is given by $s$-torsion classes, introduced in \cite{AET} for extriangulated categories with negative first extensions.
\begin{definition}[{\cite[Definition~3.1]{AET}}]
    An \emph{$s$-torsion pair} $(\T,\F)$ in $\D^{[-d+1,0]}$ is a pair of subcategories such that
    \begin{enumerate}
        \item $\D^{[-d+1,0]}=\T*\F$;
        \item $\Hom(\T,\F)=0$;
        \item $\Hom(\T,\F[-1])=0$.
    \end{enumerate}
    Here, $\T$ is called an \emph{$s$-torsion class}, and $\F$ is called an \emph{$s$-torsion-free class}.
\end{definition}

By \cite[Proposition~3.2]{AET}, every $s$-torsion pair is a positive torsion pair. Hence, each $s$-torsion pair is uniquely determined by its $s$-torsion class (or equivalently, by its $s$-torsion-free class). The following property relating positive torsion classes and $s$-torsion classes will be crucial in Section~\ref{S7}.

\begin{lemma}\label{no-middle}
Let $\T$ be an $s$-torsion class and $\T'$ a positive torsion class in $\D^{[-d+1,0]}$. Suppose that $\T' \subseteq \T$ and $\T \cap \T'^\perp = 0$. Then $\T = \T'$.
\end{lemma}

\begin{proof}
    We aim to show that $\T'^\perp \subseteq \T^\perp$, which immediately implies $\T = \T'$. Let $X\in \T'^\perp$. Since $\T$ is an $s$-torsion class, there exists a triangle 
    \[X_1\to X\to X_2\]
    with $X_1\in \T$ and $X_2\in \T^\perp\subseteq \T'^\perp$. Since $\T'$ is a positive torsion class, by Lemma~\ref{positive}, $\T'^\perp$ is closed under truncation of cocones. Hence $X_1\in \T'^\perp$. Together with the assumption $\T \cap \T'^\perp = 0$, this forces $X_1 = 0$. Therefore, $X\cong X_2\in \T^\perp$, as desired.
\end{proof}

Every $s$-torsion (resp. $s$-torsion-free) class is $d$-FAE (resp. $d$-SAE) closed. The bijections in Proposition~\ref{prop:bi1} restrict to a bijection between $s$-torsion pairs in $\D^{[-d+1,0]}$ and bounded $t$-structures with left parts lying between $\D^{\leq -d}$ and $\D^{\leq 0}$.

\begin{proposition}\label{prop:bi2}
    There is a bijection between
    \begin{itemize}
        \item the set of $s$-torsion pairs in $\D^{[-d+1,0]}$, and
        \item the set of bounded $t$-structures $(\C^{\leq 0},\C^{\geq 0})$ on $\D$ satisfying $\D^{\leq -d}\subseteq \C^{\leq 0}\subseteq \D^{\leq 0}$,
    \end{itemize}
    given by 
    \[(\T,\F)\mapsto (\D^{\leq -d}*\T, \F[1]*\D^{\geq 0})\]
    with inverse 
    \[(\C^{\leq 0},\C^{\geq 0})\mapsto (\C^{\leq 0}\cap \D^{[-d+1,0]}, \C^{\geq 1}\cap\D^{[-d+1,0]}).\]
    Moreover, let $(\T,\F)$ be an $s$-torsion pair in $\D^{[-d+1,0]}$, and let $(\C^{\leq 0},\C^{\geq 0})$ be the corresponding $t$-structure on $\D$. Then the following hold.
    \begin{enumerate}
        \item[(1)] $(\F[d],\T)$ is an $s$-torsion pair in the $d$-extended heart $\C^{[-d+1,0]}$ of $(\C^{\leq 0},\C^{\geq 0})$.
        \item[(2)] The bounded $t$-structure on $\D$ corresponding to $(\F[d],\T)$ is $(\D^{\leq-d},\D^{\geq-d})$.
        \item[(3)] $\T=\D^{[-d+1,0]}\cap\C^{[-d+1,0]}$ and $\F=\D^{[-d+1,0]}\cap\left(\C^{[-d+1,0]}[-d]\right)$.
        \item[(4)] For any $X\in \C^{[-d+1,0]}$, the truncation triangle of $X$
        \[\sigma_{\leq -d}(X)\to X\to \sigma_{\geq -d+1}(X) \]
        coincides with the canonical triangle of $X$ with respect to the $s$-torsion pair $(\F[d],\T)$, that is, $\sigma_{\leq -d}(X)\in\F[d]$ and $\sigma_{\geq -d+1}(X)\in\T$.
    \end{enumerate}
\end{proposition}

\begin{proof}
    This proposition agrees with \cite[Proposition~1.9 and Theorem~1.12]{Z}, except for the assertions in (2) and (3). The assertion in (2) can be verified by a straightforward computation: \[\C^{\leq -d}*\F[d]=(\D^{\leq -d}*\T)[d]*\F[d]=\D^{\leq 0}[d]=\D^{\leq-d}.\]
    
    For (3), it follows from (1) that $\T\subseteq\D^{[-d+1,0]}\cap\C^{[-d+1,0]}$. Conversely, let $X\in\D^{[-d+1,0]}\cap\C^{[-d+1,0]}=(\T*\F)\cap(\F[d]*\T)$. Then there exists a triangle 
    \[Y\to X\xrightarrow{f} Z\]
    with $Y\in\T$ and $Z\in\F$. Since $\Hom(\F[d],\F)=0=\Hom(\T,\F)$, we must have $f=0$. Hence $X$ is a direct summand of $Y$, and therefore $X \in \T$. This proves the first equality in (3). The second one can be shown similarly.
    
    For (4), take the canonical triangle of $X$ with respect to $(\F[d],\T)$
    \[X'\to X\to X'',\]
    with $X'\in\F[d]$ and $X''\in\T$. It follows from (2) that $X'\in\D^{\leq -d}$ and $X''\in\D^{\geq -d+1}$. Therefore, this triangle coincides with the truncation triangle of $X$.
\end{proof}

By Theorem~\ref{silting-simple}, the above bijection induces the following correspondence.

\begin{corollary}\label{bi2cor}
    There is a bijection between
    \begin{itemize}
        \item the set of $s$-torsion pairs in $\D^{[-d+1,0]}$, and
        \item the set of bounded hearts on $\D$ lying in $\D^{[-d,0]}$,
    \end{itemize}
    given by 
    \[(\T,\F)\mapsto (\D^{\leq -d}*\T)\cap(\F[1]*\D^{\geq 0})\]
    with inverse 
    \[\H'\mapsto (d\text{-}\H'\cap \D^{[-d+1,0]}, d\text{-}\H'[-d]\cap\D^{[-d+1,0]}).\]
\end{corollary}

The following proposition describes the equivalence between contravariantly finite $d$-FAE closed subcategories, contravariantly finite positive torsion classes, and $s$-torsion classes, under standard finiteness assumptions on $\D$.

\begin{proposition}\label{s-tors}
    Suppose that $\D$ is $K$-linear, $\Hom$-finite, and Krull-Schmidt. Let $\T\subseteq\D^{[-d+1,0]}$ be $d$-FAE closed. Then the following hold.
   \begin{enumerate}
       \item[(1)] $\T$ is a positive torsion class if and only if it is contravariantly finite in ${}^\perp(\T^\perp)$.
       \item[(2)] $\T$ is an $s$-torsion class if and only if it is contravariantly finite in $\D^{[-d+1,0]}$.
   \end{enumerate}
   In particular, the notions of functorially finite $d$-FAE closed subcategories, functorially finite positive torsion classes, and functorially finite $s$-torsion classes coincide.
\end{proposition}
\begin{proof}
    By definition, if $\T$ is a positive torsion class, then $\T={}^\perp(\T^\perp)$, and if $\T$ is an $s$-torsion class, then $\T$ is contravariantly finite in $\D^{[-d+1,0]}$. Thus, the forward implications in (1) and (2) are clear. 
    
    Now suppose $\T$ is a $d$-FAE closed subcategory which is contravariantly finite in ${}^\perp(\T^\perp)$. We will show that ${}^\perp(\T^\perp)\subseteq \T*\T^\perp$ which will imply that ${}^\perp(\T^\perp)\subseteq \T$, and hence, ${}^\perp(\T^\perp)= \T$.
    
    Let $Z\in {}^\perp(\T^\perp)$. Since $\T$ is contravariantly finite in ${}^\perp(\T^\perp)$ and $\D$ is $K$-linear, $\Hom$-finite, and Krull-Schmidt, there exists a minimal right $\T$-approximation $f: T\to Z$ of $Z$. By the triangulated version of Wakamatsu’s Lemma (\cite[Lemma~2.1]{J}), we have $\Hom(\T,C(f))=0$. It remains to show that $C(f)\in \D^{[-d+1,0]}$, which will imply $C(f)\in \T^\perp$ and complete the proof. 
    
    By an argument similar to that in the proof of Lemma~\ref{summand}(2), we obtain the following commutative diagram of triangles:
    \begin{center}
        \begin{tikzcd}
        T \arrow[r] \arrow[d, equal] & T' \arrow[r] \arrow[d, "f'"]& \sigma_{\leq -d}(C(f)) \arrow[d] \\
        T \arrow[r, "f"]& Z \arrow[r] \arrow[d] & C(f) \arrow[d]\\ 
        & \sigma_{\geq -d+1}(C(f)) \arrow[r, equal] & \sigma_{\geq -d+1}(C(f))        
    \end{tikzcd}
    \end{center}
    with $T'\in \Fac_{d}(\T)$. Since $\T$ is closed under $d$-factors, we have $T'\in\T$. Moreover, since $f$ is a right $\T$-approximation of $Z$, so is the induced morphism $f':T'\to Z$. As $f$ is right minimal, the morphism $f'$ decomposes as a direct sum of $f$ and a morphism of the form $V\to 0$. Hence, we have an isomorphism $\sigma_{\geq -d+1}(C(f))\cong C(f)\oplus V[1]$, which implies that $C(f)\in \D^{[-d+1,0]}$.
    
    Thus, $(\T,\T^\perp)$ is a torsion pair, and hence, by Lemma \ref{positive}, a positive torsion pair. 
    
    Now suppose $\T$ is a $d$-FAE closed subcategory which is contravariantly finite in $\D^{[-d+1,0]}$. Then by (1) it is a positive torsion class. Moreover, using the same argument as above, we can show that $\D^{[-d+1,0]}\subseteq \T*\T^\perp$. Hence, $\T$ is an $s$-torsion class. 
\end{proof}

In summary, Figure~\ref{fig:inclusions} illustrates the inclusions and equalities among the three types of torsion-like classes of $\D^{[-d+1,0]}$. Here, the equalities in the first two rows hold under the assumption that $\D$ is $K$-linear, $\Hom$-finite, and Krull-Schmidt. When $d=1$ and $\D=\D^b(\mod\La)$ equipped with its standard $t$-structure $(\D^{\leq 0}(\mod\La),\D^{\geq 0}(\mod\La))$, the inclusions in the last row also become equalities, since in this case $\D^{[-d+1,0]}\simeq \mod\La$. For general $d$, in the $d$-extended module category $\emod{d}{\La}$, the first inclusion in the last row is known to be proper (see \cite[Example~5.3]{Gu1}). Below we provide an example to show that the second inclusion is also proper. 

\begin{figure}[htpb]
    \centering
    \begin{tikzpicture}
    \draw(-4.85,2)node{$\left\{\begin{matrix}\text{functorially finite}\\ \text{$s$-torsion classes}\end{matrix}\right\}$};
    \draw(-.2,2)node{$\left\{\begin{matrix}\text{functorially finite}\\ \text{positive torsion classes}\end{matrix}\right\}$};
    \draw(5.4,2)node{$\left\{\begin{matrix}\text{functorially finite}\\ \text{$d$-FAE closed subcategories}\end{matrix}\right\}$};
    \draw(-2.8,2)node{$=$} (2.4,2)node{$=$};
    \draw(-4.85,1)node[rotate=-90]{$\subseteq$} (-.2,1)node[rotate=-90]{$\subseteq$} (5.4,1)node[rotate=-90]{$\subseteq$};
    \draw(-4.85,0)node{$\left\{\text{$s$-torsion classes}\right\}$};
    \draw(-.2,0)node{$\left\{\begin{matrix}\text{contravariantly finite}\\ \text{positive torsion classes}\end{matrix}\right\}$};
    \draw(5.4,0)node{$\left\{\begin{matrix}\text{contravariantly finite}\\ \text{$d$-FAE closed subcategories}\end{matrix}\right\}$};
    \draw(-2.8,0)node{$=$} (2.4,0)node{$=$};
    \draw(-4.85,-1)node[rotate=-90]{$=$} (-.2,-1)node[rotate=-90]{$\subseteq$} (5.4,-1)node[rotate=-90]{$\subseteq$};
    \draw(-4.85,-1.7)node{$\left\{\text{$s$-torsion classes}\right\}$};
    \draw(-.2,-1.7)node{$\left\{\begin{matrix}\text{positive torsion classes}\end{matrix}\right\}$};
    \draw(5.4,-1.7)node{$\left\{\text{$d$-FAE closed subcategories}\right\}$};
    \draw(-2.8,-1.7)node{$\subseteq$} (2.4,-1.7)node{$\subseteq$};
\end{tikzpicture}
    \caption{Inclusions and equalities among various torsion-like classes}
    \label{fig:inclusions}
\end{figure}

\begin{example}\label{d-FAE-not-positive}
    Let $Q$ be the Kronecker quiver
    \[Q=\begin{tikzcd}
    1 \arrow[r, shift left] \arrow[r, shift right] & 2
    \end{tikzcd}\]
    and $\La=KQ$. Let $\mathcal{P}$, $\mathcal{R}$, and $\mathcal{I}$ denote the preprojective, regular, and postinjective components of the Auslander-Reiten quiver of $KQ$, respectively. Let $d=2$. Using the classification of positive torsion classes of $\emod{2}{\La}$ in \cite[Example~5.3]{Gu1}, we know that $\T=\add (\mathcal{R}\bigcup\mathcal{R}[1]\bigcup\mathcal{I}[1])$ is not a positive torsion class. We claim that $\T$ is a $2$-FAE closed subcategory. Using Proposition \ref{prop:bi1}, it is enough to show that $\D^{\leq -2}*\T$ is a suspended subcategory. It is easy to see that it is closed under positive shifts. We show that it is closed under extensions. Let $\X:=\add \mathcal{R}$ and $\Z:=\D^{\leq -2}* \add (\mathcal{R}[1]\bigcup\mathcal{I}[1])$. Then we have that $\X$ and $\Z$ are closed under extensions and $\D^{\leq -2}*\T = \Z*\X$. To show that $\D^{\leq -2}*\T$ is closed under extensions, it is enough to show that $\X*\Z \subseteq \Z*\X=\add (\X\bigcup\Z)$. 

Suppose there exists a triangle $X\to Y\to Z$ with $X\in\X$ and $Z\in \Z$. Then $Y\in \D^{\leq -1}*\add (\mathcal{R}\bigcup\mathcal{I})$. Let $Y=Y_1\oplus Y_2$ with $Y_1\in \add (\mathcal{P}[1]\bigcup\mathcal{I})$ and $Y_2\in \Z*\X=\add (\X\bigcup\Z)$. Then there exists the following diagram of triangles

    \begin{center}
    \begin{tikzcd}
    & Y_2 \arrow[r, Rightarrow, no head] \arrow[d] & Y_2 \arrow[d]  \\
    X \arrow[r] \arrow[d, Rightarrow, no head]& Y_1\oplus Y_2 \arrow[r] \arrow[d] & Z \arrow[d]  \\
    X \arrow[r] & Y_1 \arrow[r]  & C 
    \end{tikzcd}
    \end{center}
Since $Y_2[1]\in \Z$, we get that $C\in \Z$. Moreover, the bottom triangle implies that $C\in Y_1*X[1]\subseteq \add(\mathcal{R}[1] \bigcup  \mathcal{I}\bigcup \mathcal{P}[1])$. Thus, we get that $C\in \add(\mathcal{R}[1])$. Since $Y_1$ is a cone of a map from $\add \mathcal{R}$ to $\add \mathcal{R}$ and $\add \mathcal{R}$ is closed under kernel and cokernels, we get that $\Ho^0(Y_1), \Ho^{-1}(Y_1)\in \add \mathcal{R}$. Since we assumed $Y_1\in \add (\mathcal{P}[1]\bigcup\mathcal{I})$, we get that $Y_1=0$ and $Y\in \Z*\X$.
\end{example}

\section{Wide subcategories and semibricks}\label{S4}

In this section, we generalize the notions of semibricks and wide subcategories, together with the correspondence between them, from abelian categories to the $d$-extended heart $\D^{[-d+1,0]}$.

Recall that a subcategory of an abelian category is said to be \emph{wide} if it is closed under kernels, cokernels, and extensions. This notion was generalized to exact categories in \cite{E}.

\begin{definition}[{\cite[Definition~2.4]{E}}]\label{wide-exact}
    Let $\C$ be an exact category. A subcategory $\W\subseteq\C$ is called \emph{wide} if the following conditions are satisfied:
    \begin{enumerate}
    \item $\W$ is closed under extensions, that is, for any conflation $L\rightarrowtail M \twoheadrightarrow N$ in $\C$, if $L$ and $N$ belong to $\W$, then so does $M$.
    \item $\W$ is an abelian category.
    \item The inclusion functor $\W \hookrightarrow \C$ is exact, that is, every short exact sequence $0\to L\to M\to N\to 0$ in $\W$ is a conflation in $\C$.
\end{enumerate}
\end{definition}

It was shown in \cite[Theorem~2.5]{E} that there exists a bijection between the class of semibricks in $\C$ and the class of length wide subcategories of $\C$. This bijection is given by assigning to a semibrick its extension closure, while its inverse sends a wide subcategory to its collection of simple objects. This result generalizes the classical correspondence for abelian categories established in \cite{R}.

An exact category $\C$ is called \emph{weakly idempotent complete} if for two composable morphisms $f$ and $g$, whenever $g\circ f$ is an inflation (resp. a deflation), then $f$ is an inflation (resp. $g$ is a deflation). Recall that a full subcategory $\B$ of $\C$ is called \emph{Serre} if for any conflation
$X\rightarrowtail Y\twoheadrightarrow Z$ in $\C$, one has $Y\in\B$ if and only if both $X\in\B$ and $Z\in\B$.

The following proposition provides a canonical construction of a wide subcategory from any exact category by considering the collection of objects to which all incoming morphisms are admissible. This result is of fundamental importance as it provides the template for constructing wide subcategories within $\D^{[-d+1,0]}$ in the next section. While we expect this general result to be known, we were unable to find a proof in the literature and thus include one here for completeness.

\begin{proposition}\label{wide-in-exact}
    Let $\C$ be a weakly idempotent complete exact category. Then the full subcategory
    \[\C':=\{X\in \C\mid \text{any morphism in $\C$ to $X$ is admissible}\}\]
    is a Serre and wide subcategory of $\C$.
\end{proposition}

\begin{proof}
    We divide the proof into several steps. Let $X\Rightarrowtail{f} Y\Twoheadrightarrow{g} Z$ be a conflation of $\C$.
    
    \textbf{(i) If $X,Z\in\C'$, then $Y\in\C'$.} Let $m: A\to Y$ be a morphism in $\C$. Since $Z\in\C'$, the morphism $g\circ m:A\to Z$ is admissible. It follows that there exist a deflation $p:A\twoheadrightarrow C$ and an inflation $p':C\rightarrowtail Z$ such that $g\circ m=p'\circ p$. Let $i=\ker(p):K\rightarrowtail A$. Then $g\circ m\circ i=0$, which implies that there exists a morphism $a: K\to X$ such that $f\circ a=m\circ i$. See the following commutative diagram, where each row is a conflation.
    \begin{center}
    \begin{tikzcd}
    K\arrow[d, "a"'] \arrow[r, "i", tail] & A \arrow[r, "p", two heads] \arrow[d, "m"] & C \arrow[d, "p'", tail] \\
    X \arrow[r, "f", tail] & Y \arrow[r, "g", two heads] & Z
    \end{tikzcd}
    \end{center}
    Since $X\in \C'$, the morphism $a$ is admissible. Hence there exist a deflation $j_1:K\twoheadrightarrow X'$ and an inflation $j_2:X'\rightarrowtail X$ such that $a=j_2\circ j_1$. Consider the pushout of the morphisms $i$ and $j_1$, as shown below, where $q$ is a deflation since it is the pushout of a deflation along an inflation (cf. the dual of \cite[Proposition~2.15]{Bu}):
    \begin{center}
    \begin{tikzcd} 	
   
    K \arrow[d, "j_1"', two heads] \arrow[r, "i", tail]                  & A \arrow[d, "q", two heads] \arrow[rdd, "m", bend left] &   \\
    X' \arrow[r, "j"', tail] \arrow[rrd, "f\circ j_2"', bend right] & M \arrow[rd, "\exists!\ t", dashed]\arrow[ul, phantom, "\ulcorner", very near start]&   \\  &   & Y
    \end{tikzcd}
    \end{center}
    Since $m\circ i=f\circ a=f\circ j_2\circ j_1$, the universal property of pushouts yields a unique morphism $t:M\to Y$ such that 
    \[t\circ j=f\circ j_2 \quad\text{and}\quad t\circ q=m.\]
    We claim that $t$ is an inflation. Indeed, since $p\circ i=0$, the universal property of pushouts yields a unique morphism $s: M\to C$ such that the following diagram commutes:
    \begin{center}
        \begin{tikzcd}
        K \arrow[d, "j_1"', two heads] \arrow[r, "i", tail]               & A \arrow[d, "q", two heads] \arrow[rdd, "p", two heads, bend left] &   &   \\
        X' \arrow[r, "j"', tail] \arrow[rrd, "0"', bend right] & M \arrow[rd, "\exists!\ s", dashed]                                &   &   \\&   & C &  
        \end{tikzcd}
    \end{center}
    By \cite[Remark~2.13]{Bu}, one has $s=\coker(j)$, so the first row of the following diagram is a conflation:
    \begin{center}
        \begin{tikzcd}
        X' \arrow[r, "j", tail] \arrow[d, "j_2"', tail] & M \arrow[r, "s", two heads] \arrow[d, "t"] & C \arrow[d, "p'", tail] \\
        X \arrow[r, "f"', tail]                       & Y \arrow[r, "g"', two heads]               & Z                    
        \end{tikzcd}
    \end{center}
    The commutativity $p'\circ s=g\circ t$ follows from
    \[p'\circ s\circ q=p'\circ p=g\circ m=g\circ t\circ q.\]
    By \cite[Corollary~3.2]{Bu}, this implies that $t$ is an inflation. Hence $m=t\circ q$ is admissible, and therefore $Y\in\C'$.
    
    \textbf{(ii) If $Y\in\C'$, then $X\in\C'$.} For any morphism $k:A\to X$, since $Y\in\C'$, the morphism $f\circ k\colon A\to Y$ is admissible. It follows that there exist a deflation $l:A\twoheadrightarrow B$ and an inflation $m:B\rightarrowtail Y$ such that $f\circ k=m\circ l$. Since $g\circ m\circ l=g\circ f\circ k=0$, it follows that $g\circ m=0$. See the following diagram.
    \begin{center}
        \begin{tikzcd}[row sep=scriptsize]
        A \arrow[r, Rightarrow, no head] \arrow[dd, "k"'] & A \arrow[d, "l", two heads]&   \\& B \arrow[d, "m", tail] \arrow[ld, "s"', dashed] &   \\
        X \arrow[r, "f"', tail] & Y \arrow[r, "g"', two heads]& Z
        \end{tikzcd}
    \end{center}
    Thus there exists a morphism $s:B\to X$ such that $f\circ s=m$. Since $\C$ is weakly idempotent complete, $s$ is an inflation. Moreover,
    $f\circ s\circ l=m\circ l=f\circ k$ implies $s\circ l=k$. Hence $k$ is admissible, and therefore $X\in \C'$.

    \textbf{(iii) If $Y\in\C'$, then $Z\in\C'$.} Given a morphism $a\colon A\to Z$, consider the pullback of morphisms $g$ and $a$. By the dual of \cite[Proposition~2.12]{Bu}, we obtain a commutative diagram 
    \begin{center}
        \begin{tikzcd}
        X \arrow[d, Rightarrow, no head] \arrow[r, "h'", tail] & P \arrow[d, "d"'] \arrow[dr, phantom, "\lrcorner", very near start] \arrow[r, "i", two heads] & A\arrow[d, "a"]\\
        X \arrow[r, "f", tail]      & Y\arrow[r, "g", two heads] & Z
        \end{tikzcd}
    \end{center}
    in which all rows are conflations, and the sequence
    \begin{center}
        \begin{tikzcd}[ampersand replacement=\&, column sep = 2em, row sep=2em]
        P  \arrow[r, "{\renewcommand{\arraystretch}{0.5}\begin{bmatrix}i\\-d\end{bmatrix}}", tail] \& A\oplus Y  \arrow[r, "{\setlength\arraycolsep{2pt}\begin{bmatrix}a&g\end{bmatrix}}", two heads] \& Z
        \end{tikzcd}
    \end{center}
    is a conflation. Since $Y\in\C'$, the morphism $d$ is admissible. It follows that there exist a deflation $b:P\twoheadrightarrow B$ and an inflation $c:B\rightarrowtail Y$ such that $d=c\circ b$. Let $j=\ker(b):L\rightarrowtail P$ and $e=\coker(c):Y\twoheadrightarrow D$. See the following diagram.
    \begin{center}
        \begin{tikzcd}[row sep=scriptsize]
        &L\arrow[d, "j", tail]&&\\
        X \arrow[dd, Rightarrow, no head] \arrow[r, "h'", tail] & P \arrow[d, "b"', two heads] \arrow[ddr, phantom, "\lrcorner", very near start] \arrow[r, "i", two heads] & A \arrow[dr, "m", dashed]\arrow[dd, "a"] & \\
        & B \arrow[d, "c", tail]& & K\arrow [dl, "k", tail, dashed]\\
        X \arrow[r, "f", tail]      & Y \arrow[r, "g", two heads] \arrow[d, "e", two heads]     & Z \arrow[dl, "l", dashed, two heads] &\\
        &D&&
        \end{tikzcd}
    \end{center}
    Then $e\circ f=e\circ c\circ b\circ h'=0$ implies the existence of a morphism $l:Z\to D$ with 
    $l\circ g=e.$ 
    Since $\C$ is weakly idempotent complete, $l$ is a deflation. Let $k=\ker(l): K\rightarrowtail Z$. Since $l\circ a\circ i=l\circ g\circ c\circ b=e\circ c\circ b=0$, we have that 
    $l\circ a=0.$ 
    This yields a morphism $m:A\to K$ with $k\circ m=a$. Consider the following commutative diagram
    \begin{center}
        \begin{tikzcd}[ampersand replacement=\&, column sep = 4em, row sep=2em] 
        L  \arrow[d, "j"', tail] \arrow[r, "i\circ j"]\& A \arrow[r, "m"] \arrow[d, "{\renewcommand{\arraystretch}{0.5}\begin{bmatrix}1\\0\end{bmatrix}}", tail] \& K  \arrow[d, "k", tail]  \\
        P  \arrow[d, "b"', two heads]\arrow[r, "{\renewcommand{\arraystretch}{0.5}\begin{bmatrix}i\\-d\end{bmatrix}}", tail] \& A\oplus Y  \arrow[d, "{\setlength\arraycolsep{2pt}\begin{bmatrix}0&1\end{bmatrix}}", two heads]\arrow[r, "{\setlength\arraycolsep{2pt}\begin{bmatrix}a&g\end{bmatrix}}", two heads] \& Z  \arrow[d, "l", two heads]    \\
        B \arrow[r, "-c"', tail] \& Y \arrow[r, "e"', two heads] \& D
        \end{tikzcd}
    \end{center}
    in which the second and third rows and all columns are conflations. Applying \cite[Corollary~3.6]{Bu}, it follows that the first row is also a conflation, which implies that $m$ is a deflation. Hence $a=k\circ m$ is admissible, and therefore $Z\in \C'$.
    
    \textbf{(iv) $\C'$ is abelian and the inclusion functor $\C' \to \C$ is exact.} By \cite[Proposition~3.1]{E}, an exact category in which every morphism is admissible is abelian, and its conflations coincide with short exact sequences. Since $\C'$ is closed under extensions in $\C$, it inherits an exact structure from $\C$. It therefore suffices to show that every morphism in $\C'$ is admissible in this induced exact structure. 
    
    Let $a: X'\to Y'$ be a morphism in $\C'$. Since $a$ is admissible in $\C$, there exist a deflation $p:X'\twoheadrightarrow M$ in $\C$ and an inflation $i:M\rightarrowtail Y'$ in $\C$ such that $a=i\circ p$. The associated conflations
    \[K\rightarrowtail X'\Twoheadrightarrow{p}M\text{ and }M\Rightarrowtail{i}Y'\twoheadrightarrow C\]
    have all their terms in $\C'$ by the closure under subobjects and quotients. Hence $p$ and $i$ are a deflation and an inflation in $\C'$, respectively. This shows that the morphism $a$ is admissible in $\C'$.
\end{proof}

In what follows, we extend the correspondence established in \cite{R} in a different direction, namely to the extriangulated category $\D^{[-d+1,0]}$. We begin by introducing the notion of wide subcategories in $\D^{[-d+1,0]}$.

\begin{definition}\label{wide}
    A subcategory $\W\subseteq\D^{[-d+1,0]}$ is said to be \emph{wide} if it satisfies the following conditions:
    \begin{enumerate}
        \item[(W1)] $\W$ is closed under extensions;
        \item[(W2)] $\Hom(\W,\W[i])=0$ for all $i<0$;
        \item[(W3)] the extriangulated structure of $\D^{[-d+1,0]}$ restricted to $\W$ is abelian.
    \end{enumerate}
    Furthermore, a wide subcategory is called \emph{length} if it is a length abelian category.
\end{definition}

\begin{remark}\label{rmk:wide}
    Condition (W1) implies that $\W$ inherits an extriangulated structure from $\D^{[-d+1,0]}$; see \cite[Remark~2.18]{NP}. By \cite[Corollary~2.7]{AET}, condition (W2) then ensures that this inherited extriangulated structure is, in fact, an exact structure. Consequently, condition (W3) means that $\W$ is an abelian category such that, for any $X,Y,Z\in \W$, a sequence
    \[X\to Y\to Z\] 
    is a short exact sequence if and only if it is a conflation in the exact category $\W$. Equivalently, $\W$ is a wide subcategory of itself in the sense of Definition~\ref{wide-exact}.
\end{remark}

\begin{remark}\label{rmk:summand}
    Since any abelian category is idempotent complete, it follows that every wide subcategory of $\D^{[-d+1,0]}$ is closed under direct summands.
\end{remark}

In the case where $d=1$, $\D^{[-d+1,0]}=\H$ is an abelian category, and our definition of wide subcategories coincides with the classical one. As mentioned above, the collection of simple objects in a wide subcategory of an abelian category forms a semibrick. Motivated by this observation, we introduce the following notion of semibricks in $\D^{[-d+1,0]}$, ensuring that the collection of simple objects in a wide subcategory remains a semibrick in this generalized setting.

\begin{definition}\label{semibrick}
    A collection of objects $\SS$ in $\D^{[-d+1,0]}$ is called a \emph{semibrick} if it satisfies the following conditions:
    \begin{enumerate}
        \item[(S1)] $\End(S)$ is a division algebra for all $S\in\SS$, that is, every $S$ is a brick.
        \item[(S2)] $\Hom(\SS,\SS[i])=0$ for all $i<0$.
        \item[(S3)] $\Hom(S,S')=0$ for any distinct $S,S'\in \SS$.
    \end{enumerate}
\end{definition}

Note that when $d=1$, the collection $\SS$ lies in the heart $\H$, and thus condition (S2) is automatically satisfied. In this case, Definition~\ref{semibrick} recovers the classical notion of a semibrick in an abelian category \cite{R}.

We consider semibricks up to isomorphism. Where appropriate, we identify a semibrick $\SS$ with the full subcategory consisting of objects in $\SS\cup\{0\}$.

For any wide subcategory $\W$ of $\D^{[-d+1,0]}$, since $\W$ is an abelian category, the collection $\ssim \W$ of its simple objects forms a semibrick. Thus, taking simple objects defines a map
\begin{equation}\label{eq:mapssim}
    \ssim\colon\{\text{wide subcategories of $\D^{[-d+1,0]}$}\} \to \{\text{semibricks in $\D^{[-d+1,0]}$}\}.
\end{equation}
Conversely, for any semibrick $\SS$ in $\D^{[-d+1,0]}$, let $W(\SS)$ denote its extension closure in $\D^{[-d+1,0]}$, that is, the smallest full subcategory of $\D^{[-d+1,0]}$ containing $\SS$ and closed under extensions. Recall that when $d=1$, the extension closure $W(\SS)$ coincides with $\Filt(\SS)$, which is well known to be a wide subcategory of $\H$; see \cite{R}.

\begin{proposition}\label{sb-wide}
    There is a bijection between the set of semibricks in $\D^{[-d+1,0]}$ and the set of length wide subcategories in $\D^{[-d+1,0]}$, given by the assignment 
    $\SS\mapsto W(\SS),$ 
    with inverse map given by 
    $\W\mapsto\ssim \W.$
\end{proposition}

\begin{proof}
    Let $\SS$ be a semibrick in $\D^{[-d+1,0]}$. By definition, $W(\SS)$ satisfies condition (W1). Since $\Hom(\SS,\SS[i])=0$ for all $i<0$, it follows by induction on the filtration length that $\Hom(W(\SS),W(\SS)[i])=0$ for all $i<0$; thus, $W(\SS)$ satisfies condition (W2). Consequently, $W(\SS)$ inherits an exact structure where the conflations are precisely those extriangles in $\D^{[-d+1,0]}$ whose terms all lie in $W(\SS)$. 
    
    In this exact setting, the extension closure $W(\SS)$ coincides with $\Filt(\SS)$. By \cite[Theorem~2.5]{E}, $\Filt(\SS)$ is a length wide subcategory of itself, with its simple objects given by $\SS$. Therefore, $W(\SS)$ is a length wide subcategory of $\D^{[-d+1,0]}$, and the equality $\ssim W(\SS) = \SS$ holds.

    Conversely, let $\W$ be a length wide subcategory in $\D^{[-d+1,0]}$. By Remark~\ref{rmk:wide}, $\W$ is a wide subcategory of itself. Hence, by \cite[Theorem~2.5]{E}, we have $\W = \Filt(\ssim \W)$. Since $\Filt(\ssim \W) = W(\ssim \W)$ by the argument above, this completes the proof.
\end{proof}

In the case where $d=1$, the bijection established above recovers the classical correspondence described in \cite[\S~1.2]{R}.

\begin{example}\label{exA2}
Let $\La$ be the path algebra $K(1\to 2)$ of a quiver of type $A_2$, and $d=2$. The Auslander-Reiten quiver of $\emod{2}{\La}$ is given as follows:
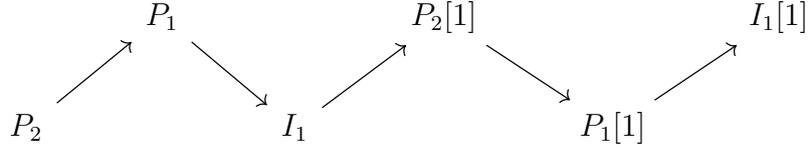
\begin{center}
    \begin{tikzcd}
    & P_1 \arrow[rd] &                & {P_2[1]} \arrow[rd] &                     & {I_1[1]} \\
    P_2 \arrow[ru] &                & I_1 \arrow[ru] &                     & {P_1[1]} \arrow[ru] &
    \end{tikzcd}
\end{center}

The poset of semibricks and their corresponding wide subcategories in $\emod{2}{\La}$ is depicted in Figure~\ref{figA2}. (The poset structure is inherited from the bijection with functorially finite, positive torsion classes; see Theorem~\ref{main}.) In the figure, semibricks are indicated by circles, while the corresponding wide subcategories are given by the additive hull of indecomposable objects in the shaded regions. In this specific example, all wide subcategories shown are either semisimple or equivalent to $\mod\La$, and hence, are abelian categories. Moreover, these constitute all possible wide subcategories; that is, every wide subcategory in this case is length. However, this property does not hold in general; see Example~\ref{exK2} for a counterexample.

\begin{figure}[htpb]
    \centering
    \begin{tikzpicture}
    \node (1) {\scalebox{0.5}{\usebox{\myboxa}}};
    \node (2)[above right = 0.5cm and 0.5cm of 1] {\scalebox{0.5}{\usebox{\myboxb}}};
    \node (3)[above left = 0.5cm and 0.5cm of 1] {\scalebox{0.5}{\usebox{\myboxc}}};
    \node (4)[above = 0.5cm of 2] {\scalebox{0.5}{\usebox{\myboxd}}};
    \node (5)[above = 0.5cm of 3] {\scalebox{0.5}{\usebox{\myboxe}}};
    \node (6)[above left = 0.5cm and 0.5cm of 3] {\scalebox{0.5}{\usebox{\myboxf}}};
    \node (7)[above = 0.5cm of 5] {\scalebox{0.5}{\usebox{\myboxg}}};
    \node (8)[above = 2cm of 4] {\scalebox{0.5}{\usebox{\myboxh}}};
    \node (9)[above = 0.5cm of 7] {\scalebox{0.5}{\usebox{\myboxi}}};
    \node (10)[above = 2.1cm of 6] {\scalebox{0.5}{\usebox{\myboxj}}};
    \node (11)[above = 0.5cm of 9] {\scalebox{0.5}{\usebox{\myboxk}}};
    \node (12)[above = 8.5cm of 1] {\scalebox{0.5}{\usebox{\myboxl}}};
    \draw  [->] (1.north east)-- (2.south west);
    \draw   [->] (1.north west)-- (3.south east);
   \draw    [->] (2.north)-- (4.south);
    \draw    [->] (3.north west)-- (6.south east); 
    \draw  [->] (3.north)-- (5.south);
    \draw   [->] (5.north west)-- (10.south);
   \draw    [->] (6.north)-- (9.south west); 
   \draw    [->] (5.north)-- (7.south); 
    \draw    [->] (7.north)-- (9.south);
     \draw    [->] (10.north east)-- (11.south west); 
      \draw    [->] (9.north)-- (11.south); 
       \draw    [->] (2.north west)-- (7.south east); 
        \draw    [->] (7.east)-- (8.south west);
         \draw    [->] (4.north)-- (8.south); 
          \draw    [->] (11.north east)-- (12.south west); 
           \draw    [->] (8.north)-- (12.south east); 
    \end{tikzpicture}
    \caption{Poset of semibricks and length wide subcategories in $\emod{2}{KA_2}$}
    \label{figA2}
\end{figure}

Consider the semibrick $\SS=\{P_2, I_1[1]\}$, where $W(\SS)=\add \SS$ is a semisimple abelian category. When $W(\SS)$ is viewed as a subcategory of $\D^b(\mod\La)$, the space $\Ext^2(I_1[1],P_2)$ is non-zero. However, when $W(\SS)$ is considered as an abelian category in its own right, the internal extension group satisfies $\Ext^2_{\mathrm{Ab}}(I_1[1],P_2)= 0$. This illustrates that, in general, higher extensions within $W(\SS)$ may not coincide with those in the ambient derived category. Note that our definition only requires an equivalence of extriangulated categories, which ensures an isomorphism on $\Ext^1$ but does not necessarily control higher-degree extensions.
\end{example}

While the correspondence between semibricks and length wide subcategories provides a classification of the simple building blocks of $\D^{[-d+1,0]}$, it is natural to ask how these wide subcategories fit into larger structural frameworks, such as those determined by finiteness or closure conditions. In the next section, we shift our focus to $d$-FAE closed subcategories and investigate their role in constructing wide subcategories, thereby further enriching the structural theory of the $d$-extended heart.

\section{From \texorpdfstring{$d$}{d}-FAE closed subcategories to wide subcategories}\label{S5}

In this section, we describe a constructive method for obtaining wide subcategories from $d$-FAE closed subcategories. To this end, we first introduce the notion of an ``exact heart'' for such subcategories, which serves as a generalization of the classical hearts of $t$-structures.

Let $\T$ be a $d$-FAE closed subcategory of $\D^{[-d+1,0]}$. We associate with $\T$ a full subcategory $\H_\T$ of $\D$ defined by:
\[\H_\T:=(\D^{\leq -d}*\T)\cap(\T^{\perp}[1]*\D^{\geq 0})\subseteq \D^{[-d,0]}.\]
We refer to $\H_\T$ as the \emph{exact heart} of $\T$.

\begin{lemma}\label{lem:H_T}
    Let $\T$ be a $d$-FAE closed subcategory of $\D^{[-d+1,0]}$. Then the following properties hold:
    \begin{enumerate}
        \item[(1)] $\Hom(\D^{\leq -d}*\T,(\T^{\perp}[1]*\D^{\geq 0})[j])=0$ for all $j<0$.
        \item[(2)] $\Hom(\H_\T,\H_\T[j])=0$ for all $j<0$.
        \item[(3)] $\H_\T$ is closed under extensions and direct summands.
        \item[(4)] $\H_\T$ is an idempotent complete exact category.
    \end{enumerate}
\end{lemma}

\begin{proof}
    For (1), since $\T^{\perp}[1]*\D^{\geq 0}\subseteq\D^{\geq -d}$, we have $\Hom(\D^{\leq -d},(\T^{\perp}[1]*\D^{\geq 0})[j])=0$ for all $j<0$. Thus, it suffices to show that
    \[\Hom(\T,\left(\T^\perp[1]*\D^{\geq 0}\right)[j])=0,\quad \text{for all } j < 0.\]
    By Proposition~\ref{prop:bi1}, the subcategory $\T^\perp[1]*\D^{\geq 0}$ is closed under negative shifts. Therefore, it is enough to verify the case $j=-1$. For any $X\in \T^\perp[1]*\D^{\geq 0}$, there exists a triangle
    \[X_1[1]\to X\to X_2,\]
    with $X_1\in \T^\perp$ and $X_2\in \D^{\geq 0}$. Applying the functor $\Hom(\T,-[-1])$ to this triangle, we obtain the exact sequence
    \[\Hom(\T,X_1)\to\Hom(\T,X[-1])\to\Hom(\T,X_2[-1]),\]
    where the first term vanishes because $X_1\in \T^\perp$, and the last term vanishes since $\T\subseteq\D^{[-d+1,0]}$ and $X_2[-1] \in \D^{\geq 1}$. This implies $\Hom(\T, X[-1]) = 0$, completing the proof of (1). Property (2) then follows directly from (1).
    
    For (3), Remark~\ref{fac->torf} ensures that $\T^{\perp}$ is a positive torsion-free class and is thus $d$-SAE closed. By Proposition~\ref{prop:bi1}, $\D^{\leq -d} * \T$ is a suspended subcategory and $\T^{\perp}[1] * \D^{\geq 0}$ is a cosuspended subcategory of $\D$. Consequently, both are closed under extensions. Furthermore, Remark~\ref{rmk:smd} implies they are closed under direct summands. It follows that their intersection $\H_\T$ is also closed under extensions and direct summands. 
    
    For (4), by combining the extension closure from (3) with the vanishing of negative extensions established in (2), by an argument analogous to the one in Remark~\ref{rmk:wide}, $\H_\T$ inherits an exact structure from the triangulated structure of $\D$. The idempotent completeness follows from the closure under direct summands established in (3).
\end{proof}

\begin{remark}\label{rmk:heart}
    If $\T$ is an $s$-torsion class, its exact heart $\H_\T$ coincides with the heart of the bounded $t$-structure associated to the $s$-torsion pair $(\T,\T^\perp)$; see Proposition~\ref{prop:bi2}. In this case, $\H_\T$ is an abelian category.
\end{remark}

Since $\T$ is closed under extensions and direct summands, the intersection $\T\cap\H_\T$ is a full subcategory of the exact category $\H_\T$ that remains closed under extensions and direct summands. It follows that $\T \cap \H_\T$ is itself an idempotent complete exact category.

\begin{definition}\label{defn:FAEtoWide}
    Let $\T\subseteq\D^{[-d+1,0]}$ be a $d$-FAE closed subcategory. We define a full subcategory $\LW(\T)$ of $\D^{[-d+1,0]}$ as follows:
    \[\LW(\T):=\{X\in \T\cap\H_\T\mid \text{every morphism in $\T\cap\H_\T$ to $X$ is admissible}\}.\]
\end{definition}

Applying Proposition~\ref{wide-in-exact} to the exact category $\T \cap \H_\T$ yields the main result of this section.

\begin{proposition}\label{prop:FAEtoWide}
    Let $\T \subseteq \D^{[-d+1,0]}$ be a $d$-FAE closed subcategory. Then $\LW(\T)$ is a wide subcategory of $\D^{[-d+1,0]}$.
\end{proposition}

\begin{proof}
    Since $\LW(\T)\subseteq \T\cap\H_\T\subseteq\H_\T$, Lemma~\ref{lem:H_T}(2) ensures that
    \[\Hom(\LW(\T), \LW(\T)[i])=0 \quad \text{for all } i < 0.\]
    As $\T \cap \H_\T$ is an idempotent complete exact category, Proposition~\ref{wide-in-exact} implies that $\LW(\T)$ is a wide subcategory of $\T \cap \H_\T$. In particular, $\LW(\T)$ is closed under extensions in $\T \cap \H_\T$, and consequently in $\D^{[-d+1,0]}$. Furthermore, since $\LW(\T)$ is a wide subcategory of itself (viewed as an exact category), it satisfies the requirements of Definition~\ref{wide}. Therefore, $\LW(\T)$ is a wide subcategory of $\D^{[-d+1,0]}$.
\end{proof}

Consequently, the assignment $\T \mapsto \LW(\T)$ defines a map
\begin{equation}\label{eq:mapLW}
    \LW \colon \left\{\begin{matrix}\text{$d$-FAE closed}\\\text{subcategories of $\D^{[-d+1,0]}$}\end{matrix}\right\}\to \left\{\begin{matrix}\text{wide subcategories}\\\text{of $\D^{[-d+1,0]}$}\end{matrix}\right\}.
\end{equation}

\begin{remark}\label{rmk:mapLW}
    In the case where $d=1$, $\D^{[-d+1,0]}=\H$ is an abelian category, and $\T$ is a full subcategory of $\H$ closed under extensions and factors. Since $\T\subseteq\H_\T$, we have $\T\cap\H_\T=\T$. In \cite[\S~2.3]{IT} (see also \cite[Exercise 8.23]{KS} as mentioned in \cite[\S~3]{MS}), the authors associated with $\T$ a wide subcategory $\alpha(\T)$ of $\H$ defined by:
    \[\alpha(\T):=\{X\in \T \mid \forall \ (g:Y\to X) \text{ with } Y \in \T,\ \ker(g)\in \T\}.\] 
    We claim that, in this case, $\LW(\T)=\alpha(\T)$, which demonstrates that the assignment $\LW$ generalizes the map $\alpha$.
    
    To see this, let $X\in \LW(\T)$ and let $g: Y \to X$ be a morphism in $\T$. Since $g$ is admissible in the exact category $\T$, it admits a factorization $g=g_2\circ g_1$ corresponding to conflations in $\T$:
    \[X_1\rightarrowtail Y \Twoheadrightarrow{g_1} X_2 \quad \text{and} \quad X_2\Rightarrowtail{g_2}X\twoheadrightarrow X_3,\]
    with $X_1, X_2, X_3\in \T$. Since $\T\subseteq \H$, these conflations are short exact sequences in $\H$. Consequently, $\ker(g)\cong\ker(g_1)= X_1\in \T$, which proves that $X \in \alpha(\T)$.
    
    Conversely, let $X\in \alpha(\T)$ and let $g \colon Y \to X$ be a morphism in $\T$. In the abelian category $\H$, we have the standard short exact sequences:
    \[ 0 \to \ker(g) \to Y \to \im(g) \to 0 \quad \text{and} \quad 0 \to \im(g) \to X \to \coker(g) \to 0. \]
    Since $\T$ is closed under factors, it follows that $\im(g),\coker(g)\in \T$. Moreover, the assumption $X \in \alpha(\T)$ ensures that $\ker(g) \in \T$. Thus, both sequences are conflations in $\T$. This implies that $g$ is an admissible morphism in $\T$, and therefore $X \in \LW(\T)$.
\end{remark}

\begin{remark}\label{W'=TcapH}
    The full subcategory $\T \cap \H_\T$ admits the following alternative characterization:
    \[\T\cap\H_\T=\{Y\in \T\mid \Hom(\T,Y[-1])=0\},\]
    which depends only on the negative extensions on $\T$, and does not require calculating the heart $\H_\T$. 
    
    The inclusion $\subseteq$ follows immediately from Lemma~\ref{lem:H_T}(1), which states that the morphism space $\Hom(\T, \H_\T[j]) = 0$ for all $j < 0$. To verify the converse inclusion, let $Y \in \T$ such that $\Hom(\T, Y[-1]) = 0$. Consider the truncation triangle
    \[ \sigma_{\geq 1}(Y[-1])[-1] \to \sigma_{\leq 0}(Y[-1]) \to Y[-1] \to \sigma_{\geq 1}(Y[-1]). \]
    Applying $\Hom(\T, -)$ to this triangle and noting that $\Hom(\T, \sigma_{\geq 1}(Y[-1])[-1]) = 0$ (as $\sigma_{\geq 1}(Y[-1])[-1] \in \D^{\geq 2}$ and $\T \subseteq \D^{\leq 0}$), we obtain $\Hom(\T, \sigma_{\leq 0}(Y[-1])) = 0$. This implies $\sigma_{\leq 0}(Y[-1]) \in \T^\perp$, and thus $Y \in \T^\perp[1] * \D^{\geq 0}$. Since $Y$ is also in $\T \subseteq \D^{\leq -d} * \T$, we conclude $Y \in \T \cap \H_\T$.
\end{remark}

\begin{example}\label{exK2}
    We continue to study the 2-extended module category $\emod{2}{\La}$ considered in Example~\ref{d-FAE-not-positive}, using the notation
    $\mathcal P$, $\mathcal R$, and $\mathcal I$ introduced there. Then the additive hull of $\mathcal{R}\cup\mathcal{I}\cup\mathcal{P}[1]$ in $\emod{2}{\La}$ is a wide subcategory that is not length.

    Now, take an arbitrary $M\in \mathcal{I}$ and consider the positive torsion class 
    \[\T:=\overline{M[1]},\]
    defined as the additive hull of $M[1]$ and all objects to its right in the Auslander-Reiten quiver of $2$-$\mod\La$ (as indicated by the shaded region in Figure~\ref{3}). Since $\overline{M[1]} \subseteq (\mod \La)[1]$, it follows from Remark~\ref{W'=TcapH} that  $\T\cap\H_{\T}=\overline{M[1]}.$ 

We claim that $\LW(\T) = \add M[1]$. By Remark~\ref{rmk:summand}, $\LW(\T)$ is closed under direct summands; thus, it suffices to show that $M[1]$ is the only indecomposable object in $\LW(\T)$. For any indecomposable object $N\in\overline{M[1]}$ not isomorphic to $M[1]$, say 
\[
(\begin{tikzcd} K^{r+1} \arrow[r, shift left] \arrow[r, shift right] & K^r \end{tikzcd})[1],
\]
there exists an irreducible surjection in $\mod\La$
\[f:(\begin{tikzcd}
K^{r+2} \arrow[r, yshift=0.5mm] \arrow[r, yshift=-0.5mm] & K^{r+1}
\end{tikzcd})\to (\begin{tikzcd}
K^{r+1} \arrow[r, yshift=0.5mm] \arrow[r, yshift=-0.5mm] & K^r
\end{tikzcd})\]
Then $f[1]$ is a morphism in $\overline{M[1]}$. If $f[1]$ were admissible in the exact category $\T \cap \H_{\T}$, it would admit a factorization $f[1] = h \circ g$, where $g$ is a deflation and $h$ is an inflation. Since $f[1]$ is irreducible, either $g$ must be a section (making $g$ an isomorphism) or $h$ must be a retraction (making $h$ an isomorphism). In the first case, $f[1]$ would be an inflation, implying its cone $C(f[1]) \in \mathcal{R}[2]$, which is not contained in $\overline{M[1]}$. In the second case, $f[1]$ would be a deflation, implying $C(f[1])[-1] \in \mathcal{R}[1]$, which is again not contained in $\overline{M[1]}$. Hence, $N \notin \LW(\T)$, and we conclude that $\LW(\T) = \add M[1]$.

This calculation provides a concrete realization of the general theory; in particular, the result aligns with the bijection established in Theorem~\ref{main2}, identifying $\overline{M[1]}$ as the positive torsion class $T(\{M[1]\})$.

\begin{figure}[htpb]
    \centering
     \begin{tikzcd}[column sep=tiny, execute at end picture={\begin{scope}[on background layer]\draw[rounded corners,fill=gray!20, draw=none] (A.south west)--(B.north west)--(C.north east)--(D.south east)--cycle;\end{scope}
     \draw[fill=none](A) circle (12.5 pt);}]
    & P_1 \arrow[rdd, shift left] \arrow[rdd] &        & \vdots \arrow[d, bend right]&      & \circ \arrow[rdd] \arrow[rdd, shift left] &        & \vdots \arrow[d, bend right]                      &  & |[alias=B]|\circ \arrow[rdd, shift left] \arrow[rdd] &    &     &                                            & |[alias=C]| {I_1[1]} \\ & & \cdots & \circ \arrow[d, bend right] \arrow[u, bend right] & \cdots & & \cdots & \circ \arrow[d, bend right] \arrow[u, bend right] & \cdots  &  & & \cdots \arrow[rd, shift left] \arrow[rd] &     &          \\
P_2 \arrow[ruu, shift right] \arrow[ruu] & & \circ  & \circ \arrow[u, bend right] \arrow[d, bend right] & \circ \arrow[ruu, shift left] \arrow[ruu] &    & \circ  & \circ \arrow[d, bend right] \arrow[u, bend right] & |[alias=A]|{M[1]} \arrow[ruu] \arrow[ruu, shift left] &  & \circ \arrow[ru, no head] \arrow[ru, Rightarrow, no head] &   & |[alias=D]|\circ \arrow[ruu, shift right] \arrow[ruu] &          \\ & &        & \circ \arrow[u, bend right]   &  &    &        & \circ \arrow[u, bend right]                       &   &  &    &       &    &         
\end{tikzcd}
    \caption{An example of $\LW(\T)$ for a positive torsion class $\T$}
    \label{3}
\end{figure}
\end{example}

We now investigate the simple objects in the wide subcategory $\LW(\T)$, or equivalently, the composition of the maps $\LW$ from \eqref{eq:mapLW} and $\ssim$ from \eqref{eq:mapssim}. For any $d$-FAE closed subcategory $\T$ of $\D^{[-d+1,0]}$, we denote by $\ssim \H_\T$ the set of (isoclasses of) simple objects in the exact heart $\H_\T$.

\begin{lemma}\label{incl}
    Let $\T \subseteq \D^{[-d+1,0]}$ be a $d$-FAE closed subcategory. Then any simple object in the wide subcategory $\LW(\T)$ is also a simple object in the exact heart $\H_\T$. In particular, we have 
    \[\ssim \LW(\T)\subseteq\left(\ssim \H_\T\right)\cap\D^{[-d+1,0]}.\]
\end{lemma}
    
\begin{proof}
    Since $\LW(\T)\subseteq\T\cap\H_\T\subseteq\D^{[-d+1,0]}$, it suffices to prove the first statement. Let $S \in \ssim \LW(\T)$, and consider an arbitrary conflation
    \[X \Rightarrowtail{f} S \twoheadrightarrow Y\]
    in $\H_\T$. We must show that either $X=0$ or $Y=0$. Suppose $X\neq 0$. It follows from the above conflation that
    \[X\in \H_\T\cap(Y[-1]* S)\subseteq \D^{[-d,0]}\cap\D^{[-d+1,1]}=\D^{[-d+1,0]}.\]
    Moreover, as $X\in\H_\T\subseteq \D^{\leq -d}*\T$, we have $X\in\T$. Thus $X \in \T \cap \H_\T$. 
    
    Since $S \in \LW(\T)$, the morphism $f \colon X \to S$ is admissible in the exact category $\T \cap \H_\T$. Thus, there exist conflations 
    \[X'\rightarrowtail X\Twoheadrightarrow{f_1}X'' \quad \text{and} \quad X''\Rightarrowtail{f_2}S\twoheadrightarrow X'''\]
    in $\T\cap\H_\T$ such that $f=f_2\circ f_1$. By Proposition~\ref{wide-in-exact}, $\LW(\T)$ is closed under subobjects and factors in $\T\cap\H_\T$; therefore, both $X''$ and $X'''$ belong to $\LW(\T)$. Since $S$ is simple in $\LW(\T)$, we must have either $X''=0$ or $X'''=0$.
    
    The assumption $X \neq 0$ implies $f \neq 0$, which forces $f_2 \neq 0$ and $X'' \neq 0$. This implies $X''' = 0$, meaning the map $f_2 \colon X'' \rightarrowtail S$ is an isomorphism. Consequently, $f$ is equivalent to the deflation $f_1$ in $\T \cap \H_\T$, with kernel $X' \cong Y[-1]$. However, since $X' \in \T \cap \H_\T \subseteq \H_\T$ and $Y \in \H_\T$, the condition $X' \cong Y[-1]$ is only possible if $Y = 0$, as $\Hom(\H_\T, \H_\T[-1]) = 0$. This completes the proof.
\end{proof}

So far, we have constructed a map \eqref{eq:mapLW} from $d$-FAE closed subcategories to wide subcategories via the exact heart construction. Specifically, we have shown that the semibrick associated with such a wide subcategory is contained in the set of simple objects of the corresponding exact heart. To complete the structural correspondence, it is natural to investigate the inverse direction. In the following section, we shall study the construction that leads from semibricks back to the torsion-like classes introduced in Section~\ref{S3}, such as $d$-FAE closed subcategories and positive torsion classes. This development will culminate in the main classification theorems of this paper.

\section{From semibricks to torsion-like classes}\label{S6}

For module categories of finite-dimensional algebras, it was shown in \cite[Lemma~2.8]{A} that the assignment sending a semibrick $\SS$ in $\mod\La$ to the smallest torsion class containing $\SS$ is injective. Recall from Section~\ref{S3} that we introduced three notions of torsion-like classes in $\D^{[-d+1,0]}$. While these notions coincide with the classical concept of torsion classes in module categories, they are generally not equivalent for $d > 1$. Among these, $d$-FAE closed subcategories and positive torsion classes are closed under intersections. Consequently, we can define the following two maps:
\begin{equation}\label{eq:twomaps}
\xymatrix{
    & {}\{\text{positive torsion classes in }\D^{[-d+1,0]}\}\\
    \{\text{semibricks in }\D^{[-d+1,0]}\}\ar[ur]^{T\qquad}\ar[dr]_{\phi\qquad}\\
    & \left\{\text{$d$-FAE closed subcategories in $\D^{[-d+1,0]}$}\right\}}
\end{equation}
These two maps send a semibrick $\SS$ in $\D^{[-d+1,0]}$ to the smallest positive torsion class $T(\SS)$ containing $\SS$, and the smallest $d$-FAE closed subcategory $\phi(\SS)$ containing $\SS$, respectively. As in Section~\ref{S4}, we identify a semibrick $\SS$ with the full subcategory consisting of its objects and zero. Since every positive torsion class is $d$-FAE closed, the inclusion $\phi(\SS) \subseteq T(\SS)$ holds.

We remark that if $\D$ is $K$-linear, $\Hom$-finite, and Krull--Schmidt, and if the $d$-extended heart $\D^{[-d+1,0]}$ contains only finitely many indecomposable objects, then every subcategory of $\D^{[-d+1,0]}$ closed under direct sums and summands is functorially finite. In such cases, Proposition~\ref{s-tors} implies that $\phi(\SS) = T(\SS)$. However, whether this equality holds in general is not pursued in this paper.

In this section, we will prove that both $\phi$ and $T$ are injective. Furthermore, we will show that the simple objects in the wide subcategory $\LW(\phi(\SS))$, constructed from $\phi(\SS)$ as in Section~\ref{S5}, are precisely the elements of $\SS$. Although it remains unclear whether a similar characterization holds for $T(\SS)$ in general, we show that $T(\SS)$ satisfies this property whenever $T(\SS)$ is an $s$-torsion class (Proposition~\ref{idforphi}).

Let $\SS$ be a semibrick in $\D^{[-d+1,0]}$. We denote by $\operatorname{Tor}(\SS)$ the set of $d$-FAE closed subcategories $\T$ of $\D^{[-d+1,0]}$ that contain $\SS$ and satisfy $\T^\perp = \SS^{\perp_{\leq 0}}$.

\begin{lemma}\label{same-tor-free}
Both $T(\SS)$ and $\phi(\SS)$ belong to $\operatorname{Tor}(\SS)$. Moreover, for any $\T\in\operatorname{Tor}(\SS)$, the inclusions $\phi(\SS)\subseteq\T\subseteq T(\SS)$ hold.
\end{lemma}

\begin{proof}
    According to Proposition~\ref{perp:tor}, $(T(\SS), \SS^{\perp_{\leq 0}})$ is a positive torsion pair. This implies $T(\SS)^\perp = \SS^{\perp_{\leq 0}}$ and ${}^\perp(\SS^{\perp_{\leq 0}}) = T(\SS)$. The first equality ensures that $T(\SS) \in \operatorname{Tor}(\SS)$, while the second equality implies that any $\T \in \operatorname{Tor}(\SS)$, being a subcategory of ${}^\perp(\T^\perp) = {}^\perp(\SS^{\perp_{\leq 0}})$, must be contained in $T(\SS)$.
    
    Next, since $\phi(\SS) \subseteq T(\SS)$, we have $\phi(\SS)^\perp \supseteq T(\SS)^\perp = \SS^{\perp_{\leq 0}}$. Recall that $\phi(\SS)$ is $d$-FAE closed, so Remark~\ref{fac->torf} yields $\phi(\SS)^\perp = \phi(\SS)^{\perp_{\leq 0}}$. Furthermore, as $\SS \subseteq \phi(\SS)$, it follows that $\SS^{\perp_{\leq 0}} \supseteq \phi(\SS)^{\perp_{\leq 0}}$. Combining these inclusions, we obtain $\phi(\SS)^\perp = \SS^{\perp_{\leq 0}}$, which means $\phi(\SS) \in \operatorname{Tor}(\SS)$. Finally, for any $\T \in \operatorname{Tor}(\SS)$, since $\T$ is by definition a $d$-FAE closed subcategory containing $\SS$, it must contain the smallest such subcategory $\phi(\SS)$. This establishes the inclusion $\phi(\SS) \subseteq \T$.
\end{proof}

Recall from Section~\ref{S5} that, for a $d$-FAE closed subcategory $\T$ of $\D^{[-d+1,0]}$, the exact heart of $\T$ is defined as 
$\H_\T:=(\D^{\leq -d}*\T)\cap (\T^\perp[1]*\D^{\geq 0}).$ 
We denote the set of (isoclasses of) simple objects in $\H_\T$ by $\ssim \H_\T$. The following proposition guarantees that we can recover the original semibrick $\SS$ from the exact heart of any $\T \in \operatorname{Tor}(\SS)$.

\begin{proposition}\label{sim-in-heart}
    Let $\SS$ be a semibrick in $\D^{[-d+1,0]}$ and let $\T\in\operatorname{Tor}(\SS)$. Then \[\SS=(\ssim \H_{\T})\cap\D^{[-d+1,0]}.\]
\end{proposition}

\begin{proof}
    We divide the proof into three steps.
    
    \textbf{Step 1: $\boldsymbol{\SS\subseteq \H_{\T}}$.} Since $\SS\subseteq \T\subseteq \D^{\leq -d}*\T$, it suffices to show $\SS\subseteq \T^\perp[1]*\D^{\geq 0}$. Let $X\in\SS$. Applying the functor $\Hom(\SS,-)$ to the truncation triangle 
    \[\sigma_{\leq -1}(X)\to X\to \sigma_{\geq 0}(X),\]
    we obtain the exact sequence
    \[\Hom(\SS,(\sigma_{\geq 0}(X))[j])\to \Hom(\SS,(\sigma_{\leq -1}(X))[j+1])\to \Hom(\SS, X[j+1]).\] 
    When $j\leq -2$, the first term vanishes because $\SS\subseteq\D^{[-d+1,0]}$ and $\sigma_{\geq 0}(X)[j] \in \D^{\geq 2}$, while the last term vanishes since $j+1 \leq -1$ and $\SS$ is a semibrick. Consequently, the middle term also vanishes, which implies $(\sigma_{\leq -1}(X))[-1]\in \SS^{\perp_{\leq 0}}=\T^\perp$. Thus, $X\in \sigma_{\leq -1}(X) * \sigma_{\geq 0}(X)\subseteq T^\perp[1]*\D^{\geq 0}$, and we conclude $\SS\subseteq \H_{\T}$.
    
    \textbf{Step 2: $\boldsymbol{\SS\subseteq\ssim \H_{\T}}$.} Let $X\in \SS$ and let $L\Rightarrowtail{f} X\twoheadrightarrow M$ be a conflation in $\H_{\T}$. We show that either $L=0$ or $M=0$. Suppose $L\neq 0$. Since $L\in \H_{\T}\subseteq\D^{\leq -d}*\T$, there exists a triangle $L_1\xrightarrow{g} L\xrightarrow{h} L_2$ with $L_1\in \D^{\leq -d}$ and $L_2\in \T$. Since $X\in \T$ and $\Hom(\D^{\leq -d},\T)=0$, we have $f\circ g=0$. Thus, $f$ factors through $h$, yielding the following commutative diagram of triangles via the octahedral axiom:
    \begin{center}
    \begin{tikzcd}
        L_1 \arrow[d, "g'"'] \arrow[r, Rightarrow, no head] & L_1 \arrow[d, "g"] & \\
        M[-1] \arrow[d] \arrow[r] & L \arrow[d, "h"] \arrow[r, "f"] & X \arrow[d, Rightarrow, no head] \\
        W \arrow[r] & L_2 \arrow[r] & X                               
    \end{tikzcd}  
    \end{center}
    Given $M[-1]\in \D^{[-d+1,1]}$ and $L_1 \in \D^{\leq -d}$, we have $g'=0$, which implies $g=0$. Consequently, $L_2\cong L_1[1]\oplus L$, which forces $L_1=0$ and $L\cong L_2\in \T$. Since $L\neq 0$, it follows that $L\notin\T^\perp=\SS^{\perp_{\leq 0}}$. Therefore, there exists a non-zero morphism $f':X'\to L[m]$ for some $X'\in\SS$ and $m\leq 0$.
    
    Since both $X'$ and $L$ belong to $\H_{\T}$, Lemma~\ref{lem:H_T}(2) implies $\Hom(\H_{\T}, \H_{\T}[m]) = 0$ for $m < 0$, forcing $m = 0$. Moreover, $f \circ f'$ is non-zero; otherwise, by the triangle in the second row of the diagram above, $f'$ would factor through $M[-1] \in \H_{\T}[-1]$, contradicting $\Hom(\H_{\T}, \H_{\T}[-1]) = 0$ from Lemma~\ref{lem:H_T}(2). Since $\SS$ is a semibrick and $X', X \in \SS$, the non-zero morphism $f \circ f' \colon X' \to X$ must be an isomorphism. It follows that $f$ is a retraction in the exact category $\H_{\T}$. As $f$ is also an inflation (being part of a conflation), it must be an isomorphism, which implies $M = 0$. This proves $X \in \ssim \H_{\T}$.
    
    \textbf{Step 3: $\boldsymbol{\left(\ssim \H_{\T}\cap\D^{[-d+1,0]}\right)\setminus\SS=\emptyset}$.} Suppose there exists $X \in (\ssim \H_{\T} \cap \D^{[-d+1,0]}) \setminus \SS$. Then $X \in \H_\T \cap \D^{[-d+1,0]} \subseteq \T$. If there exists a non-zero $f \colon S \to X$ with $S \in \SS$, then $C(f) \in \H_\T$ (proof below). This makes $S \xrightarrow{f} X \to C(f)$ a conflation in $\H_\T$. As $X$ is simple, $C(f)=0$ and $f$ is an isomorphism, contradicting $X \notin \SS$. If no such $f$ exists, then $\Hom(\SS, X) = 0$. Since $\SS\subseteq\H_\T$ and $X\in\H_\T$, Lemma~\ref{lem:H_T}(2) implies that $\Hom(\SS, X[i]) = 0$ for all $i < 0$. Consequently, $X \in \SS^{\perp_{\leq 0}} = \T^\perp$. Since $X$ is also an element of $\T$, we must have $X = 0$, which contradicts the fact that $X$ is a simple object in $\H_\T$.
 
    Finally, we prove $C(f) \in \H_\T$. Lemma~\ref{summand}(2) ensures that $\sigma_{\geq -d+1}(C(f)) \in \T$, which implies $C(f)\in\D^{\leq -d}*\T$. Let $S'\in \SS$. Applying the functor $\Hom(S',-)$ to the triangle $S\xrightarrow{f}X\to C(f)$, we obtain the following exact sequence for any $j \leq 0$:
    \[\Hom(S',X[j-1])\to \Hom(S',C(f)[j-1])\to \Hom(S',S[j])\to \Hom(S',X[j]).\]
    The first term vanishes since $S', X \in \H_\T$ and $j-1 < 0$ (applying Lemma~\ref{lem:H_T}(2)).
    \begin{itemize}
    \item For $j=0$: If $S' \neq S$, then $\Hom(S', S) = 0$ as $\SS$ is a semibrick. If $S' = S$, the last map $\Hom(S', f)$ in the sequence is injective because $f \neq 0$ and $S$ is a brick. In either case, we obtain $\Hom(S', C(f)[-1]) = 0$.
    \item For $j < 0$: The third term $\Hom(S', S[j])$ vanishes because $\SS$ is a semibrick. 
    \end{itemize}
    It follows that $\Hom(S', C(f)[k]) = 0$ for all $k \leq -1$. We then apply $\Hom(S', -)$ to the truncation triangle
    \[ \sigma_{\leq 0}(C(f)[-1]) \to C(f)[-1] \to \sigma_{\geq 1}(C(f)[-1]),\]
    yielding the exact sequence
    \[\Hom(S', \sigma_{\geq 1}(C(f)[-1])[j])\to \Hom(S',\sigma_{\leq 0}(C(f)[-1])[j+1])\to \Hom(S',C(f)[j]).\]
    For all $j \leq -1$, the first term vanishes because $\sigma_{\geq 1}(C(f)[-1]) \in \D^{\geq 1}$ and $S' \in \D^{\leq 0}$, while the last term vanishes by our previous calculation. This implies that $\sigma_{\leq 0}(C(f)[-1]) \in \SS^{\perp_{\leq 0}} = \T^\perp$. Consequently, $C(f)[-1] \in \T^\perp * \D^{\geq 1}$, or equivalently, $C(f) \in \T^\perp[1] * \D^{\geq 0}$. Combined with $C(f) \in \D^{\leq -d} * \T$, we conclude $C(f) \in \H_\T$.
\end{proof}

Applying Proposition~\ref{sim-in-heart} to the specific cases where $\T = T(\SS)$ and $\T = \phi(\SS)$, we obtain the following consequence.

\begin{corollary}\label{inj}
    The maps $T$ and $\phi$ defined in \eqref{eq:twomaps} are both injective.
\end{corollary}

\begin{proof}
    Suppose $\SS$ and $\SS'$ are two semibricks in $\D^{[-d+1,0]}$ such that $T(\SS) = T(\SS')$. According to Lemma~\ref{same-tor-free}, we have $T(\SS) \in \operatorname{Tor}(\SS)$ and $T(\SS') \in \operatorname{Tor}(\SS')$. Proposition~\ref{sim-in-heart} then implies that
    \[ \SS = (\ssim \H_{T(\SS)}) \cap \D^{[-d+1,0]} = (\ssim \H_{T(\SS')}) \cap \D^{[-d+1,0]} = \SS'. \]
    This establishes the injectivity of $T$. The injectivity of $\phi$ follows from an analogous argument.
\end{proof}

The following is another immediate consequence of Proposition~\ref{sim-in-heart}, which will be used later.

\begin{lemma}\label{no-neg-map}
Let $\SS$ be a semibrick in $\D^{[-d+1,0]}$ and let $\T\in\operatorname{Tor}(\SS)$. Then \[\Hom(\T,\SS[j])=0\text{ for all } j<0.\]
\end{lemma}

\begin{proof}
    By Proposition~\ref{sim-in-heart}, we have $\SS\subseteq\H_\T$. The required vanishing then follows directly from Lemma~\ref{lem:H_T}(1) and the fact that $\T \subseteq \D^{\leq -d} * \T$.
\end{proof}

Recall from Lemma~\ref{union} that there is an explicit description of the smallest $d$-FAE closed subcategory containing a given subcategory. From this, we can deduce the following property of $\phi(\SS)$. In the classical case, this recovers \cite[Lemma~2.7(1)]{A}.

\begin{lemma}\label{closed-cocone}
For any non-zero morphism $f \colon X\to S$ with $X\in \phi(\SS)$ and $S\in \SS$, the cocone $C(f)[-1]$ belongs to $\phi(\SS)$.
\end{lemma}

\begin{proof}
By Lemma~\ref{union}, $X$ belongs to $\Phi^s(\SS)$ for some $s \geq 0$, where $\Phi^0(\SS)=\SS$ and 
\[\Phi^i(\SS)=\Fac_{d}\left(\Phi^{i-1}(\SS)\right)*\Fac_{d}\left(\Phi^{i-1}(\SS)\right),\quad \text{for any $i>0$}.\]
We proceed by induction on $s$ to prove the result.

The base case $s=0$ implies $X \in \SS$. Since $\SS$ is a semibrick, the non-zero morphism $f$ is an isomorphism, hence $C(f)[-1] = 0 \in \phi(\SS)$. Now, assume the result holds for all $s \leq m$, and let $s = m+1$.

We first establish a claim: for any non-zero morphism $g \colon Z \to S$ with $Z \in \Fac_d(\Phi^{m}(\SS))$, the cocone $C(g)[-1]$ belongs to $\phi(\SS)$. By the definition of $d$-factors, there exist $d$ triangles
\[Z_i\xrightarrow{f_i} Y_i\xrightarrow{g_i} Z_{i-1}\xrightarrow{h_i} Z_i[1],\ 1\leq i\leq d,\]
with $Z_0=Z$, $Z_1,\dots,Z_d\in\D^{[-d+1,0]}$, and $Y_1,\dots,Y_d\in\Phi^m(\SS)\subseteq\phi(\SS)$. We claim that the composition $g\circ g_1: Y_1\to S$ is non-zero. Otherwise, $g$ would factor through $h_1$ as shown in the following diagram:
\begin{center}
\begin{tikzcd}
Z_1 \arrow[r, "f_1"] & Y_1 \arrow[r, "g_1"] & Z_0 \arrow[r, "h_1"] \arrow[d, "g"] & {Z_1[1]} \arrow[ld, "\alpha_1", dashed] \\ &  & S         & 
\end{tikzcd}
\end{center}
By Lemma~\ref{no-neg-map}, $\Hom(Y_2, S[-1]) = 0$, which implies that the composition $\alpha_1[-1] \circ g_2 \colon Y_2 \to S[-1]$ vanishes. Thus, $\alpha_1[-1]$ factors through $h_2$:
\begin{center}
\begin{tikzcd}
Z_2 \arrow[r, "f_2"] & Y_2 \arrow[r, "g_2"] & Z_1 \arrow[r, "h_2"] \arrow[d, "{\alpha_1[-1]}"'] & {Z_2[1]} \arrow[ld, "\alpha_2", dashed] \\&                      & {S[-1]}                                           &                        \end{tikzcd}
\end{center}
Continuing this process, we obtain a factorization
\[g=\alpha_d[d-1]\circ h_d[d-1]\circ\cdots\circ h_2[1]\circ h_1,\]
where $\alpha_d: Z_d[1]\to S[-d+1]$ has to be zero because both $Z_d$ and $S$ belong to $\D^{[-d+1,0]}$. This implies $g=0$, a contradiction. Thus, $g\circ g_1\neq 0$. By the octahedral axiom, we have the following commutative diagram:
\begin{center}
\begin{tikzcd}
Z_1 \arrow[d, Rightarrow, no head] \arrow[r] & V \arrow[r] \arrow[d]& C(g)[-1] \arrow[d]\\
Z_1 \arrow[r, "f_1"] & Y_1 \arrow[d, "g\circ g_1"] \arrow[r, "g_1"] & Z_0 \arrow[d, "g"] \\
& S \arrow[r, Rightarrow, no head]             & S                 
\end{tikzcd}
\end{center}
By the inductive hypothesis applied to $Y_1\in\Phi^m(\SS)$, the object $V=C(g\circ g_1)[-1]$ belongs to $\phi(\SS)\subseteq\D^{[-d+1,0]}$. It follows that
\[C(g)[-1]\in (V*Z_1[1])\cap (S[-1]*Z_0)\subseteq \D^{[-d+1,0]}.\]
The first row of the diagram then shows $C(g)[-1] \in \Fac_d(\phi(\SS)) \subseteq \phi(\SS)$, which proves the claim.

For the general case $X \in \Phi^{m+1}(\SS)$, there exists a triangle $X_1 \xrightarrow{\alpha} X \xrightarrow{\beta} X_2$ with $X_1, X_2 \in \Fac_d(\Phi^{m}(\SS))$. We distinguish two cases:
\begin{enumerate}
    \item[(1)] $f \circ \alpha \neq 0$: The octahedral axiom gives the following diagram:
    \begin{center}
        \begin{tikzcd}
        M \arrow[r] \arrow[d] & C(f)[-1] \arrow[r] \arrow[d] & X_2 \arrow[d, Rightarrow, no head] \\
        X_1 \arrow[r, "\alpha"] \arrow[d, "f\circ \alpha\neq 0"'] & X \arrow[d, "f"] \arrow[r, "\beta"] & X_2 \\
        S \arrow[r, Rightarrow, no head] & S & 
        \end{tikzcd}
    \end{center}
    By the claim, $M = C(f \circ \alpha)[-1] \in \phi(\SS)$. Since $X_2 \in \phi(\SS)$, we conclude $C(f)[-1] \in \phi(\SS)$.
    \item[(2)] $f \circ \alpha = 0$: Then $f$ factors through $\beta$ as $f = \gamma \circ \beta$ for some $\gamma \colon X_2 \to S$. The octahedral axiom yields:
    \begin{center}
    \begin{tikzcd}
    X_1 \arrow[r] \arrow[d, Rightarrow, no head] & C(f)[-1] \arrow[r] \arrow[d] & N \arrow[d]\\
    X_1 \arrow[r, "\alpha"]  & X \arrow[d, "f"] \arrow[r, "\beta"] & X_2 \arrow[d, "\gamma\neq 0"] \\
    & S \arrow[r, Rightarrow, no head]    & S                   \end{tikzcd}
    \end{center}
    By the claim, $N = C(\gamma)[-1] \in \phi(\SS)$. Since $X_1 \in \phi(\SS)$, it follows that $C(f)[-1] \in \phi(\SS)$.
\end{enumerate}
    This completes the induction and the proof.
\end{proof}

For $T(\SS)$, an analogous statement holds only under a certain additional condition and with a different approach.

\begin{lemma}\label{closed-cocone-pos}
    Let $\SS$ be a semibrick in $\D^{[-d+1,0]}$ such that $T(\SS)$ is an $s$-torsion class. For any non-zero morphism $f: Y\to X$ with $Y\in T(\SS)$ and $X\in \SS$, we have that $C(f)[-1]\in T(\SS)$.
\end{lemma}

\begin{proof}
    By Proposition~\ref{prop:bi2}~(1), $T(\SS)$ is an $s$-torsion-free class in the $d$-extended heart
    \[d\text{-}\H_{T(\SS)}=\H_{T(\SS)}[d-1]*\cdots*\H_{T(\SS)}[1]*\H_{T(\SS)}.\] 
    Thus, there exists a triangle 
    \[Y_1\to Y\to Y_2\to Y_1[1],\]
    with $Y_1\in \H_{T(\SS)}[d-1]*\cdots*\H_{T(\SS)}[1]$ and $Y_2\in \H_{T(\SS)}$. Since by Proposition~\ref{sim-in-heart}, $\SS\subseteq\ssim \H_{T(\SS)}\subseteq\H_{T(\SS)}$, the morphism space $\Hom(Y_1,S)$ vanishes. It follows that the morphism $f$ factors through $Y_2$. Using the octahedral axiom, we get the following commutative diagram of triangles. 
    \begin{center}
        \begin{tikzcd}
        Y_1 \arrow[r] \arrow[d, Rightarrow, no head] & C(f)[-1] \arrow[r] \arrow[d]            & N \arrow[d]   \\
        Y_1 \arrow[r]    & Y \arrow[r] \arrow[d, "f"]            & Y_2 \arrow[d,"g"] \\
        & X \arrow[r, Rightarrow, no head] & X            
        \end{tikzcd}
    \end{center}
    Since $f\neq 0$, so is $g$. Since $Y_2\in \H_{T(\SS)}$ and $\SS\subseteq\ssim \H_{T(\SS)}$, the morphism $g$ is an epimorphism in $\H_{T(\SS)}$, which, together with the triangle in the third column of the above diagram implies $N\in \H_{T(\SS)}$. Therefore, by the triangle in the first row, we have
    \[C(f)[-1]\in Y_1*N\subseteq \H_{T(\SS)}[d-1]*\cdots*\H_{T(\SS)}[1]*\H_{T(\SS)}=d\text{-}\H_{\T(\SS)}.\]
    Since $T(\SS)$ is an $s$-torsion-free class in $d\text{-}\H_{\T(\SS)}$, due to $X,Y\in T(\SS)$, by Lemma~\ref{positive}, we get that $C(f)[-1]\in T(\SS)$. 
\end{proof}

Recall from Definition~\ref{defn:FAEtoWide} and Proposition~\ref{prop:FAEtoWide} that, for any $d$-FAE closed subcategory $\T$ of $\D^{[-d+1,0]}$, we associate a wide subcategory $\LW(\T)$.

\begin{proposition}\label{idforphi}
    Let $\SS$ be a semibrick in $\D^{[-d+1,0]}$. Then $\SS=\ssim \LW(\phi(\SS))$. If $T(\SS)$ is an $s$-torsion class in $\D^{[-d+1,0]}$, then $\SS=\ssim \LW(T(\SS))$.
\end{proposition}

\begin{proof}
    Let $\T=\phi(\SS)$ or $T(\SS)$. By Lemma~\ref{same-tor-free}, $\T\in\operatorname{Tor}(\SS)$. Applying Lemma~\ref{incl} we have the inclusion:
    \[\ssim \LW(\T)\subseteq(\ssim \H_{\T})\cap\D^{[-d+1,0]}.\]
    Applying Proposition~\ref{sim-in-heart} yields:
    \begin{equation}\label{eq:S=sim}
        \SS=\left(\ssim \H_{\T}\right)\cap\D^{[-d+1,0]}.
    \end{equation}
    Consequently, it suffices to show that $\SS \subseteq \ssim \LW(\T)$.
    
    From \eqref{eq:S=sim}, we know that $\SS \subseteq \H_{\T}$, which implies $\SS \subseteq \T \cap \H_{\T}$. To establish $\SS \subseteq \LW(\T)$, we must show that any morphism $f \colon Y \to S$ in the exact category $\T \cap \H_{\T}$ with $S \in \SS$ is admissible. Since the zero morphism is always admissible, we assume $f \neq 0$. In this case, Lemmas~\ref{closed-cocone}~and~\ref{closed-cocone-pos} ensure that $C(f)[-1] \in \T$. Applying $\Hom(\T, -)$ to the triangle $Y \xrightarrow{f} S \to C(f)$, we obtain the exact sequence:
    \[\Hom(\T,S[-2])\to \Hom(\T,C(f)[-2])\to \Hom(\T,Y[-1]).\]
    The first and last terms vanish by Lemma~\ref{lem:H_T}(1). Hence, $\Hom(\T, C(f)[-2]) = 0$, and by the characterization in Remark~\ref{W'=TcapH}, we conclude that $C(f)[-1] \in \T \cap \H_{\T}$.
    
    Thus, the sequence $C(f)[-1]\to Y\xrightarrow{f}S$ is a conflation in the exact category $\T \cap \H_{\T}$. In particular, $f$ is a deflation and is therefore admissible. This proves $\SS \subseteq \LW(\T)$. Furthermore, since the elements of $\SS$ are simple objects in $\H_{\T}$ by \eqref{eq:S=sim}, they are necessarily simple in the subcategory $\LW(\T)$ as well.
\end{proof}

For a given $d$-FAE closed subcategory $\T$, the resulting wide subcategory $\LW(\T)$ may not be a length category in general. To establish a self-consistent framework within the set of length wide subcategories, which are in one-to-one correspondence with semibricks, we introduce the modified map
\[\LWlen = W \circ \ssim \circ \LW \colon \left\{\begin{matrix}\text{$d$-FAE closed}\\\text{subcategories of $\D^{[-d+1,0]}$}\end{matrix}\right\}\to \left\{\begin{matrix}\text{length wide}\\\text{subcategories of $\D^{[-d+1,0]}$}\end{matrix}\right\}.\]
This composition refines the wide subcategory $\LW(\T)$ by first extracting its simple objects and then taking their extension closure, thereby ensuring that the output is a length wide subcategory.

For any wide subcategory $\W$, similarly as for a semibrick, we denote by $\phi(\W)$ (resp. $T(\W)$) the smallest $d$-FAE closed  subcategory (resp. positive torsion class) of $\D^{[-d+1,0]}$ containing $\W$. By Propositions~\ref{sb-wide}~and~\ref{idforphi}, we have the following consequence, which in classical case recovers \cite[Proposition~3.3]{MS}. 

\begin{corollary}\label{cor:idforphi}
    Let $\W$ be a length wide subcategory of $\D^{[-d+1,0]}$. Then $\LWlen(\phi(\W))=\W$. If $T(\W)$ is an $s$-torsion class, then $\LWlen(T(\W))=\W$.
\end{corollary}

The results obtained up to this point are summarized in Figure~\ref{fig:outerbijs}, with two identities $\ssim\circ\LW\circ\phi=\id$ and $\LWlen\circ\phi=\id$. 

\begin{figure}[htpb]
    \centering
    \begin{tikzpicture}
        \draw(0,-3)node{$\left\{\begin{matrix}\text{positive torsion}\\\text{classes in $\D^{[-d+1,0]}$}\end{matrix}\right\}$};
        \draw(0,-5)node{$\left\{\begin{matrix}\text{$d$-FAE closed}\\\text{subcategories of $\D^{[-d+1,0]}$}\end{matrix}\right\}$};
        \draw(0,-4)node[rotate=-90]{$\subseteq$};
        \draw(-3,0)node{$\left\{\begin{matrix}\text{semibricks}\\\text{in $\D^{[-d+1,0]}$}\end{matrix}\right\}$};
        \draw(3,0)node{$\left\{\begin{matrix}\text{length wide}\\\text{subcategories}\\ \text{of $\D^{[-d+1,0]}$}\end{matrix}\right\}$};
        \draw(3,2)node{$\left\{\begin{matrix}\text{wide subcategories}\\\text{of $\D^{[-d+1,0]}$}\end{matrix}\right\}$};
        \draw[->,thick] (1.4,-.2) tonode[below]{$\ssim$} (-1.6,-.2);
        \draw[<-,thick] (1.4,.2) tonode[above]{$W$} (-1.6,.2);
        \draw(-.1,0)node{1-1} (3,1.2)node[rotate=90]{$\subseteq$};
        \draw[->,thick] (1,2) tonode[above left]{$\ssim$} (-3,.5);
        \draw[->,thick,commutative diagrams/hook](-3,-.5) tonode[right]{$T$} (-1.5,-2.2);
        \draw[->,thick,commutative diagrams/hook](-3.3,-.5) tonode[left]{$\phi$} (-2,-4.2);
        \draw[->,thick](3.9,-.7) tonode[right]{$\phi$} (2.6,-4.2);
        \draw[->,thick](2.4,-4.2) tonode[left]{$\LWlen$} (3.7,-.7);
        \draw[->,thick](2.7,-5) to[out=30,in=-60]node[right]{$\LW$} (5,1.5);
    \end{tikzpicture}
    \caption{A summary of the maps}
    \label{fig:outerbijs}
\end{figure}

Although Proposition~\ref{idforphi} and Corollary~\ref{cor:idforphi} suggest that the map $\phi$ behaves more favorably than $T$ for arbitrary semibricks and length wide subcategories, we shall see in the following section that $T$ becomes more advantageous upon restriction to certain subsets. Compare Figure~\ref{fig:outerbijs} with Figure~\ref{fig:main2}.

\section{Left/right-finite semibricks and wide subcategories}\label{S7}

For the remainder of this paper, we fix the triangulated category $\D=\D^b(\mod\La)$, and the $t$-structure 
$(\D^{\leq 0}, \D^{\geq 0})=(\D^{\leq 0}(\mod\La), \D^{\geq 0}(\mod\La)).$ 
In this setting, the $d$-extended heart is given by $\D^{[-d+1,0]} = \emod{d}{\La}$. In this section, we define left/right-finite semibricks and wide subcategories in $\emod{d}{\La}$ and establish their bijections with functorially finite positive torsion pairs in $\emod{d}{\La}$ and $(d+1)$-term simple-minded collections in $\D^b(\mod\La)$. We begin by restricting the bijection introduced in Proposition~\ref{prop:bi2} as follows.

\begin{proposition}\label{tors-pair-t-str}
    There is a bijection between 
    \begin{itemize}
        \item the set of functorially finite positive torsion pairs in $\emod{d}{\La}$, and 
        \item the set of bounded $t$-structures $(\C^{\leq 0},\C^{\geq 0})$ on $\D^b(\mod\La)$ with length heart satisfying $\D^{\leq -d}\subseteq \C^{\leq 0}\subseteq \D^{\leq 0}$.
    \end{itemize}
    The map is given by
    \[(\T,\F)\mapsto (\D^{\leq -d}*\T, \F[1]*\D^{\geq 0})\]
    with the inverse map defined by
    \[(\C^{\leq 0},\C^{\geq 0})\mapsto (\C^{\leq 0}\cap \emod{d}{\La}, \C^{\geq 1}\cap \emod{d}{\La}).\]
\end{proposition}

\begin{proof}
    First, recall that functorially finite positive torsion pairs in $\emod{d}{\La}$ coincide with functorially finite $s$-torsion pairs, as shown in Proposition~\ref{s-tors} (see also  \cite[Proposition~4.9]{Gu1}). By \cite[Theorem~4.1]{Gu1}, there is a bijection between basic $(d+1)$-term silting objects in $\K^b(\proj\La)$ and functorially finite positive torsion pairs in $\emod{d}{\La}$. Although \cite{Gu1} assumes throughout that the base field is algebraically closed, the proof of the cited theorem does not use this assumption, and hence the bijection remains valid over an arbitrary field. Under this bijection, a silting object $P$ is sent to the pair $(\T(P),\F(P))$, where
    \[\begin{array}{c}
         \T(P)=\{X\in \emod{d}{\La}\mid \Hom(P,X[i])=0,\ \forall i>0\},  \\
         \F(P)=\{X\in \emod{d}{\La}\mid \Hom(P,X[i])=0,\ \forall i\leq 0\}. 
    \end{array}\]
    On the other hand, Theorem~\ref{silting-simple} provides a bijection between basic $(d+1)$-term silting objects and bounded $t$-structures $(\C^{\leq 0},\C^{\geq 0})$ on $\D^b(\mod\La)$ with length hearts and satisfying $\D^{\leq -d}\subseteq\C^{\leq 0}\subseteq\D^{\leq 0}$, sending a silting object $P$ to $(\D^{\leq 0}(P),\D^{\geq 0}(P))$. Observing that $\T(P)=\D^{\leq 0}(P)\cap\emod{d}{\La}$ and $\F(P)=\D^{\geq 0}(P)[-1]\cap\emod{d}{\La}$, we obtain the desired bijection by combining these results and comparing them with Proposition~\ref{prop:bi2}.
\end{proof}

Let $\T$ be a functorially finite positive torsion class in $\emod{d}{\La}$. Its exact heart 
\[\H_\T=(\D^{\leq -d}*\T)\cap(\T^{\perp}[1]*\D^{\geq 0})\subseteq\D^{[-d,0]}\]
coincides with the heart of the corresponding bounded $t$-structure. Dually, for any functorially finite positive torsion-free class $\F$, its heart $\H^\F$ is defined as the heart of the corresponding bounded $t$-structure:
\[\H^\F:=(\D^{\leq -d}*{}^\perp\F)\cap(\F[1]*\D^{\geq 0})\subseteq\D^{[-d,0]}.\]
Thus, by Proposition~\ref{tors-pair-t-str} and Corollary~\ref{bi2cor}, we obtain the following two bijections:
\begin{equation}\label{eq:PhiLR}
\xymatrix{
    \left\{{\begin{matrix}\text{functorially finite positive}\\\text{torsion classes in $\emod{d}{\La}$}\end{matrix}}\right\}\ar[dr]^{\Theta_L} \\
    & \left\{{\begin{matrix}\text{length bounded hearts of}\\\text{$\D^b(\mod\La)$ lying in }\D^{[-d,0]}\end{matrix}}\right\}\\
    \left\{{\begin{matrix}\text{functorially finite positive}\\
    \text{torsion-free classes in $\emod{d}{\La}$}\end{matrix}}\right\}\ar[ur]_{\Theta_R}
    }
\end{equation}
These maps are given by $\Theta_L(\T)=\H_\T$ and $\Theta_R(\F)=\H^\F$, and their respective inverses are given by:
\[\Theta_L^{-1}(\H)=d\text{-}\H\cap \emod{d}{\La}\text{ and }\Theta_R^{-1}(\H)=d\text{-}\H[-d]\cap \emod{d}{\La}.\]

Recall that, according to Theorem~\ref{silting-simple}, there exists a bijection
\begin{equation}\label{eq:ssim}
    \ssim\colon\left\{{\begin{matrix}\text{length bounded hearts of}\\\text{$\D^b(\mod\La)$ lying in }\D^{[-d,0]}\end{matrix}}\right\}\to\left\{\begin{matrix}\text{$(d+1)$-term simple-minded}\\\text{collections in $\D^b(\mod\La)$}\end{matrix}\right\},
\end{equation}
which sends a length heart $\H$ to the set $\ssim \H$ of (isoclasses of) its simple objects.

\subsection{From simple-minded collections to semibricks}\label{subsec:sb}

Comparing Definitions~\ref{defn:smc} and \ref{semibrick}, one observes that for any $(d+1)$-term simple-minded collection $\X$, the intersections $\X\cap\emod{d}{\La}$ and $\X[-1]\cap\emod{d}{\La}$ are semibricks in $\emod{d}{\La}$. Consequently, we define two maps
\[\begin{array}{rccc}
    \Pi_L, \,\Pi_R\colon & \left\{\begin{matrix}\text{$(d+1)$-term simple-minded}\\\text{collections in $\D^b(\mod\La)$}\end{matrix}\right\} & \to & \{\text{semibricks in $\emod{d}{\La}$}\}
\end{array}\]
given by $\Pi_L(\X)=\X\cap\emod{d}{\La}$ and $\Pi_R(\X)=\X[-1]\cap\emod{d}{\La}$, respectively. In the sequel, we will prove that these maps are injective and provide a precise characterization of their images.

A consequence of Proposition~\ref{tors-pair-t-str} is the following useful lemma.

\begin{lemma}\label{neg-ext}
Let $\T$ be a functorially finite positive torsion class in $\emod{d}{\La}$, and let $\SS_\T=(\Pi_L\circ\ssim\circ\Theta_L)(\T)$ be the associated semibrick. Then  $\T \cap \SS^{\perp_{\leq 0}}_{\T}=0$.
\end{lemma}

\begin{proof}
    Let $X\in \T \cap \SS^{\perp_{\leq 0}}_{\T}$, and suppose $X\neq 0$. By Proposition~\ref{prop:bi2}(1), $\T$ is contained in the $d$-extended heart of $\H_\T$: \[ \T \subseteq d\text{-}\H_\T := \H_\T[d-1] * \H_\T[d-2] * \cdots * \H_\T. \]
    Thus, there exists an integer $0\leq l\leq d-1$ such that $X$ admits a triangle
    \[X'[l]\xrightarrow{h} X\to X'',\]
    where $0\neq X'\in \H_\T$ and $X''\in \H_\T[l-1]*\cdots\H_\T[1]*\H_\T$. Since $\H_\T$ is length by Proposition~\ref{tors-pair-t-str}, there exists a conflation in $\H_\T$:
    \[S \Rightarrowtail{f} X' \twoheadrightarrow Y,\]
    where $S\in\ssim \H_\T$. Applying the octahedral axiom to the composition $g=h\circ f[l]$, we obtain the following diagram of triangles:
    \begin{center}
         \begin{tikzcd}
        {S[l]} \arrow[d, "{f[l]}"'] \arrow[r, Rightarrow, no head] & {S[l]} \arrow[d, "g"] &     \\{X'[l]} \arrow[d] \arrow[r,"h"]           & X \arrow[d] \arrow[r] & X'' \arrow[d] \\
        {Y[l]} \arrow[r]                      & L \arrow[r]           & X'' \end{tikzcd}
    \end{center}
    By the triangle in the third row of the diagram, we have
    \[L\in Y[l]*X''\subseteq \H_\T[l]*\cdots\H_\T[1]*\H_\T\subseteq\D^{\geq -d-l},\]
    as $\H_\T \subseteq \D^{\geq -d}$. Then, the triangle in the second row implies:
    \[S[l]\in L[-1]*X\subseteq \D^{\geq -d-l+1}*\T\subseteq\D^{\geq -d-l+1},\]
    where we use the fact that $\T \subseteq \D^{\geq -d+1}\subseteq\D^{\geq -d-l+1}$. This implies $S \in \D^{\geq -d+1}$. Moreover, since $S \in \H_\T \subseteq \D^{\leq 0}$, we conclude that $S \in \emod{d}{\La}$, which implies $S \in \SS_\T$.
    
    Because $X \in \SS^{\perp_{\leq 0}}_{\T}$ and $S \in \SS_\T$, the morphism $g \colon S[l] \to X$ must be zero. This forces $f[l]$ to factor through $X''[-1]$. However, $S[l] \in \H_\T[l]$ and $X''[-1]\in \H_\T[l-2]*\cdots*\H_\T*\H_\T[-1]$. By the property of the heart, $\Hom(\H_\T, \H_\T[j]) = 0$ for $j<0$. Thus $\Hom(S[l], X''[-1]) = 0$, which implies $f[l] = 0$ and hence $f = 0$. This contradicts the fact that $S$ is a non-zero simple object. We conclude that $X = 0$.
\end{proof}

The following proposition states that, analogous to the classical case $d=1$ \cite[Remark~4.11]{BY}, a $(d+1)$-term simple-minded collection $\X$ is precisely the union of $\Pi_L(\X)$ and $\Pi_R(\X)[1]$. However, unlike the classical case, this union is not necessarily disjoint. For instance, consider the simple-minded collection $\X = \{P_2, I_1[1]\}$ from Example~\ref{exA2}, where $d=2$. In this case, $\Pi_L(\X) = \X$ and $\Pi_R(\X)[1] = \{I_1[1]\}$, illustrating the overlap.

\begin{proposition}\label{smc-to-semibrick-pair}
    Let $\X$ be a $(d+1)$-term simple-minded collection. Then 
    \[\X=\Pi_L(\X)\cup \Pi_R(\X)[1].\]
\end{proposition}

\begin{proof}
    The inclusion $\Pi_L(\X) \cup \Pi_R(\X)[1] \subseteq \X$ follows directly from the definitions of $\Pi_L$ and $\Pi_R$. For the reverse inclusion, let $X \in \X$. Suppose $X \notin \Pi_L(\X)$, which means $X \notin \emod{d}{\La}$. We aim to show that $X[-1] \in \emod{d}{\La}$.
    
    Let $\T$ be the functorially finite positive torsion class in $\emod{d}{\La}$ associated with $\X$, that is, $\X=\ssim \H_\T$. Since $\H_\T\subseteq\D^{[-d,0]}$ and $X\notin \emod{d}{\La}=\D^{[-d+1,0]}$, the truncation of $X$ with respect to the standard $t$-structure yields a triangle:
    \begin{equation}\label{eq:trunX}
        \sigma_{\leq -d}(X)\xrightarrow{f}X\to \sigma_{\geq -d+1}(X),
    \end{equation}
    where $f\neq 0$ and $\sigma_{\leq -d}(X)\in\D^{[-d,-d]}$. By Proposition~\ref{prop:bi2}(4), we have $\sigma_{\leq -d}(X)\in\T^\perp[d]\text{ and }\sigma_{\geq -d+1}(X)\in\T.$
    
    Recall from Proposition~\ref{prop:bi2}(1) that $\T^\perp[d] \subseteq \H_\T[d-1] * \H_\T[d-2] * \cdots * \H_\T$. Consequently, there exists a triangle
    \[Z\xrightarrow{g} \sigma_{\leq -d}(X)\xrightarrow{h} Y\]
    with $Z \in \H_\T[d-1] * \cdots * \H_\T[1]$ and $Y \in \H_\T$. Since $X\in\H_\T$, we have $\Hom(Z,X)=0$, which implies $f\circ g=0$. Hence, there exists a morphism $f':Y\to X$ such that $f=f'\circ h$. Thus, by the octahedral axiom, we obtain the following commutative diagram:
    \begin{center}
        \begin{tikzcd}
        Z \arrow[r, Rightarrow, no head] \arrow[d]    & Z \arrow[d, "g"]   & \\
        {(\sigma_{\geq -d+1}(X))[-1]} \arrow[r] \arrow[d] & \sigma_{\leq -d}(X) \arrow[d, "h"] \arrow[r, "f"] & X \arrow[d, Rightarrow, no head] \\ Y' \arrow[r] & Y \arrow[r, "f'"'] & X         \end{tikzcd}
    \end{center}
    Since $f \neq 0$, the induced map $f' \colon Y \to X$ is also non-zero. Given that $X$ is a simple object in $\H_\T$ and $Y \in \H_\T$, the map $f'$ must be an epimorphism in the abelian category $\H_\T$, which implies that the triangle in the bottom row is a short exact sequence in $\H_\T$. Thus, in particular, we have $Y' \in \H_\T$.
    
    From the triangle in the left column, we have $\sigma_{\geq -d+1}(X) \in Z[1] * Y'[1]$. Since $Z[1] \in \H_\T[d] * \dots * \H_\T[2]$ and $Y'[1] \in \H_\T[1]$, we have 
    \[ \sigma_{\geq -d+1}(X) \in \H_\T[d] * \dots * \H_\T[1] \subseteq \D^{\leq -1}.\]
    It follows that $\sigma_{\geq -d+1}(X)\in\D^{[-d+1,-1]}$. Therefore, the triangle~\eqref{eq:trunX} shows that $X\in\D^{[-d,-1]}$, which forces $X[-1]\in \emod{d}{\La}$. This completes the proof.
\end{proof}

Recall from Section~\ref{sec:smc} that the left (resp. right) mutation of a $(d+1)$-term simple-minded collection $\X$ at an element $X_i$ exists, if the resulting simple-minded collection $\mu_i^+(\X)$ (resp. $\mu_i^-(\X)$) is again $(d+1)$-term.

\begin{proposition}\label{smc-mutation}
Let $\X=\{X_i\}_{i=1}^n$ be a $(d+1)$-term simple-minded collection.
\begin{enumerate}
    \item[(1)] The left mutation of $\X$ at $X_i$ exists if and only if $X_i\in \Pi_L(\X)$.
    \item[(2)] The right mutation of $\X$ at $X_i$ exists if and only if $X_i\in \Pi_R(\X)[1]$.
\end{enumerate}
Moreover, for any $1\leq i\leq n$, at least one of the left and right mutations of $X$ at $X_i$ exists.
\end{proposition}

\begin{proof}
Suppose that the left mutation of $\X$ at $X_i$ exists. By definition, all elements of the mutated collection $\mu_i^+(\X)$ must lie in the $(d+1)$-extended heart $\D^{[-d,0]}$. In particular, $X_i[1] \in \D^{[-d,0]}$, which implies $X_i \in \D^{\geq -d+1}$. Since $X_i \in \X \subseteq \D^{\leq 0}$, it follows that $X_i \in \emod{d}{\La}$, and thus $X_i \in \Pi_L(\X)$.

Conversely, suppose $X_i \in \Pi_L(\X)$. The collection $\mu_i^+(\X)$ consists of $\{X'_j\}_{j \neq i} \cup \{X_i[1]\}$. Since $X_i \in \emod{d}{\La}=\D^{[-d+1,0]}$, we have $X_i[1] \in \D^{[-d, -1]}\subseteq \D^{[-d, 0]}$. For each $j \neq i$, the object $X'_j$ is defined by a triangle
\[Y_j \to X'_j \to X_j \to Y_j[1],\]
where $Y_j$ lies in the extension closure of $X_i$. Since both $X_j$ and $X_i$ are in $\D^{[-d, 0]}$, it follows that $X'_j \in \D^{[-d, 0]}$. Thus, $\mu_i^+(\X)$ is a $(d+1)$-term simple-minded collection. The argument for the right mutation is analogous.

Finally, by Proposition~\ref{smc-to-semibrick-pair}, every $X_i$ is either in $\Pi_L(\X)$ or in $\Pi_R(\X)[1]$. Therefore, at least one of the left mutation and the right mutation of $\X$ at $X_i$ exists.
\end{proof}

\subsection{Left/right-finite semibricks}

Recall that for any semibrick $\SS$ in $\D^{[-d+1,0]}$, we denote by $T(\SS)$ the smallest positive torsion class containing $\SS$. Dually, we denote by $F(\SS)$ the smallest positive torsion-free class containing $\SS$. Thus, $T$ and $F$ define maps from the set of semibricks in $\emod{d}{\La}$ to the set of positive torsion classes and positive torsion-free classes in $\emod{d}{\La}$, respectively. We are particularly interested in those semibricks $\SS$ whose associated classes $T(\SS)$ or $F(\SS)$ satisfy the functorial finiteness condition.

\begin{definition}\label{def:l.f.semibrick}
    A semibrick $\SS\subseteq \emod{d}{\La}$ is called \emph{left-finite} (resp. \emph{right-finite}) if $T(\SS)$ (resp. $F(\SS)$) is functorially finite.
\end{definition}

Our first main result is summarized in the following theorem.

\begin{theorem}\label{main}
    There exists a commutative diagram of bijections shown in Figure~\ref{fig:main1}, such that each oriented digon composes to the identity.
    \begin{figure}[htpb]
        \centering
        \begin{tikzpicture}
    \draw(0,1)node{$\left\{\begin{matrix}\text{$(d+1)$-term simple-minded}\\\text{collections in $\D^b(\mod\La)$}\end{matrix}\right\}$};
    \draw(0,-1)node{$\left\{\begin{matrix}\text{length bounded hearts}\\\text{lying in $\D^{[-d,0]}(\mod\La)$}\end{matrix}\right\}$};
    \draw[thick,->] (0,-.5) tonode[left, font=\tiny]{$\ssim$} (0,.5);
    \draw(4.5,4)node{$\left\{\begin{matrix}\text{right-finite semibricks}\\\text{in $\emod{d}{\La}$}\end{matrix}\right\}$};
    \draw(-4.5,4)node{$\left\{\begin{matrix}\text{left-finite semibricks}\\\text{in $\emod{d}{\La}$}\end{matrix}\right\}$};
    \draw(4.5,-4)node{$\left\{\begin{matrix}\text{functorially finite}\\\text{positive torsion-free}\\\text{classes in $\emod{d}{\La}$}\end{matrix}\right\}$};
    \draw(-4.5,-4)node{$\left\{\begin{matrix}\text{functorially finite}\\\text{positive torsion}\\\text{classes in $\emod{d}{\La}$}\end{matrix}\right\}$};
    \draw[->, thick] (-.5,1.5) tonode[sloped, midway, above, font=\tiny]{$\X\cap\emod{d}{\La}\mapsfrom\X$}node[sloped, midway, below, font=\tiny]{$\Pi_L$} (-4,3.5);
    \draw[->, thick] (.5,1.5) tonode[sloped, midway, above, font=\tiny]{$\X\mapsto\X[-1]\cap\emod{d}{\La}$}node[sloped, midway, below, font=\tiny]{$\Pi_R$} (4,3.5);
    \draw[->, thick] (-4.5,3.5)tonode[right, font=\tiny]{$T$} (-4.5,-3.2);
    \draw[->, thick] (4.5,3.5)tonode[left, font=\tiny]{$F$} (4.5,-3.2);
    \draw[->, thick] (-3.8,-3.2) tonode[sloped, midway, below, font=\tiny]{$\Theta_L\colon\T\mapsto\H_\T$} (-.3,-1.5);
    \draw[<-, thick] (-4.2,-3.2) tonode[sloped, midway, above, font=\tiny]{$d\text{-}\H\cap \emod{d}{\La}\mapsfrom\H\quad$} (-.7,-1.5);
    \draw[<-, thick] (4.2,-3.2) tonode[sloped, midway, above, font=\tiny]{$\qquad\H\mapsto d\text{-}\H[-d]\cap \emod{d}{\La}$} (.7,-1.5);
    \draw[->, thick] (3.8,-3.2) tonode[sloped, midway, below, font=\tiny]{$\H^\F\mapsfrom\F\colon\Theta_R$} (.3,-1.5);
    \draw[->,thick] (-2.4,-3.9) tonode[above, font=\tiny]{$\T\mapsto\T^\perp$} (2.4,-3.9);
    \draw[<-,thick] (-2.4,-4.1) tonode[below, font=\tiny]{$^\perp\F\mapsfrom\F$} (2.4,-4.1);
\end{tikzpicture}
        \caption{A commutative diagram of bijections}
        \label{fig:main1}
    \end{figure}
\end{theorem}

\begin{proof}
    The bijections $\Theta_L$, $\Theta_R$ (and their inverses) and $\ssim$ were established in \eqref{eq:PhiLR} and \eqref{eq:ssim} respectively. The lower triangles commute since, for any functorially finite positive torsion pair $(\T, \F)$, we have the identity $\H_\T = \H^{\F}$ by construction.
    
    Let $\T$ be a functorially finite positive torsion class in $\emod{d}{\La}$, and let $\SS_\T = (\Pi_L \circ \ssim \circ \Theta_L)(\T)$. By definition, we have
    \[\SS_\T=\ssim \H_\T\cap \emod{d}{\La}\subseteq\H_\T\cap\emod{d}{\La}\subseteq d\text{-}\H_\T\cap\emod{d}{\La}=\T.\]
    Since $\T$ is a positive torsion class containing $\SS_\T$, it follows that $T(\SS_\T) \subseteq \T$. Combining Lemma~\ref{same-tor-free} and Lemma~\ref{neg-ext}, we obtain $\T \cap T(\SS_\T)^\perp = \T \cap \SS^{\perp_{\leq 0}}_\T = 0$. Since $\T$ is an $s$-torsion class and $T(\SS_\T)$ is a positive torsion class in $\emod{d}{\La}$ with $T(\SS_\T) \subseteq \T$, Lemma~\ref{no-middle} implies $\T = T(\SS_\T)$. This proves the identity
    \[ T \circ \Pi_L \circ \ssim \circ \Theta_L = \mathrm{id}. \]
    
    Consequently, the restriction of the map $T$ to the set of left-finite semibricks is surjective. Moreover, because $\ssim$ and $\Theta_L$ are bijections, the map $\Pi_L$ is injective and its image is contained within the set of left-finite semibricks. The identity above also ensures the commutativity of the left part of the diagram. By Corollary~\ref{inj}, $T$ is injective, and thus its restriction is a bijection. Finally, this implies that the image of $\Pi_L$ coincides exactly with the set of left-finite semibricks.
    
    A dual argument applies to the right-hand side.
\end{proof}

The following result is an immediate consequence of Theorem~\ref{main} and \cite[Corollary~5.5]{KY}. For the case $d=1$, it recovers \cite[Corollary~2.10]{A}.

\begin{corollary}
    Let $\SS$ be a semibrick in $\emod{d}{\La}$. If $\SS$ is either left-finite or right-finite, then the cardinality of $\SS$ is at most the rank of $\La$.
\end{corollary}

\subsection{Left/right-finite wide subcategories}

Recall that for any subcategory $\Y$ of $\emod{d}{\La}$, we denote by $T(\Y)$ (resp. $F(\Y)$) the smallest positive torsion class (resp. torsion-free class) containing $\Y$.

\begin{definition}\label{l.f.wide}
    A wide subcategory $\W \subseteq \emod{d}{\La}$ is called \emph{left-finite} (resp. \emph{right-finite}) if the smallest positive torsion class $T(\W)$ (resp. torsion-free class $F(\W)$) containing $\W$ is functorially finite.
\end{definition}

Note that for any semibrick $\SS$ in $\emod{d}{\La}$, we have the following identities
\begin{equation}\label{eq:T(S)=T(W)}
    T(\SS)=T(W(\SS))\text{ and }F(\SS)=F(W(\SS)),
\end{equation}
where $W$ is the bijection established in Proposition~\ref{sb-wide}; specifically, $W(\SS)$ is the extension closure of $\SS$. This implies that a semibrick $\SS$ is left-finite (resp. right-finite) if and only if $W(\SS)$ is left-finite (resp. right-finite).

\begin{proposition}\label{finite-wide-fin-len}
    Any left-finite or right-finite wide subcategory of $\emod{d}{\La}$ is a length category.
\end{proposition}

\begin{proof}
    Let $\W$ be a left-finite wide subcategory. By definition,  $\T:=T(\W)$ is a functorially finite positive torsion class. We claim that $\W\subseteq \H_{\T}$. If this claim holds, $\W$ becomes a wide subcategory of the length abelian category $\H_\T$ (as per Proposition~\ref{tors-pair-t-str}), which directly implies that $\W$ is itself a length category.
    
    To prove the claim, we follow a strategy similar to Step 1 of the proof of Proposition~\ref{sim-in-heart}. First, the inclusion $\W \subseteq \T \subseteq \D^{\leq -d} * \T$ is immediate from the definition of $\T$. Next, let $X \in \W$. For the truncation triangle $\sigma_{\leq -1}(X) \to X \to \sigma_{\geq 0}(X) \to \sigma_{\leq -1}(X)[1]$, applying the functor $\Hom(\W, -)$ yields the following exact sequence for any $j \in \mathbb{Z}$:
    \[\Hom(\W,(\sigma_{\geq 0}(X))[j])\to \Hom(\W,(\sigma_{\leq -1}(X))[j+1])\to \Hom(\W, X[j+1]).\]
    For $j \leq -2$, the first term $\Hom(\W, \sigma_{\geq 0}(X)[j])$ vanishes because $\W \subseteq \D^{\leq 0}$ and $\sigma_{\geq 0}(X)[j] \in \D^{\geq -j} \subseteq \D^{\geq 2}$. The last term $\Hom(\W, X[j+1])$ also vanishes because $\W$ is a wide subcategory and $j+1 < 0$. This implies that $\Hom(\W, \sigma_{\leq -1}(X)[k]) = 0$ for all $k \leq -1$, or equivalently, $(\sigma_{\leq -1}(X))[-1] \in \W^{\perp_{\leq 0}}$.
    
    Recall that $\W^{\perp_{\leq 0}} = \T^\perp$ by Proposition~\ref{perp:tor}. Thus, $(\sigma_{\leq -1}(X))[-1] \in \T^\perp$, which means $\sigma_{\leq -1}(X) \in \T^\perp[1]$. Consequently, 
    \[ X \in \sigma_{\leq -1}(X) * \sigma_{\geq 0}(X) \subseteq \T^\perp[1] * \D^{\geq 0}. \]
    Combined with $\W \subseteq \D^{\leq -d} * \T$, we conclude $\W \subseteq \H_\T$. The case for right-finite wide subcategories is analogous.
\end{proof}

Recall from \eqref{eq:mapLW} that we defined a map $\LW$, which sends a $d$-FAE closed subcategory of $\emod{d}{\La}$ to a wide subcategory of $\emod{d}{\La}$. Our second main result establishes the following commutative diagram, compared with Figure~\ref{fig:outerbijs}.

\begin{theorem}\label{main2}
    There is a commutative triangle of bijections shown in Figure~\ref{fig:main2}. 
    \begin{figure}[htpb]
        \centering
        \begin{tikzpicture}
        \draw(-4,1.5)node{$\left\{{\begin{matrix}\text{left-finite semibricks}\\\text{in $\emod{d}{\La}$}\end{matrix}}\right\}$};
        \draw(4,1.5)node{$\left\{{\begin{matrix}\text{left-finite wide}\\\text{subcategories in $\emod{d}{\La}$}\end{matrix}}\right\}$};
        \draw(0,-1.5)node{$\left\{{\begin{matrix}\text{functorially finite, positive}\\\text{torsion classes in $\emod{d}{\La}$}\end{matrix}}\right\}$};
        \draw[->,thick] (-1.7,1.5) tonode[above]{$W$} (1.4,1.5);
        \draw[->,thick] (-4,.9) tonode[below left]{$T$} (-.3,-.9);
        \draw[->,thick] (.3,-.9) tonode[below right]{$\LW$} (4, .9);
    \end{tikzpicture}
        \caption{A commutative triangle of bijections}
        \label{fig:main2}
    \end{figure}
\end{theorem}

\begin{proof}
    By Theorem~\ref{main}, the map $T$ restricts to a bijection between left-finite semibricks and functorially finite positive torsion classes. Moreover, Proposition~\ref{finite-wide-fin-len} and the identities in \eqref{eq:T(S)=T(W)} ensure that the bijection $W$ from Proposition~\ref{sb-wide} restricts to a bijection between the set of left-finite (resp. right-finite) semibricks and the set of left-finite (resp. right-finite) wide subcategories.

    For any left-finite semibrick $\SS$, the torsion class $T(\SS)$ is functorially finite and positive by definition. According to Proposition~\ref{s-tors}, such a class is necessarily an $s$-torsion class. This allows us to apply Proposition~\ref{idforphi}, which yields the identity $\SS = \ssim \LW(T(\SS))$. Since $W(\SS)$ is also characterized by its simple objects (i.e., $\SS = \ssim W(\SS)$), the commutativity of the triangle $W=\LW\circ T$ follows. Consequently, the map $\LW$ restricts to the desired bijection as shown in the diagram.
\end{proof}

\subsection{Realizing left/right-finite wide subcategories as module categories}\label{subset:realasmod}

In the classical case $d=1$, left-finite wide subcategories are known to be equivalent to module categories of certain quotient algebras \cite[Theorem~3.15]{A}. In this subsection, we show that this property extends to the $d$-setting.

Let $P=\bigoplus_{i=1}^n P_i$ be a basic $(d+1)$-term silting complex, and let 
\begin{equation}\label{eq:heart}
    \H(P)=\{X\in\D^b(\mod\La)\mid\Hom(P,X[i])=0,\ \forall \ i\neq 0\}
\end{equation}
be the heart of the corresponding $t$-structure (see Theorem~\ref{silting-simple}). We recall the equivalence of categories from \cite{KY}:
\begin{equation}\label{eq:equiv}
\Hom(P,-):\H(P)\to \mod C,
\end{equation}
where $C:=\End(P)$. Let $\X=\ssim \H(P)=\{X_1,\dots,X_n\}$ be the corresponding $(d+1)$-term simple-minded collection, ordered such that $\Hom(P_i,X_j)=0$ for $i\neq j$.

By Proposition~\ref{smc-to-semibrick-pair}, the simple-minded collection $\X$ is the union of its left and (shift of) right semibrick parts. Without loss of generality, we may assume that the left-finite semibrick is
\[\SS_L=\Pi_L(\X)=\{X_1,\dots, X_l, X_{l+1}, \dots, X_m\}\]
and the right-finite semibrick is
\[\SS_R=\Pi_R(\X)=\{X_{l+1}[-1],\dots, X_m[-1], X_{m+1}[-1],\dots, X_n[-1]\},\]
where the indices reflect the possible overlap between $\Pi_L(\X)$ and $\Pi_R(\X)[1]$. Let $\W = W(\SS_L)$ and $\W' = W(\SS_R)$ be the corresponding wide subcategories in $\emod{d}{\La}$.

Let $f_i\in C$ be the primitive idempotent corresponding to $P_i$, and set $f=\sum_{i=m+1}^n f_i$. Dually, considering the Nakayama functor $\nu \colon \K^b(\proj \La) \to \K^b(\inj \La)$, the algebra $C$ is isomorphic to $C' := \End(\nu P)$, and the functor
\[D\Hom(-, \nu P): \H(P) \to \mod C'\]
is an equivalence. Let $f'_i \in C'$ be the primitive idempotent corresponding to $\nu P_i$, and set $f' = \sum_{i=1}^l f'_i$.

\begin{theorem}
Any left-finite (resp. right-finite) wide subcategory of $\emod{d}{\La}$ is equivalent to the module category of a finite-dimensional algebra. More precisely:
\begin{enumerate}
    \item The equivalence $\Hom(P,-)$ restricts to an equivalence
    \[\W\simeq \mod (C/\langle f\rangle).\]
    \item The equivalence $D\Hom(-, \nu P)$ restricts to an equivalence
    \[\W'[1]\simeq \mod (C'/\langle f'\rangle).\]
\end{enumerate}
\end{theorem}

\begin{proof}
(1) The equivalence $\H(P) \simeq \mod C$ maps the wide subcategory $\W$ to the extension closure of the simple $C$-modules $\{\Hom(P, X_i) \mid i=1, \dots, m\}$. In $\mod C$, such an extension closure is precisely the category of modules whose composition factors are among these $m$ simples, which is equivalent to $\mod (C / \langle \sum_{i=m+1}^n f_i \rangle)$. The reader is referred to \cite[Theorem~3.15]{A} for further details. The proof of (2) follows dually.
\end{proof}

\section{Mutations of \texorpdfstring{$(d+1)$}{(d+1)}-term silting complexes}\label{S8}

In this final section, we provide a criterion to determine which mutations of a $(d+1)$-term silting complex $P$ remain $(d+1)$-term. We achieve this by establishing a relationship between the elements of the simple-minded collection associated with $P$ and certain objects in the additive subcategory $\add P$.

Throughout this section, let $P$ be a basic $(d+1)$-term silting complex in $\K^b(\proj \La)$. We begin with the following technical lemma concerning the existence of a specific sequence of triangles.

\begin{lemma}
There exist $d$ triangles in $\K^b(\proj \La)$ of the form
\begin{equation}\label{triangles}
    \begin{array}{rcccl}
         Z_0=\La & \xrightarrow{f_0} & Q_0 & \xrightarrow{g_1} & Z_1 \\
         Z_1 & \xrightarrow{f_1} & Q_1 & \xrightarrow{g_2} & Z_2\\
         &&\vdots&&\\
         Z_{d-1} & \xrightarrow{f_{d-1}} & Q_{d-1} & \xrightarrow{g_d} & Q_d=Z_d
    \end{array}
\end{equation}
such that for each $0 \leq j \leq d$, the following conditions hold:
\begin{enumerate}
    \item $Q_j\in \add P$;
    \item $f_j$ is a minimal left $(\add P)$-approximation of $Z_j$, with $f_d=\id_{Z_d}$;
    \item $\Hom(Z_j, P[l])=0$ for all $l\geq 1$;
    \item $\Hom(P,Z_j[l])=0$ for all $l\geq d-j+1$.
\end{enumerate}
\end{lemma}

\begin{proof}
    We construct the triangles by recursively taking $f_j\colon Z_j \to Q_j$ to be a minimal left $(\add P)$-approximation of $Z_j$ and define $Z_{j+1}:=C(f_j)$ for $0 \leq j \leq d-1$. By this construction, (1) and (2) hold for all $0\leq j\leq d-1$.
    
    We proceed to show that conditions (3) and (4) hold by induction on $j$.
    
    For the starting case $j=0$, we have $Z_0 = \La$. Since $P$ is a $(d+1)$-term silting complex, it follows that $\Hom(Z_0, P[l]) = 0$ for all $l \geq 1$ and $\Hom(P, Z_0[l]) = 0$ for all $l \geq d+1$, satisfying the requirements.
    
    Assume that (3) and (4) hold for all $j \leq t < d$. For $j = t+1$, applying the functor $\Hom(P, -)$ to the $(t+1)$-th triangle in \eqref{triangles} yields the exact sequence
    \[\Hom(P,Q_t[l])\to \Hom(P,Z_{t+1}[l])\to \Hom(P,Z_t[l+1]).\]
    By the induction hypothesis, the third term vanishes for all $l+1 \geq d-t+1$. Since $P$ is silting and $Q_t \in \add P$, the first term vanishes for all $l > 0$. Thus, $\Hom(P, Z_{t+1}[l]) = 0$ for any $l \geq d-t$, which confirms (4) for $j = t+1$.
    
    Next, applying $\Hom(-, P)$ to the same triangle, we obtain the exact sequence
    \[\Hom(Q_t, P[l])\to \Hom(Z_t, P[l])\to \Hom(Z_{t+1}, P[l+1])\to \Hom(Q_t,P[l+1]).\]
    For $l\geq 1$, the second term vanishes for $l\geq 1$ by the induction hypothesis. For $l\geq 0$, the last term vanishes for $l\geq 0$ since $P$ is silting. For $l = 0$, the first map is surjective because $f_t$ is a left $(\add P)$-approximation. We conclude that $\Hom(Z_{t+1}, P[l+1]) = 0$ for all $l \geq 0$, which is equivalent to $\Hom(Z_{t+1}, P[l]) = 0$ for any $l \geq 1$. Thus, (3) holds for $j = t+1$.
    
    Finally, we show that $Z_d \in \add P$. Since $Z_d \in \K^b(\proj \La)$ and $P$ is a silting complex, by \cite[Proposition~2.23~(b)]{AI}, there exists a triangle
    \[X\xrightarrow{f} Z_d\xrightarrow{g}Y,\]
    where $X \in P[-s] * \cdots * P$ and $Y \in P[1] * \cdots * P[s]$ for some $s > 0$. Applying (3) to $Z_d$, we have $\Hom(Z_d, P[l]) = 0$ for any $l \geq 1$, which implies $\Hom(Z_d, Y) = 0$. Thus, $g = 0$, and $Z_d$ is a direct summand of $X$. Since $P[-s] * \cdots * P$ is closed under direct summands \cite[Proposition~2.1]{IY}, we have $Z_d \in P[-s] * \cdots * P$.
    
    Consequently, there exists a triangle $V \xrightarrow{h} Z_d \to Z$ such that $Z \in \add P$ and $V \in P[-s] * \cdots * P[-1]$. Applying (4) to $Z_d$, we have $\Hom(P, Z_d[l]) = 0$ for all $l \geq 1$, which implies $\Hom(V, Z_d) = 0$. Thus $h = 0$, and we conclude that $Z_d$ is a direct summand of $Z$ and hence belongs to $\add P$.
\end{proof}

Let $h_j := f_j \circ g_j \colon Q_{j-1} \to Q_j$ for $1\leq j\leq d$. Let $\H(P)$ be the heart of the bounded $t$-structure associated with $P$; see \eqref{eq:heart}. For brevity, we denote $\Hom(X, Y)$ by $(X, Y)$ in the following proofs.

\begin{lemma}\label{distinct-top}
    For each $X\in \ssim \H(P)$, we have $\Hom(h_j,X)=0$ for all $1\leq j\leq d$.
\end{lemma}

\begin{proof}
    Since $f_{j-1}$ is a minimal left approximation, the map $g_j \in \rad(Q_{j-1}, Z_j)$, which implies that the composition $h_j$ lies in $\rad(Q_{j-1}, Q_j)$. We consider the following commutative diagram:
    \begin{center}
        \begin{tikzcd}
        {(Q_j,X)} \arrow[d, "\cong"', sloped] \arrow[rrr, "{(h_j,X)}"] & & & {(Q_{j-1},X)} \arrow[d, "\cong", sloped] \\
        {((P,Q_j),(P,X))} \arrow[rrr, "{((P,h_j),(P,X))}"']       & & & {((P,Q_{j-1}),(P,X))}           
        \end{tikzcd}
    \end{center}
    The vertical isomorphisms are induced by the Yoneda lemma, as $X \in \H(P)$ and $Q_j, Q_{j-1} \in \add P$. Since $h_j \in \rad(Q_{j-1}, Q_j)$, the equivalence of categories 
    \[ \add P \xrightarrow{(P, -)} \proj \End(P) \]
    ensures that $(P, h_j) \in \rad((P, Q_{j-1}), (P, Q_j))$. Note that $(P, Q_{j-1})$ and $(P, Q_j)$ are projective modules in $\mod \End(P)$. Consequently, the image of $(P, h_j)$ is contained in $\rad (P, Q_j)$.
    
    It follows from the equivalence \eqref{eq:equiv} that $(P, X)$ is a simple module in $\mod \End(P)$. Any morphism from a projective module to a simple module must vanish on the radical of the projective module. Therefore, the morphism $((P, h_j), (P, X))$ must be zero, which implies $(h_j, X) = 0$ as required. 
\end{proof}

Suppose $P = \bigoplus_{i=1}^n P_i$ with $P_i$ indecomposable. Let $\X = \{X_i\}_{i=1}^n = \ssim \H(P)$ be the associated simple-minded collection, ordered such that $(P_i, X_j[m]) = 0$ for $i \neq j$ or $m \neq 0$, and $(P_i, X_i) \cong \End(X_i)$. The following proposition demonstrates that for each $0 \leq j \leq d$, $\Ho^{-j}(X_i)$ is nonzero if and only if $P_i \in \add Q_j$.

\begin{proposition}\label{cohom-of-simples}
    For any $X\in \ssim \H(P)$ and each $0\leq j\leq d$, there is an isomorphism
    \[\Hom(\La[j],X)\cong \Hom(Q_j, X).\]
\end{proposition}

\begin{proof}
    \textbf {Claim 1: $\boldsymbol{(Z_j[l],X)=0}$ for any $\boldsymbol{l<0}$ and $0 \leq j \leq d$.} This is immediate for $j=d$ since $Z_d = Q_d \in \add P$ and $X \in \H(P)$. For $0 \leq j \leq d-1$, applying $(-, X)$ to the triangles in \eqref{triangles} yields the exact sequence
    \begin{equation}\label{eq:extseq}
        (Z_{j+1}[l],X) \to (Q_j[l],X) \to (Z_j[l],X) \to (Z_{j+1}[l-1],X) \to (Q_j[l-1],X), \quad l\in\mathbb{Z}.
    \end{equation}
    For any $l < 0$, the terms involving $Q_j$ vanish because $Q_j \in \add P$ and $X \in \H(P)$. Thus, we obtain the following isomorphism chain:
    \[(Z_j[l], X) \cong (Z_{j+1}[l-1], X) \cong \cdots \cong (Q_d[j+l-d], X) \cong 0.\]
    
    \textbf{Claim 2: $\boldsymbol{(f_j,X):(Q_j,X)\to (Z_j,X)}$ is an isomorphism for $\boldsymbol{0\leq j\leq d}$, and $\boldsymbol{(g_j,X):(Z_j,X)\to (Q_{j-1},X)}$ is zero for $\boldsymbol{1\leq j\leq d}$.} We proceed by induction on $j$ from $d$ down to $0$. For the starting case $j=d$, $f_d = \mathrm{id}_{Z_d}$ and $g_d=h_d$, so it is trivial that $(f_j,X)$ is an isomorphism, and $(g_d, X) = 0$ follows directly from Lemma~\ref{distinct-top}. Now, assume the statement holds for $j+1$ and consider $0 \leq j \leq d-1$. The exact sequence \eqref{eq:extseq} for $l=0$ is:
    \[(Z_{j+1}, X) \xrightarrow{(g_{j+1}, X)}(Q_{j},X)\xrightarrow{(f_j,X)} (Z_{j},X)\to (Z_{j+1}[-1],X).\]
    The last term vanishes by Claim 1, and the first map $(g_{j+1}, X)$ is zero by the induction hypothesis. Thus, $(f_j, X)$ is an isomorphism. Furthermore, for $j \geq 1$, Lemma~\ref{distinct-top} implies $(h_j, X) = (g_j, X) \circ (f_j, X) = 0$. Since $(f_j, X)$ is an isomorphism, we must have $(g_j, X) = 0$.
    
    \textbf{Claim 3: $\boldsymbol{(Z_{j}[1],X)\cong (Z_{j+1},X)}$ for all $\boldsymbol{0\leq j\leq d-1}$.} The exact sequence \eqref{eq:extseq} for $l=1$ is 
    \[(Q_j[1], X) \to (Z_j[1], X) \to (Z_{j+1}, X) \xrightarrow{(g_{j+1}, X)} (Q_j, X). \]
    Since $(Q_j[1], X) = 0$ and $(g_{j+1}, X) = 0$ (from Claim 2), the middle map is an isomorphism.
    
    \textbf{Claim 4: $\boldsymbol{(\La[j],X)\cong (Z_{j-1}[1], X)}$ for all $\boldsymbol{1\leq j\leq d}$.} Observing that for $l>1$, the terms in the exact sequence~\eqref{eq:extseq} involving $Q_j$ also vanish as $Q_j \in \add P$ and $X\in\H(P)$. Thus, we obtain the isomorphism chain:
    \[(\La[j],X)=(Z_0[j],X)\cong (Z_{1}[j-1],X)\cong\cdots\cong (Z_{j-1}[1],X).\]
    
    \textbf{Conclusion:} By combining these claims, for all $1 \leq j \leq d$, we have:
    \[ (\La[j], X) \stackrel{\text{C4}}{\cong} (Z_{j-1}[1], X) \stackrel{\text{C3}}{\cong} (Z_j, X) \stackrel{\text{C2}}{\cong} (Q_j, X).\]
    For $j=0$, the isomorphism $(\La, X) = (Z_0, X) \cong (Q_0, X)$ follows directly from Claim 2. This completes the proof.
\end{proof}

Recall from Section~\ref{sec:sil} that the left (resp. right) mutation of a $(d+1)$-term silting object $P$ at a direct summand $P_i$ exists if $\mu_i^+(P)$ (resp. $\mu_i^-(P)$) remains a $(d+1)$-term silting object. We obtain the following criterion for the existence of such mutations.

\begin{corollary}
    Let $P=\bigoplus_{i=1}^n P_i$ be a basic $(d+1)$-term silting complex with $P_i$ indecomposable. 
    \begin{enumerate}
    \item The left mutation of $P$ at $P_i$ exists if and only if $P_i \notin \add Q_d$.
    \item The right mutation of $P$ at $P_i$ exists if and only if $P_i \notin \add Q_0$.
\end{enumerate}
Furthermore, for any indecomposable summand $P_i$, at least one of the left and right mutations of $P$ at $P_i$ exists.
\end{corollary}

\begin{proof}
    Let $\X = \{X_i\}_{i=1}^n = \ssim \H(P)$ be the simple-minded collection associated with $P$. By \cite[Theorem~7.12]{KY}, the left mutation of $P$ at $P_i$ exists if and only if the left mutation of $\X$ at $X_i$ exists. According to Proposition~\ref{smc-mutation}, this is equivalent to $X_i \in \Pi_L(\X)$, which means $X_i \in \emod{d}{\La}$.
    
    By Proposition~\ref{cohom-of-simples}, we have $\Hom(\La[d], X_i) \cong \Hom(Q_d, X_i)$. Since $X_i \in \D^{[-d, 0]}$, the condition $X_i \in \emod{d}{\La} = \D^{[-d+1, 0]}$ is equivalent to the vanishing of $\Hom(\La[d], X_i)$. In light of the isomorphism $\Hom(Q_d, X_i) \cong \Hom(P_i, X_i)^{\oplus m_i}$ where $m_i$ is the multiplicity of $P_i$ in $Q_d$, this vanishes if and only if $P_i \notin \add Q_d$.
    
    Analogously, the right mutation of $P$ at $P_i$ exists if and only if $X_i \in \Pi_R(\X)[1] = \D^{[-d, -1]}\cap\X$, which is equivalent to $\Hom(\La, X_i) \cong \Hom(Q_0, X_i) = 0$. This occurs if and only if $P_i \notin \add Q_0$. The existence of at least one mutation follows directly from Proposition~\ref{smc-mutation}.
\end{proof}
\sloppy
\printbibliography

@book{KS, 
    series={Grundlehren der Mathematischen Wissenschaften},
    title={Categories and sheaves}, 
    %DOI={10.1007/3-540-27950-4}, 
    publisher={Springer Berlin, Heidelberg}, 
    author={Kashiwara, Masaki and Schapira, Pierre},
    year={2006}
}

@article{R,
author = {Claus Michael Ringel},
title = {Representations of K-species and bimodules},
journal = {Journal of Algebra},
volume = {41},
number = {2},
pages = {269-302},
year = {1976}
}

@article{E,
    author = {Enomoto, Haruhisa},
    title = "{Schur's lemma for exact categories implies abelian}",
    journal = {Journal of Algebra},
    volume = {584},
    pages = {260-269},
    year = {2021}

}

@article{Bu,
    author = {Theo Bühler},
    title = "{Exact categories}",
    journal = {Expositiones Mathematicae},
    volume = {28},
    pages = {1-69},
    year = {2010}, 
    number={1}
}

@article{Bri,
    author = {Tom Bridgeland},
    title = "{Stability conditions on triangulated categories}",
    journal = {Annals of Mathematics},
    volume = {166},
    pages = {317-345},
    year = {2007}
}

@article{LC,
    author = {Jue Le and Xiao-Wu Chen},
    title = "{Karoubianness of a triangulated category}",
    journal = {Journal of Algebra},
    volume = {310},
    number = {1},
    pages = {452-457},
    year = {2007}

}

@article{DF,
    author = {Harm Derksen and Jiarui Fei},
    title = "{General presentations of algebras}",
    journal = {Advances in Mathematics},
    volume = {278},
    pages = {210-237},
    year = {2015}

}

@article{AI,
    author = {Aihara, Takuma and Iyama, Osamu},
    title = "{Silting mutation in triangulated categories}",
    journal = {Journal of the London Mathematical Society},
    volume = {85},
    number = {3},
    pages = {633-668},
    year = {2012}

}

@incollection{BY,
author = {Brüstle, Thomas and Yang, Dong},
title="{Ordered exchange graphs}",
booktitle={Advances in Representation Theory of Algebras},
publisher={European Mathematical Society},
year={2014},
pages={135--193},
series={EMS Series of Congress Reports}
}

@article{KY,
    author    = "Steffen Koenig and Dong Yang",
    title     = "Silting objects, simple-minded collections, $t$-structures and co-$t$-structures for finite-dimensional algebras",
    year      = "2014",
    %addendum = "(accessed: 31.03.2020)",
    journal   = "Doc. Math.",
    volume={19},
    pages ={403-438}
}

@article{AIR,
    title={$\tau $-tilting theory}, 
    volume={150},
    number={3},
    journal={Compositio Mathematica},
    publisher={London Mathematical Society},
    author={Adachi, Takahide and Iyama, Osamu and Reiten, Idun},
    year={2014},
    pages={415–452}
}

@article{AET,
    title={Intervals of $s$-torsion pairs in extriangulated categories with negative first extensions}, 
    volume={174},
    number={3}, 
    journal={Mathematical Proceedings of the Cambridge Philosophical Society}, 
    publisher={Cambridge University Press}, 
    author={Adachi, Takahide and Enomoto, Haruhisa and Tsukamoto, Mayu}, 
    year={2023}, 
    pages={451–469}
}

@article{NP,
  title = {{Extriangulated categories, Hovey twin cotorsion pairs and model structures}},
  author = {Nakaoka, Hiroyuki and Palu, Yann},
  journal = {Cahiers de topologie et g{\'e}om{\'e}trie différentielle cat{\'e}goriques},
  volume = {LX},
  NUMBER = {2},
  PAGES = {117--193},
  YEAR = {2019}
}

@article{F,
  TITLE = {Representations in Abelian Categories},
  AUTHOR = {Freyd, Peter},
  JOURNAL = {Proceedings of the Conference on Categorical Algebra},
  PAGES = {95--120},
  YEAR = {1966}
}

@article{PZ,
author = {Pauksztello, David and Zvonareva, Alexandra},
title = {Co-$t$-structures, cotilting and cotorsion pairs},
volume = {175},
number={1},
journal = {Mathematical Proceedings of the Cambridge Philosophical Society},
year={2023},
pages={89–106}}

@article{IJ,
author = {Iyama, Osamu and Jin, Haibo},
title = {Positive Fuss–Catalan Numbers and Simple-Minded Systems in Negative Calabi–Yau Categories},
volume = {2023},
number={8},
journal = { International Mathematics Research Notices},
year={2023},
pages={6624-6647}}

@article{J,
    title={Auslander-Reiten triangles in subcategories}, 
    author={Peter J{\o}rgensen},
    year={2009},
    journal={Journal of K-Theory},
    volume={3},
    pages={583-601}
}

@article{Ga1,
    title={On thick subcategories of the category of projective presentations}, 
    author={Garcia, Monica},
    journal={Mathematische Zeitschrift},
    year={2026},
    volume={312},
    pages={article number 110}
}

@misc{Gu1,
      title={On $d$-term silting objects, torsion classes, and cotorsion classes}, 
      author={Esha Gupta},
      year={2024},
      eprint={2407.10562v3},
      archivePrefix={arXiv},
      primaryClass={math.RT},
}

@article{MP,
    title={On the Auslander--Reiten theory for extended hearts of proper connective DG algebras}, 
    author={Mochizuki, Nao and Plogmann, Marvin},
    journal = {Bulletin of the London Mathematical Society},
    year = {2026},
    volume={58},
    number={7},
    pages={e70413}
}

@article{AT,
author = {Adachi, Takahide and Tsukamoto, Mayu},
title = {Hereditary cotorsion pairs and silting subcategories in extriangulated categories},
volume = {594},
journal = {Journal of Algebra},
year={2022},
pages={109-137}}

@article{MS,
    author = {Marks, Frederik and Šťovíček, Jan},
    title = {Torsion classes, wide subcategories and localisations},
    journal = {Bulletin of the London Mathematical Society},
    year = {2017},
    volume = {49},
    number={3},
pages={405-416}}

@article{IT, 
title={Noncrossing partitions and representations of quivers}, 
volume={145}, 
number={6}, 
journal={Compositio Mathematica}, 
author={Ingalls, Colin and Thomas, Hugh}, 
year={2009}, 
pages={1533–1562}}

@article{BBDG, 
    title={Faisceaux pervers}, 
    author={Beilinson, Alexander and Bernstein, Joseph and Deligne, Pierre and Gabber, Ofer},
    journal={Astérisque},
    volume={100 Second edition},
    year={2018},
    pages={vi+180 pp}}

@article{P, 
title={Compact corigid objects in triangulated categories and co-$t$-structures}, 
volume={6}, 
number={1},
journal={Central European Journal of Mathematics}, 
author={Pauksztello, David}, 
year={2008}, 
pages={25-42}}

@article{D,
 author = {Dickson, Spencer E.},
 journal = {Transactions of the American Mathematical Society},
 number = {1},
 pages = {223--235},
 title = {A Torsion Theory for Abelian Categories},
 volume = {121},
 year = {1966}
}

@article{XY,
 author = {Xu, Jinde and Yang, Yichao},
 journal = {Archiv der Mathematik},
 number = {4},
 pages = {383--389},
 title = {A Bongartz-type lemma for silting complexes over a hereditary algebra},
 volume = {114},
 year = {2020}
}

@article{IY,
   title={Mutation in triangulated categories and rigid Cohen–Macaulay modules},
   volume={172},
   number={1},
   journal={Inventiones mathematicae},
   author={Iyama, Osamu and Yoshino, Yuji},
   year={2008},
    pages={117--168}
}

@article{A,
    author = {Asai, Sota},
    title = {Semibricks},
    journal = {International Mathematics Research Notices},
    volume = {2020},
    year = {2020},
    number = {16},
    pages = {4993-5054}}

@misc{Z,
      title={Tilting theory for extended module categories}, 
      author={Yu Zhou},
      year={2024},
      eprint={2411.15473},
      archivePrefix={arXiv},
      primaryClass={math.RT}
}

@article{KN,
    author = {Keller, Bernhard and Nicolás, Pedro},
    title = {Weight Structures and Simple dg Modules for Positive dg Algebras},
    journal = {International Mathematics Research Notices},
    volume = {2013},
    number = {5},
    pages = {1028-1078},
    year = {2013}
}

@article{MSSS,
    author = {Mendoza Hernández, Octavio and Sáenz Valadez, Edith Corina and Santiago Vargas, Valente and Souto Salorio, María José},
    title = {Auslander–Buchweitz Context and Co-$t$-structures},
    journal = {Applied Categorical Structures},
    volume = {21},
    number = {5},
    pages = {417--440},
    year = {2013}
}

@article{IJY,
    author = {Osamu Iyama and Peter Jørgensen and Dong Yang},
    title = {Intermediate co-$t$-structures, two-term silting objects, $\tau$-tilting modules, and torsion classes},
    journal = {Algebra and Number Theory},
    volume = {8},
    number = {10},
    pages = {2413--2431},
    year = {2014}
}

@article{AMV,
    author = {Angeleri Hügel, Lidia and Marks, Frederik and Vitória, Jorge},
    title = {Partial silting objects and smashing subcategories},
    journal = {Mathematische Zeitschrift},
    volume = {296},
    pages = {887--900},
    year = {2020}
}

@article{Buan,
    author = {Buan, Aslak Bakke},
    title = {Subcategories of the derived category and cotilting complexes},
    journal = {Colloquium Mathematicum},
    volume = {88},
    pages = {1--11},
    year = {2001}
}
\end{document}